\documentclass[12pt]{amsart}
\usepackage{epsfig}

\usepackage{graphicx}

\newcommand{\R}{{\mathbb R}}

\newcommand{\Q}{{\mathbb Q}}
\newcommand{\Z}{{\mathbb Z}}
\newcommand{\N}{{\mathbb N}}

\newcommand{\lcal}{{\mathcal L}}
\newcommand{\gcal}{{\mathcal G}}
\newcommand{\acal}{{\mathcal A}}

\renewcommand{\phi}{\varphi}

\newcommand{\DB}{\Delta^{\Omega}_B}

\newtheorem{theo}{{\sc Theorem}}[section]
\newtheorem{cor}[theo]{{\sc Corollary}}

\newtheorem{lem}[theo]{{\sc Lemma}}
\newtheorem{sublem}[theo]{{\sc Sublemma}}
\newtheorem{prop}[theo]{{\sc Proposition}}
\newtheorem{defn}[theo]{{\sc Definition}}
\newtheorem{rem}[theo]{{\sc Remark}}

\title[Inverse spectral problem ]{Inverse spectral problem for analytic domains II:
\newline$\Z_2$-
 symmetric domains}

\author{Steve Zelditch}
\address{Department of Mathematics, Johns Hopkins University, Baltimore, MD
21218, USA}
\email{
zelditch@@math.jhu.edu }

\thanks{Research partially supported by  NSF grants \#DMS-0071358 and \#DMS 0302518.}

\date{September 2, 2006}

\begin{document}

\begin{abstract} This paper  develops and implements a new algorithm
for calculating wave trace invariants of a bounded  plane domain
around a periodic billiard orbit. The algorithm is based on a new
expression for the localized wave trace as a special multiple
oscillatory integral over the boundary,  and on a Feynman
diagrammatic analysis of the stationary phase expansion of the
oscillatory integral. The algorithm is particularly effective  for
Euclidean plane  domains possessing  a $\Z_2$ symmetry which
reverses the orientation of a bouncing ball orbit. It is also very
effective for domains with dihedral symmetries. For simply
connected analytic Euclidean plane domains in either symmetry
class, we prove that the domain is determined within the class by
either its Dirichlet or Neumann spectrum. This improves and
generalizes the best prior inverse result (cf. \cite{Z1,Z2, ISZ})
that simply connected analytic plane domains with two symmetries
are spectrally determined within that class.
 \end{abstract}

\maketitle
\section{Introduction}  This paper is part of a series (cf. \cite{Z5, Z4}) devoted
to the inverse spectral problem for simply connected analytic
Euclidean  plane domains $\Omega$. The motivating problem is
whether generic analytic Euclidean drumheads are determined by
their spectra. All known counterexamples to the question, `can you
hear the shape of a drum?', are  plane domains with corners
\cite{GWW1}, so it is possible, according to current knowledge,
that analytic drumheads are spectrally determined. Our main
results give the strongest evidence to date  for this conjecture
by proving it for two classes of analytic drumheads: (i) those
with an up/down symmetry, and (ii) those with a dihedral symmetry.
This  improves and generalize  the best prior results that simply
connected analytic domains with the symmetries of an ellipse and a
 bouncing ball orbit of prescribed length $L$ are spectrally determined within this class
\cite{Z1, Z2, ISZ}.

The proofs of  the  inverse results involve three new ingredients.
The first is a simple and precise expression  (cf. Theorem
\ref{MAINCOR2})   for  the localized trace
 of the wave group (or dually the resolvent), up to a
given order of singularity, as a finite sum of special
 oscillatory integrals over the boundary $\partial
\Omega$ of the domain with transparent dependence on the boundary
defining function.  Theorem \ref{MAINCOR2} is a general result
combining   the Balian-Bloch approach to the wave trace expansion
of \cite{Z5} with a reduction to boundary integral operators
explained in \cite{Z4}. Presumably it could be obtained by other
methods, such as the monodromy operator method  of Iantchenko,
Sj\"ostrand and Zworski \cite{SZ, ISZ}. Aside from this initial
step,
 this  paper is self-contained.

The next and most substantial ingredient is a stationary phase
analysis of the special  oscillatory integrals in Theorem
\ref{MAINCOR2}. To bring order into the profusion of terms in the
wave trace (or resolvent trace) expansion, we use a Feynman
diagrammatic method to enumerate the terms in the expansion.
Diagrammatic analyses have been previously used in \cite{AG} (see
also  \cite{Bu}) to compute the sub-principal wave invariant. A
novel aspect of the diagrammatic analysis in this paper is its
focus on the  diagrams whose amplitudes involve the maximum number
of derivatives of the boundary in a given order of wave invariant.
A key result, Theorem \ref{SETUPA}, is   that only one term, the
{\it principal} term in Theorem \ref{MAINCOR2}, contributes such
highest derivative terms. That is, the stationary phase expansion
of the principal term generates all terms of the $j$th order wave
invariant (for all $j$) which depend on the maximal number $2j -
2$ of derivatives of the curvature of the boundary at the
reflection points.
 In the principal term, the  `transparent dependence' of the phase
  and amplitude on the boundary  is encapsulated
in the  simple properties of the phase and amplitude stated in the
display in Theorem \ref{SETUPA}. Only  these properties are used
to make the key calculations of the wave invariants stated in
Theorem \ref{BGAMMAJ}.

This focus on highest derivative terms in each wave invariant
turns out to be crucial for the inverse spectral problem on
domains with the symmetries studied in this article. The third key
ingredient  is the analysis  in \S \ref{PFONESYM} of these highest
order derivative terms in the case of domains in our two symmetry
classes.  The main result is that the other terms in the wave
invariants are redundant, and further that the domain can be
determined from the wave invariants within these symmetry classes.
These results are
 based on the use of the finite Fourier transform to diagonalize
the Hessian matrix of the length function, and an analysis of
Hessian power sums.

As this outline suggests, we take a direct approach to calculating
wave trace invariants and   do not employ Birkhoff normal forms as
in  \cite{G, Z1, Z2, Z3, ISZ}. We do this because  the classical
normal form of the first return map does not contain sufficient
information to determine domains with only one symmetry. Therefore
one would need to use the full quantum Birkhoff normal form.  But
we found the calculations based on the Balian-Bloch approach
simpler than those involved in the full quantum Birkhoff normal
form.

\subsection{\label{SR} Statement of results}

Let us now state the results more precisely.  We recall that  the
inverse spectral problem for plane domains is to determine a
domain $\Omega$ as much as possible from the spectrum of its
Euclidean Laplacian $\DB$ in $\Omega$ with boundary conditions
$B$:

\begin{equation}\label{EIG}
\left\{ \begin{array}{l} \DB \phi_j(x)  = \lambda_j^2
\phi_j(x),\;\;\;\langle \phi_i, \phi_j \rangle = \delta_{ij},\;\;
(x  \in \Omega) \\ \\
 B \phi_j (q) = 0, \;\;\; q \in
\partial \Omega
\end{array}
\right.
\end{equation}
The  boundary conditions could be either Dirichlet $B \phi = \phi
|_{\partial \Omega}$,  or Neumann $B \phi = \partial_{\nu} \phi
|_{\partial \Omega}$ where $\partial_{\nu}$ is the interior unit
normal.

We briefly introduce some other notation and terminology,
referring to \S 2 and to  \cite{KT} -\cite{PS} for further
background and definitions regarding billiards. By $Lsp(\Omega)$
we denote the length spectrum of $\Omega$, i.e. the set of lengths
of closed trajectories of its billiard flow. By a bouncing ball
orbit $\gamma$ is meant a 2-link periodic trajectory of the
billiard flow.
 The orbit $\gamma$ is a curve in $S^*\Omega$ which
 projects to an `extremal diameter' under the natural
projection $\pi: S^*\Omega \to \Omega,$   i.e. a line segment in
the interior of $\Omega$
 which intersects $\partial \Omega$ orthogonally at both
 boundary points.
 For simplicity of
notation, we often refer to $\pi(\gamma)$ itself as a bouncing
ball orbit and denote it as well by $\gamma$.  By rotating and
translating $\Omega$ we may assume that $\gamma$ is vertical, with
endpoints at  $A = (0,\frac{L}{2})$ and $B = (0, -\frac{L}{2})$.
In a strip $T_{\epsilon}(\overline{AB})$ of width epsilon around
$\gamma$, we may locally express
 $\partial \Omega = \partial \Omega^+ \cup \partial \Omega^-$ as the union of two graphs
over the $x$-axis, namely
\begin{equation}\label{GRAPHS}  \partial \Omega^+ = \{ y = f_+(x),\;\;\; x \in (-\epsilon, \epsilon )\},\;\; \partial \Omega^- = \{ y = f_-(x),\;\;\; x \in
(-\epsilon, \epsilon) \}. \end{equation}

\begin{figure}
\centerline{\includegraphics{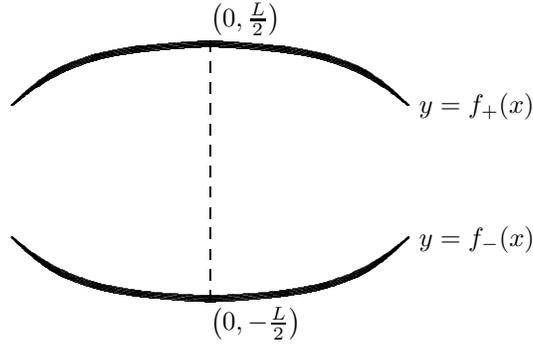}} \caption{$\partial
\Omega$ as a pair of local graphs }
\end{figure}

Our inverse  results  pertain to the following two classes of
drumheads: (i) the class $\mathcal D_{1, L}$ of drumheads with one
symmetry $\sigma$ and a  bouncing ball orbit of length $2L$ which
is reversed by $\sigma$; and (ii) the class $\mathcal D_{m, L}$
$(m \geq 2)$ of drumheads with the dihedral symmetry group $D_m$
and an invariant $m$-link reflecting ray.
 Let us define the classes more precisely and
state the results.

\subsubsection{Domains with one symmetry}

The class
  $\mathcal D_{1, L}$ consists of simply connected real-analytic plane domains $\Omega$ satisfying:
\begin{itemize}

\item  (i) There exists an  isometric involution $ \sigma$ of
$\Omega$ which  `reverses' a  non-degenerate bouncing ball orbit $
\gamma \to \gamma^{-1}$ of length $L_{\gamma} = 2 L$. Hence
$f_+(x) =   - f_-(x)$;

\item (ii) The lengths $2 r L$ of  all iterates $\gamma^r$ ($r =
1, 2, 3, \dots$) have  multiplicity one in $Lsp(\Omega)$, and in
the elliptic case,  the eigenvalues $e^{i \alpha}$  of the linear
Poincare map $P_{\gamma}$ satisfy that $a = - 2 \cos
\frac{\alpha}{2} $ does not belong to the `bad set' ${\mathcal B}
= \{a = 0, -1, 2, -2\}$.

\item (iv) The endpoints of $\gamma$ are not vertices of $\partial
\Omega$.

\end{itemize}

\begin{figure}
\centerline{\includegraphics{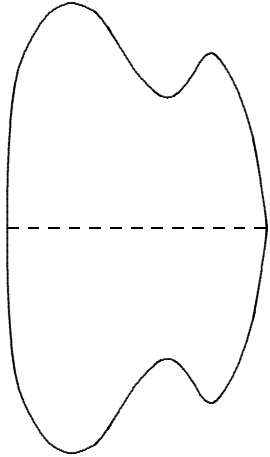}}

 \caption{A domain in
${\mathcal D}_{1, L}$}
\end{figure}

Let Spec$_B(\Omega)$ denote the spectrum of the Laplacian
$\Delta_{\Omega}$  of the domain $\Omega$ with  boundary
conditions $B$ (Dirichlet or Neumann).

\begin{theo} \label{ONESYM} For Dirichlet (or Neumann)  boundary
conditions $B$, the map Spec$_B: {\mathcal D}_{1, L} \mapsto
\R_+^{{\bf N}}$ is 1-1. \end{theo}

Let us clarify the assumptions and consider related problems on
$\Z_2$-symmetric domains:
\medskip

\noindent(a) Under the up-down  symmetry assumption,
 $f_+(x) =   - f_-(x)$ (see Figure (2)). Hence there is 'only one' analytic function $f$ to determine.
  It is quite a different
problem if $\sigma$ preserves orientation of $\gamma$ (i.e. flips
the domain left-right rather than up-down), which amounts to
saying that $f_{\pm}$ are even functions but does not give a
simple relation between them.

\noindent(b) Condition (ii) on the  multiplicity of $2 L$  means
that $\gamma$ is the only closed billiard orbit of  length $2L$.
Since $\gamma = \gamma^{-1}$ for a bouncing ball orbit, the
multiplicity is one rather than two. The method we use to
calculate the trace combines the interior and exterior problems,
and so one might think it necessary to assume that no exterior
closed billiard trajectory (in the complement $\Omega^c$ of
$\Omega$) has length $2L$. However, it is  known that there exists
a purely interior wave trace (cf. \S \ref{OV}) and that the wave
trace invariants at $\gamma$ are spectral invariants; we use the
interior/exterior combination only to simplify the calculation.
Therefore, it is not necessary to exclude exterior closed orbits
of length $L$. When making stationary phase calculations, we only
consider the interior closed orbits.

\noindent(c) The linear Poincar\'e map $P_{\gamma}$ is defined in
\S \ref{BACKGROUND}. In the elliptic case, its eigenvalues
$\{e^{\pm i \alpha}\}$ are of modulus one and we require that $a =
- 2 \cos \frac{\alpha}{2}$ lies outside the bad set ${\mathcal
B}$. In the hyperbolic case, its eigenvalues $\{e^{\pm \alpha}\}$
are real and they are never roots of unity in the non-degenerate
case. These are generic conditions in the class of analytic
domains. We refer to the angles $\alpha$ as Floquet angles. The
set ${\mathcal B}$ consists of  angle parameters where  certain
functions fail to be independent as one `iterates' the geodesic
$\gamma$. The role of this set will be described more precisely in
\S \ref{IRST}.

\noindent(d) Assumption (iii) is equivalent to $f_{\pm}^{(3)}(0)
\not= 0$. The third derivatives $f_{\pm}^{(3)}(0)$ of $f_{\pm}$ at
the endpoints of the bouncing ball orbit appear as coefficients of
certain terms in the  wave invariants, and we make assumption (iv)
to ensure that the corresponding term does not vanish.
Geometrically, $f_{\pm}^{(3)}(0) = 0$  only if the  endpoints of
the bouncing ball orbit are  vertices of $\partial \Omega$, i.e.
critical points of the curvature.  This is a technical condition
which we believe can be  removed by an extension  of the argument,
as will be discussed at the end of the proof. We do not give a
complete argument for the sake of brevity.

As a corollary, we of course have the main result of \cite{Z1, Z2,
ISZ} that a simply connected analytic domain with the symmetries
of an ellipse and with one axis of a prescribed length $L$ is
spectrally determined within this class.

\begin{cor} \label{ONESYMCOR} Let ${\mathcal D}_{2}$ be the class of analytic convex domains with
central symmetry, i.e. the symmetries of an ellipse. Assume that
$\{r L_{\gamma}\}$ are of multiplicity one in $Lsp(\Omega)$ up to
time reversal ($r = 1, 2, 3, \dots$). Then Spec$_B$: ${\mathcal D}
\mapsto \R_+^{{\bf N}}$ is 1-1.
\end{cor}

We give a new proof at the start of \S \ref{PFONESYM} since it is
much simpler than the one-symmetry case and since the proof is
simpler than the ones in \cite{Z1, Z2}.

This inverse result is also true for non-convex simply connected
analytic domains  with the symmetries of the ellipse if we assume
one axis has length $L$ and is of multiplicity one.  We stated the
result only for convex domains because, by a recent result of M.
Ghomi \cite{Gh}, the shortest closed trajectory of a
centrally-symmetric convex domain
 is automatically a bouncing ball orbit, hence it is not necessary to mark the
  length $L$ of an invariant bouncing ball orbit.

Theorem (\ref{ONESYM}) removes  the (left/right) symmetry from the
conditions on the domains considered in \cite{Z1, Z2}.  The
situation for analytic plane domains is now quite analogous to
that for analytic surfaces of revolution \cite{Z3}, where the
rotational symmetry  implies that the profile curve is up/down
symmetric but not necessarily left/right symmetric.

Theorem \ref{ONESYM} admits a  generalization to the special
piecewise analytic mirror symmetric domains with corners which are
formed by reflecting the graph of an  analytic function $y = f(x)$
around the $x$-axis. More precisely, let $ f(x)$ be an analytic
function on an interval  $[-a, a]$ (for some $a$) such that $f(a)
= f(-a) = 0$ and that $f$ has no other zeros in $[-a,a]$. Then
consider the domain $\Omega_f$ bounded by the union of the graphs
$y = \pm f(x)$.

\begin{figure}
\centerline{\includegraphics{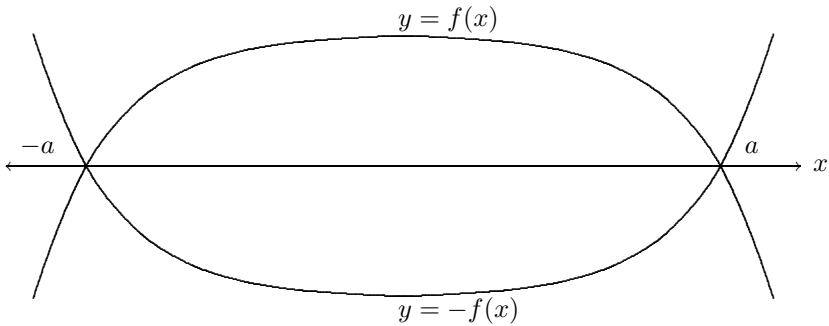}} \caption{$\Z_2$ symmetric
domain with corners}
\end{figure}

 Let ${\mathcal F}$ be the class of real analytic functions with the stated properties,
and consider those $f$ for which  precisely one critical value of
$f$ equals $L/2$. The vertical line through $(x, \pm L/2)$ is then
a bouncing ball orbit. We further  impose the same  generic
conditions on $\Omega_{f}$ as in Theorem \ref{ONESYM}. We denote
the resulting class of real analytic graphs by ${\mathcal F}_L$.

\begin{theo} \label{ONESYMCORN} Up to translation (i.e. choice of $a$), the Dirichlet (or Neumann)
spectrum of $\Omega_f$ determines $f$ within ${\mathcal F}_L$,
i.e.: Spec$: {\mathcal F}_{ L} \mapsto \R_+^{{\bf N}}$ is 1-1.
\end{theo}

The proof is identical to that of Theorem \ref{ONESYM} once it is
established that there exists a wave trace expansion around the
length $t = 2 L$ of the bouncing ball orbit for domains in
${\mathcal F}$ with the same coefficients as in the smooth case.
This fact follows from work of A. Vasy \cite{V} on the Poisson
relation for manifolds with corners. In other words, the presence
of corners does not affect the wave trace expansion at the
bouncing ball orbit.

\subsubsection{Dihedrally symmetric domains}

The second class of domains is the class
$\mathcal D_{m, L}$ of dihedrally symmetric analytic drumheads $\Omega$, i.e. domains
 satisfying:
\begin{itemize}

\item  (i) $\tau \Omega = \Omega $ for all  $\tau \in D_m$;

\item (ii) $D_m$ leaves invariant at least one  $m$-link periodic reflecting ray $  \gamma$  of  length  $ 2L$;

\item (iii)  The lengths $ 2r L$  have  multiplicity one in
$Lsp(\Omega)$

\end{itemize}

\begin{figure}
\centerline{\includegraphics{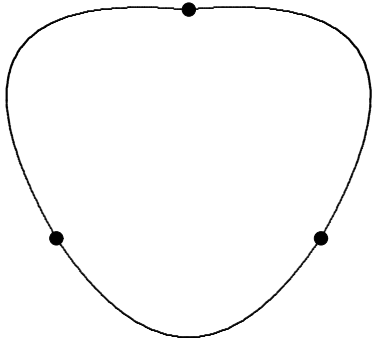}} \caption{A $D_3$-symmetric
domain}
\end{figure}

We then have:

\begin{theo} \label{DISYM} For any $m \geq 2$,
Spec$_B: {\mathcal D}_{m, L} \mapsto \R_+^{{\bf N}}$ is 1-1.
\end{theo}

  We recall  that  $D_m$ is the group generated by
elements $\{\sigma, R_{2 \pi /m}\}$ where $R_{2 \pi /m}$ is
counter-clockwise rotation through the angle $2 \pi /m$ and where
$\sigma^2 = 1$, with the relations $\sigma R_{2 \pi /n} \sigma =
R_{- 2 \pi /n}.$ Also, by an $m$-link periodic reflecting ray we
mean a periodic billiard trajectory with $m$ points of transversal
reflection off $\partial \Omega.$ It is easy to see that such a
ray exists if $\Omega$ is convex. In general, it is a non-trivial
additional assumption. With this proviso, Theorem (\ref{DISYM}) is
a second kind of generalization of the inverse spectral  result of
\cite{Z1, Z2} for the class ${\mathcal D}_{2, L}$ of
`bi-axisymmetric domains'.  That result obviously covers the
classes  ${\mathcal D}_{2n, L}$, but the general case is new. For
any prime $p$, the result for ${\mathcal D}_{p, L}$ is independent
of any other case where $p$ does not divide $n$.

\subsection{\label{OV} Overview}  Let us give a brief overview of
the proofs.

We denote by
$$ E_B^{\Omega} (t,x,y) = \sum_{j} \cos t
\lambda_j \phi_j(x) \phi_j(y)$$  the kernel of the even part of
the
 wave group $\cos t \sqrt{\Delta_B^{\Omega}}$,  generated by the Laplacian  $\Delta_B^{\Omega}$
 of (\ref{EIG}) with either Dirichlet $B u =
u|_{\partial \Omega}$ or Neumann $Bu =
\partial_{\nu} u |_{\partial
\Omega}$ boundary conditions.    Its distribution trace is defined
by
\begin{equation} \label{WAVETRACE} Tr 1_{\Omega}  E_{B}^{\Omega} (t)
:= \int_{\Omega} E_{B}^{\Omega}(t, x,x) dx = \sum_{j = 1}^{\infty}
\cos t \lambda_j \end{equation}

When $L_{\gamma}$ is the length of a non-degenerate periodic
reflecting ray $\gamma$ of the generalized  billiard flow, and
when the only periodic orbits of length $L_{\gamma}$ are $\gamma$
and $\gamma^{-1}$ (the time-reversal of $\gamma$),  then $Tr
1_{\Omega} E_B^{\Omega}(t)$ is a Lagrangian distribution in the
interval $(L_{\gamma} - \epsilon, L_{\gamma} + \epsilon)$ for
sufficiently small $\epsilon$, and has the following expansion in
terms of homogeneous singularities: (see \cite{GM}, Theorem 1, and
also page 228; see also \cite{PS} Theorem 6.3.1).
\medskip

{\it Let $\gamma$ be a non-degenerate billiard trajectory whose
length $L_{\gamma}$ is isolated and of multiplicity one in $Lsp
(\Omega)$.  Then for $t$ near $L_{\gamma}$, the trace of the even
part of the  wave group has the singularity expansion
\begin{equation} \label{Wavetrace} Tr 1_{\Omega}  E_B^{\Omega} (t) \sim \Re\{ a_{\gamma}
(t - L_{\gamma} + i0)^{-1} + a_{\gamma 0} \log (t - L_{\gamma} + i
0) + \sum_{k = 1}^{\infty}a_{\gamma k}(t - L_{\gamma} + i 0)^k
\log (t - L_{\gamma} + i 0)\}, \end{equation} where the
coefficients $a_{\gamma k}$ (the wave trace invariants) are
calculated by the stationary phase method from a microlocal
parametrix for $E_B^{\Omega}$ at $\gamma$. }
\medskip

Here, $a_{\gamma}$ is a sum of the contributions from $\gamma$ and
$\gamma^{-1}$, which are the same.
 In general, the contribution at $t = L_{\gamma}$ is the sum
over all periodic orbits of length $L_{\gamma}$. The sum to the
right of $\Re$ is the trace of  the wave group $e^{i t
\sqrt{\Delta_B^{\Omega}}}$; the trace of the even part $
E_B^{\Omega} (t)$  of the wave group   equals the real part of
that trace.

 In \cite{Z5}, \S 3.1,  this expansion was reformulated
  in terms of a regularized trace of the
interior resolvent $R_B^{\Omega}(k + i \tau) = - (\DB +( k + i
\tau)^2)^{-1} : H^{s}(\Omega) \to H^{s + 2}(\Omega)$, with $k \in
\R, \tau > 0$ and with  boundary condition $B$. The Schwartz
kernel or  {\it Green's kernel} $G_B^{\Omega}(k + i \tau, x, y)
\in \mathcal D'(\Omega \times \Omega)$ of the resolvent  is the
unique solution of the boundary problem:
\begin{equation}\label{GREENK}
\left\{ \begin{array}{l}  - (\DB + (k + i \tau)^2)
G_{B}^{\Omega}(k + i \tau, x, y) = \delta_{y}(x),\;\;\;
(x, y \in \Omega) \\ \\
 B G_{B}^{\Omega}(k + i \tau, x, y) = 0, \;\;\; x \in
\partial \Omega \end{array}
\right.
\end{equation}

  Let  $\hat{\rho} \in C_0^{\infty}(L_{\gamma} - \epsilon, L_{\gamma} + \epsilon)$ be a cutoff,  equal to one on an interval  $(L_{\gamma} - \epsilon/2, L_{\gamma} + \epsilon/2)$ which contains
no other lengths in Lsp$(\Omega)$ occur in its support, and define
the smoothed (and localized) resolvent with a choice of boundary
conditions  by
\begin{equation}\label{RRHO} R_{B \rho}^{\Omega}(k + i \tau):= \int_{\R}
 \rho (k - \mu) (\mu + i \tau)  R_B^{\Omega}(\mu + i \tau) d \mu.  \end{equation}
 The definition is chosen so that
 \begin{equation}\label{RTWT} R_{B \rho}^{\Omega}(k + i \tau) =
 \int_0^{\infty} \hat{\rho}(t) e^{i(k + i \tau)t} E_B^{\Omega}(t)
 dt. \end{equation}
Then the smoothed resolvent trace admits an asymptotic expansion
of the form
\begin{equation} \label{PR} Tr 1_{\Omega} R_{B \rho}^{\Omega}
(k + i \tau) \sim  {\mathcal D}_{B, \gamma}(k + i \tau)  \sum_{j =
0}^{\infty} B_{\gamma, j}  k^{-j},\;\;\; k \to
\infty,\end{equation} where \begin{itemize}

\item ${\mathcal D}_{B, \gamma}(k + i \tau)$ is the {\it
symplectic pre-factor} $${\mathcal D}_{B, \gamma}(k + i \tau)  =
C_0 \; \epsilon_B(\gamma) \frac{e^{ i (k + i \tau) L_{\gamma}}
e^{i \frac{\pi}{4} m_{\gamma}}}{\sqrt{|\det (I - P_{\gamma})|}}$$

\item $P_{\gamma}$ is the Poincar\'e map associated to $\gamma$
(see \S \ref{BACKGROUND} for background);

\item $\epsilon_B(\gamma)$ is the signed  number of intersections
of $\gamma$ with  $\partial \Omega$ (the sign depends on the
boundary conditions; $\pm 1$ for each bounce for Neumann/Dirichlet
boundary conditions);

\item  $m_{\gamma}$ is the Maslov  index of $\gamma$;

\item $C_0$ is a universal constant (e.g. factors of $2 \pi$)
which it is not necessary to know for the proof of Theorem
\ref{ONESYM}.

\end{itemize}
The resolvent trace (or Balian-Bloch)  coefficient $B_{\gamma, j}$
associated to a periodic orbits $\gamma, \gamma^{-1}$ is easily
related to the wave trace coefficient $a_{\gamma, k}$.
   We
henceforth work solely with the expansion (\ref{PR}), which we
term the `Balian-Bloch expansion' after \cite{BB2}. In fact, we
actually analyze the closely related resolvent trace asymptotics
along logarithmic curves $k + i \tau \log k$ in the upper half
plane. It is clear that the `Balian-Bloch coefficients'
$B_{\gamma, j}$ are spectral invariants and it is these invariants
we use in our inverse spectral results.

As mentioned above,  the inverse results have  three  main
ingredients, which we now describe in  detail  as a guide to the
paper and its connections to \cite{Z4, Z5}.

\subsubsection{Reduction to boundary oscillatory integrals of the
wave trace}

 The first step  (Theorem \ref{MAINCOR2}) is a reduction to the boundary of
the wave trace. This reduction was largely achieved in \cite{Z5,
Z4} by means of a rigorous version of the Balian-Bloch approach to
the Poisson relation between spectrum and closed billiard orbits
\cite{BB1, BB2}.   It expresses the wave trace localized at the
length of  a periodic reflecting ray,  up to a given order of
singularity, as a finite sum of oscillatory integrals  $ I_{M,
\rho}^{\sigma, w}(k + i \tau) $ over the boundary (see
(\ref{EXPRESSION}). It is related in spirit to the monodromy
operator approach of \cite{SZ, ISZ}.

\subsubsection{Feynman diagram analysis and proof of Theorem
\ref{SETUPA}}  The second ingredient is a  stationary phase
analysis of the oscillatory integral expressions for the wave
invariants at transversally reflecting periodic orbits. The key
role is played by  a (Feynman) diagrammatic analysis of the
stationary phase expansions, which has not previously been used in
inverse spectral theory (see \cite{AG} for prior use in calculated
the sub-principal invariant). As reviewed in \S \ref{SPFD}, the
terms of stationary phase expansion correspond to labelled graphs
$\Gamma$ and the coefficients of the stationary phase expansion
can be expressed as `Feynman amplitudes' determined by the graphs
$\Gamma$. The Euler characteristic of $\Gamma$ corresponds to the
power $k^{-j}$ of $k$ in the wave trace expansion.

The inverse spectral problem involves a novel point of the
diagrammatic analysis: namely, to separate out the (labelled
graphs) of Euler characteristic $-j$ whose amplitudes contain the
maximum numbers ($2j + 2, 2j - 1$) of derivatives of $\partial
\Omega.$ In Theorem \ref{SETUPA} we prove that  the terms in a
given wave invariant which contain the maximal number of
derivatives of $\partial \Omega$ only arise in the stationary
phase expansion of
 one {\it principal term} and its time reversal,  whose amplitudes have
special properties stated in table in Theorem \ref{SETUPA}. The
principal terms are defined
 in Definition \ref{PRINCIPAL}. Only the special properties of the phase and
amplitude are used in the calculation of the wave trace
invariants.

The analysis leads to the explicit formulae for the top derivative
parts of the  wave invariants at iterates of bouncing ball orbits
in Theorem \ref{BGAMMAJ}. For instance, in the symmetric bouncing
ball case there is only one important diagram for the even
derivatives $f^{(2j)}(0)$ and two important diagrams for the odd
derivatives $f^{(2j - 1)}(0)$. Modulo terms involving $\leq 2j -
2$ derivatives, the wave trace (or more precisely resolvent trace)
invariants (cf. (\ref{Wavetrace})-(\ref{PR}) $B_{\gamma^r, j - 1}$
 take the form (cf. Corollary \ref{BGAMMAJSYM}):
\begin{equation}  \label{BGAMMAJSYMpre}  \begin{array}{lll} B_{\gamma^r, j -
 1} & = &    (4 L r) {\mathcal A}_r(0) i^{j - 1} \{ 2 (w_{{\mathcal G}_{1, j}^{2j, 0}})\;
 (h_{2r}^{11})^j f^{(2j)}(0) \\ && \\ && + 4 (w_{{\mathcal G}_{2, j + 1 }^{2j - 1, 3,
 0}}) \;
 (h_{2r}^{11})^j \frac{1}{2 - 2 \cos \alpha/2} ( f^{(3)}(0) f^{(2j - 1)}(0)) \\ && \\
 && +  4 (w_{\widehat{{\mathcal G}}_{2, j + 1 }^{2j - 1, 3, 0}})\; (h_{2r}^{11})^{j - 2}
 \sum_{q = 1}^{2r} (h_{2r}^{1 q})^3 ( f^{(3)}(0) f^{(2j - 1)}(0))\}.\end{array}  \end{equation}

 Here and throughout the paper we use the following notational
 conventions:
 \begin{itemize}

\item $h_{2r}^{pq}$ are the matrix elements of the inverse of the
Hessian $H_{2r}$ of the length function ${\mathcal L}$ in
Cartesian graph coordinates at $\gamma^r$ (cf. \S
\ref{BACKGROUND}).

\item ${\mathcal A}_r(0)$ is an $\Omega$-independent  (non-zero)
 constant obtained from amplitude of  the principal terms at
the critical bouncing ball orbit.

\item $w_{{\mathcal G}_{1, j}^{2j, 0}}$ (etc.) are certain
non-zero combinatorial constants  associated to  Feynman graphs
denoted here by ${\mathcal G}_{1, j}^{2j, 0}$ etc.  For a given
graph $\gcal$,  $w_{\gcal} = \frac{1}{|Aut(\gcal)|}$ where
$|Aut(\gcal)|$ is the order of the symmetry group of the graph;
see the discussion after (\ref{MSPJTERM}).

\end{itemize}

The amplitude value ${\mathcal A}_r(0)$ and the Wick constants may
be evaluated explicitly. However it is not necessary for the proof
of Theorem \ref{ONESYM} to do so and it seems more illuminating to
specify the origins, rather than their values, of the various
constants. We note that the $h_{2r}^{ij}$ depend on, and only on,
$r$ and  the eigenvalues of the Poincar\'e map $P_{\gamma}$ (i.e.
on the Floquet angles) and on the length of $\gamma$. We also note
that $\gamma = \gamma^{-1}$ when $\gamma$ is a bouncing-ball orbit
(such an orbit is called reciprocal).

 The analysis
shows  that the non-principal oscillatory integrals only give rise
to sub-maximal  derivative terms in the wave invariants,
completing the proof of Theorem \ref{SETUPA}.

\subsubsection{\label{IRST} Inverse results}

The third ingredient is the analysis of the top derivative terms
in the wave trace invariants  in the symmetry classes above. The
key point is determine the $2j-1$st and $2j$th Taylor coefficients
of the curvature at each reflection point from the $j-1$st wave
trace invariant for $\gamma$ and its iterates $\gamma^r.$

We note that the previously known inverse result for
 analytic domains with the symmetry of an ellipse drops out
 immediately from (\ref{BGAMMAJSYMpre}), since the odd Taylor
 coefficients are zero. On the other hand, there is an obstruction
 to recovering the Taylor coefficients of $f$ when there is only
 one symmetry: namely, we must recover two Taylor coefficients
 $f^{(2j)}(0), f^{(2j - 1)}(0)$ for each new value of $j$ (the
 degree of the singularity). This is the principal obstacle to
 overcome.

We overcome it in \S \ref{PFONESYM} as follows: The expression
(\ref{BGAMMAJSYMpre}) for the Balian-Bloch invariants of $\gamma,
\gamma^2, \dots$ consists of two types of  terms, in terms of
their dependence on the iterate $r$. They have a common factor of
$2 r L (h_{2r}^{11})^{j - 2} {\mathcal A}_r(0) $, and after
factoring it out we obtain one term
$$ (h_{2r}^{11})^2 \{
(w_{{\mathcal G}_{1, j}^{2j, 0}}) f^{(2j)}(0) + \frac
{(w_{{\mathcal G}_{2, j + 1 }^{2j - 1, 3,
 0}})}{2 - 2 \cos
\alpha/2} f^{(3)}(0) f^{(2j - 1)}(0)\}  $$
 which  depends on the iterate $r$ through the coefficient
  $ (h_{2r}^{11})^2$, and one $$ (w_{\widehat{{\mathcal G}}_{2, j + 1 }^{2j - 1, 3, 0}}) \left(
\sum_{q = 1}^{2r} (h_{2r}^{1 q})^3 f^{(3)}(0) f^{(2j - 1)}(0)
\right)$$
 which depends on $r$ through the cubic sums $ \sum_{q = 1}^{2r}
(h_{2r}^{1 q})^3$ of inverse Hessian matrix elements
$h_{2r}^{pq}$. In order to `decouple' the even and odd
derivatives, it suffices to show that the functions
$(h_{2r}^{11})^2 $ and $\sum_{q = 1}^{2r} (h_{2r}^{1 q})^3$ are,
at least for `most' Floquet angles $\alpha$,  linearly independent
as functions of $ r\in \Z$, i.e. that $(h_{2r}^{11})^{-2} \sum_{q
= 1}^{2r} (h_{2r}^{1 q})^3$ is a non-constant function of $r$. It
is convenient to use the parameter $a = - 2 \cos \frac{\alpha}{2}$
and we write the dependence as $h_{2r}^{ij}(a)$.

 We therefore define the `bad' set of Floquet angles by
\begin{equation} \label{BAD} {\mathcal B} = \{a:\; \;\mbox{the sequence }\; \{(h_{2r}^{11}(a))^{-2}
\sum_{q = 1}^{2r} (h_{2r}^{1 q}(a))^3, \;\; r = 1, 2, 3, \dots \}
\; \mbox{is constant in }\; r\}. \end{equation} Using facts about
the finite Fourier transform and circulant matrices, we compute
that $ {\mathcal B} = \{0, 1, \pm 2\}.$ Since the
 proof is  computational,
  we also present a simple conceputal argument (cf. Proposition
  \ref{FINITE})that ${\mathcal B}$ is finite, although the proof
  only gives the poor
 estimate $3^{20}$ on its number of elements.  For Floquet angles outside of ${\mathcal B}$,
  we can determine
 all Taylor coefficients $f_+^{(j)}(0)$ from the wave
 invariants and hence the analytic domain.

We use a similar strategy in the dihedral $D_n$-case in \S
\ref{PFDISYM}.  Due to the extra symmetries, the inverse results
in the dihedral  case require much less information about the wave
invariants than in the one symmetry case.

\subsection{Related results} \medskip

\noindent{\bf (i)}  We have already mentioned the prior result
that analytic drumheads with up/down and left/right symmetries are
spectrally determined in that class \cite{Z1, Z2}. Previously, it
was proved by  Colin de Verdiere \cite{CV} that such domains are
spectrally rigid.  To our knowledge, the only other prior result
giving a `large' class of spectrally domains is that of
Marvizi-Melrose \cite{MM1}, in which members of  a spectrally
determined two-parameter family of convex plane domains are
determined among generic convex domains by their spectra. \medskip

\noindent{\bf (ii)} In  \cite{Z4}, we extend the inverse result to
the exterior problem of determining a $\Z_2$-symmetric
configuration of analytic obstacles from its scattering phase (or
resonance poles). Our result may be stated as follows: Let $\Omega
= \R^2 - \{{\mathcal O} \cup \tau_{x, L} {\mathcal O}\}$ where
${\mathcal O}$ is a convex analytic obstacle, where $x \in
{\mathcal O}$ and where $\tau_{x,L}$ is the mirror reflection
across the orthogonal line segment of length $L$ from $x$. Thus,
$\{ {\mathcal O} \cup \tau_{x, L} ({\mathcal O})\}$ is a
$\Z_2$-symmetric obstacle consisting of two components. Let
$\Delta_{\Omega}$ denote the Dirichlet Laplacian on $\Omega.$ We
have:

\begin{theo} \cite{Z4} With the same genericity assumptions
as in Theorem \ref{ONESYM},  the resonance poles of
$\Delta_{\Omega}$ determine ${\mathcal O}$ within the class of
$\Z_2$ symmetric analytic obstacles.
\end{theo}

\subsection{Future directions}

 An obvious future  direction is to  study the wave
invariants without any symmetry assumptions. As will become clear
from the calculations in this article (cf. Theorems \ref{SETUPA}
and \ref{MAINCOR2}), symmetries make `lower order derivative data'
in wave invariants  redundant and allow one to concentrate on
terms in a given wave invariant with maximal numbers of
derivatives. Lacking symmetries, the lower order derivative data
is no longer redundant and one has to navigate a complicated
jungle of terms to determine which combinations are spectral
invariants. It is plausible that one cannot work with just one
orbit but must combine information from two bouncing ball orbits
(they always exist in a convex plane domain). The main problem is
then to extract from the wave invariants of the iterates of each
bouncing ball orbit sufficient Taylor series data at the endpoints
to determine the domain. To do this, it seems necessary to analyze
 how  Feynman amplitudes of labelled diagrams behave
as a function of the iterate $r$ of the orbits. The graphs
themselves do not depend on $r$, so the dependence comes from the
labelling.

\subsection{Acknowledgements}  The first draft of this article was
posted in  2001 (arXiv math.SP/0111078) as part of the series now
published as \cite{Z4, Z5}. In the intervening period,  some of
the computational details of this article have received
independent confirmation.
 R. Bacher found an independent proof corroborating the result of
\S \ref{FD} that only five graphs in the Feynman diagrammatics
have the right form to contribute to the highest order data
\cite{B}. C. Hillar did a numerical study of the circulant sums of
\S \ref{PFONESYM} to corroborate that the nonlinear power sums
have  nonlinear $r$ dependence as  $r \to \infty$. A. Vasy and J.
Wunsch gave advice on  \cite{V}.  We thank Y. Colin de Verdi\`ere
and the spectral theory seminar at Grenoble for their patience in
listening to earlier versions of this article and for their
comments. Especially,  we thank the referee for many remarks and
corrections. The referee pointed out that   `bad set' ${\mathcal
B}$ could be explicitly calculated, and the calculations were done
by H. Hezari. We thank H. Hezari as well  for carefully reading
the final version.

\tableofcontents

\section{\label{BACKGROUND} Billiards and the length functional}

We begin by  establishing   notation
 on plane billiards and length functions.  After recalling basic
 notions, we calculate the Hessian of the length functional at
 iterates of a
 critical bouncing ball orbit
in  Cartesian coordinates adpated to the orbit.

We denote by  $\Omega$ a simply connected analytic plane domain
with  boundary $\partial \Omega$ of length $2\pi$. The billiard
flow $\Phi^t$ of $\Omega$ is the broken geodesic of the Euclidean
metric on $\Omega$. That is, for $(x, \xi) \in T^* \Omega^o$, the
trajectory $\Phi^t(x, \xi)$ follows the Euclidean straight line in
the interior $\Omega^o$ of $\Omega$ and reflects from the boundary
by Snell's law of equal angles. By the billiard map $\beta$ of
$\Omega$ we mean the  map on $B^*
\partial \Omega$  induced
by $\Phi^t$: we add a multiple of the inward unit normal $\nu_q$
to $(q, \eta) \in B^*(\partial \Omega)$ to obtain an inward
pointing unit vector $v$ at $q$. We then follow the billiard
trajectory $\Phi^t(q, v)$ until it hits the boundary, and then
define $\beta(q, \eta)$ to be its tangential projection. We refer
 to \cite{PS, KT, Z5} for details and discussions of
 the billiard flow on domains in $\R^2$.

 It is natural at first to  parametrize
$\partial \Omega$ by arclength,
\begin{equation} q: {\bf T} \to \partial \Omega \subset \R^2, \end{equation}
 starting at some
point $q_0 \in \partial \Omega$. Here, ${\bf T} = \R \backslash 2
\pi \Z$ denotes the unit circle. By an $m$-link  {\it periodic
reflecting ray}  of $\Omega$ we mean a periodic billiard
trajectory $\gamma$ which intersects $\partial \Omega$
transversally at $m$ points $q(\phi_1), \dots, q(\phi_m)$, and
reflects off $\partial \Omega$ at each point according to Snell's
law
\begin{equation}\label{SN} \frac{q(\phi_{j + 1}) - q(\phi_j)}{|q(\phi_{j + 1}) - q(\phi_j)|} \cdot \nu_{q(\phi_j)} =
 \frac{q(\phi_{j }) - q(\phi_{j - 1})}{|q(\phi_{j }) - q(\phi_{j - 1})|} \cdot \nu_{q(\phi_j)}. \end{equation}
Here,  $\nu_{q(\phi)}$ is the inward  unit normal to $\partial \Omega$ at $q(\phi)$. We refer to the segments
 $q(\phi_{j }) - q(\phi_{j - 1})$ as the {\it links} of the trajectory.
We  denote the acute angle between the link $q(\phi_{j + 1 }) -
q(\phi_{j })$ and the inward unit normal $\nu_{q(\phi_{j})}$ by
$\angle ( q(\phi_{j + 1}) - q(\phi_{j }), \nu_{q(\phi_j)})$
 and that between  $q(\phi_{j }) - q(\phi_{j - 1 })$ and
the  inward unit  normal  at $q(\phi_{j })$ by $\angle (q(\phi_{j
 }) - q(\phi_{j - 1 }), \nu_{q(\phi_{j })})$, i.e. we put
\begin{equation}  \frac{q(\phi_{j + 1}) - q(\phi_j)}{|q(\phi_{j + 1}) - q(\phi_j)|} \cdot \nu_{q(\phi_j)} =
\cos \angle (q(\phi_{j + 1 }) - q(\phi_{j }), \nu_{q(\phi_{j })}).
\end{equation} For notational simplicity we often do not
distinguish between a billiard trajectory in $S^*\Omega$ and its
projection to $\Omega$.

We  define the length functional on ${\bf T}^M$ by:
\begin{equation}\label{LENGTH}    L(\phi_1, \dots, \phi_M) =  |q(\phi_1) - q(\phi_2)| + \dots
 + |q(\phi_{M - 1} ) -q(\phi_M)| + |q(\phi_M) - q(\phi_1)|.\end{equation}
 We often use cyclic index notation where $q(\phi_{M + 1}) =
 q(\phi_1). $
 It is clear  that $L$ is a smooth function away from
the `large diagonals' $\Delta_{j, j + 1}:= \{\phi_j = \phi_{j +
1}\}$, where it   has  $|x|$   singularities.  We have:
\begin{equation}\label{Oneder} \begin{array}{l} \frac{\partial}{\partial \phi_j} |q(\phi_j) - q(\phi_{j - 1})| =
 - \sin \angle ( q(\phi_{j }) - q(\phi_{j -1 }), \nu_{q(\phi_j)}),\\ \\
\frac{\partial}{\partial \phi_{j }} |q(\phi_j) - q(\phi_{j + 1})|
 =  \sin \angle ( q(\phi_{j + 1}) - q(\phi_{j }), \nu_{q(\phi_{j})}) \\ \\ \implies
 \frac{\partial}{\partial \phi_j} L  = \sin \angle ( q(\phi_{j + 1}) - q(\phi_{j }), \nu_{q(\phi_{j})}) - \sin \angle ( q(\phi_{j }) - q(\phi_{j -1 }), \nu_{q(\phi_j)}) .\end{array}\end{equation}
  Hence, the  condition that $\frac{\partial}{\partial \phi_j} L
= 0$ is the same as (\ref{SN}) for  the $2$-link defined by the
triplet $(q(\phi_{ j- 1}), q(\phi_j), q(\phi_{j + 1}))$.

Let $\gamma$ denote a periodic reflecting ray of $\Omega$.
 The linear
Poincare map $P_{\gamma}$ of $\gamma$ is the derivative at
$\gamma(0)$ of the first return map to a transversal to $\Phi^t$
at $\gamma(0).$ By a non-degenerate periodic reflecting ray
$\gamma$ we mean one whose linear Poincar\'e map $P_{\gamma}$ has
no eigenvalue equal to one (cf.  \cite{PS, KT}). The following
relates $P_{\gamma}$ and the Hessian of the length functional in
angular coordinates:

\begin{prop} \label{KTa} (\cite{KT} (Theorem 3)) Let $H_n^a$ denote the Hessian of $L$  in angular coordinates
$\phi_j$ at a critical point $\gamma$, and let $b_j =
\frac{\partial^2 | q(\phi_{j + 1}) - q(\phi_j)|}{\partial \phi_j
\partial \phi_{j + 1}}.$ Then
$$\det (I - P_{\gamma}) = - \det (- H_n^a) \cdot (b_1 \cdots
b_n)^{-1}.$$
\end{prop}
This identity may be proved by expressing both sides in terms of
bases of horizontal and vertical Jacobi fields.

\subsection{\label{BBO} Cartesian coordinates around bouncing ball orbits}

We now specialize to the case where $\gamma$ is a  bouncing ball
orbits (i.e. $2$-link periodic reflecting rays). As in the
Introduction, we orient $\Omega$ so that the bouncing ball orbit
is along the $y$-axis with endpoints $A = (0,\frac{L}{2}), B = (0,
-\frac{L}{2})$ and parametrize $\partial \Omega$ near $A$ by $y =
f_+(x)$ and near $B$ by $y = f_-(x)$. We do not assume the domain
is up-down symmetric.

We denote by $R_A$, resp. $R_B$, the radius of curvature of
$\Omega$ at the endpoints $A, B$.  When $\gamma$ is elliptic, the
eigenvalues of $P_{\gamma}$ are of the form $\{e^{\pm i \alpha}\}$
($\alpha \in \R$) while in the hyperbolic case they are of the
form $\{e^{\pm \alpha}\}$ ($\alpha \in \R)$. They are given by the
same formulae in both elliptic and hyperbolic cases:
\begin{equation} \label{COSALPHA} \left\{ \begin{array}{ll}
 \cos (\alpha/2) = \sqrt{(1 - \frac{L}{R_A}) (1 - \frac{L}{R_B})}, & \;\;(\mbox{elliptic
case}),\\ \\ \cosh (\alpha/2) = \sqrt{(1 - \frac{L}{R_A}) (1 -
\frac{L}{R_B})}, & \;\;(\mbox{hyperbolic case}). \end{array}
\right.
\end{equation}

We  define the  length functionals in Cartesian coordinates for
the two possible orientations of the  $r$th iterate of a bouncing
ball orbit by
\begin{equation} \label{LPM} {\mathcal L}_{\pm } (x_1, \dots, x_{2r}) = \sum_{j = 1}^{2r} \sqrt{(x_{j + 1} - x_j)^2 +
 (f_{w_{\pm}(j + 1) }(x_{j +1}) - f_{w_{\pm}(j)}(x_j))^2}. \end{equation}
 Here, $w_{\pm}: \Z_{2r} \to \{\pm\}$, where $w_{+}(j)$ (resp. $w_-(j))$ alternates sign starting with
 $w_+(1) = +$ (resp. $w_-(1) = -$). Also, we use cylic index
 notation where $x_{2r + 1} = x_1$.

We have:
\begin{equation} \label{crit} \begin{array}{lll} \frac{\partial {\mathcal L}_{\pm}}{\partial
x_j} &= &
 \frac{(x_j - x_{j + 1}) + (f_{w_{\pm}(j)}(x_j) - f_{w_{\pm}(j +
 1)}
 (x_{j + 1})) f_{w_{\pm}(j)}'(x_j)}{\sqrt{(x_j - x_{j + 1})^2 +
(f_{w_{\pm}(j)}(x_j) - f_{w_{\pm}(j + 1)}(x_{j + 1}) )^2}} \\ && \\
&& - \frac{(x_{j - 1} - x_{j }) + (f_{w_{\pm} (j - 1)} (x_{j - 1})
- f_{w_{\pm}(j)}(x_{j}) ) f_{w_{\pm} (j)}'(x_j)}{\sqrt{(x_j - x_{j
- 1})^2 + (f_{w_{\pm}(j)}(x_j) - f_{w_{\pm}(j - 1)}(x_{j - 1})
)^2}}.
\end{array} \end{equation}

We will  need formulae for the entries of the  Hessian of
${\mathcal L}_{+}$ at its critical point $(x_1, \dots, x_{2r}) =
0$ in Cartesian coordinates corresponding to the $r$th repetition
of a bouncing ball orbit.

\begin{prop} \label{KT} Put  $$a = - 2(1 + L f_+''(0)) = - 2 (1 - \frac{L}{R_A}),\; b = - 2(1 - L
f_-''(0)) = -2(1 -  \frac{L}{R_B}). $$  Then the Hessian $H_{2 r}$
of ${\mathcal L}_+$ at $x = 0$ in Cartesian graph coordinates has
the form $H_2 = \frac{- 1}{L}
\begin{pmatrix} a & 2 \\ & \\ 2 & b \end{pmatrix}$ for $r = 1$ and
for $r \geq 2$,
$$ \label{HBB} H_{2r} = \frac{-1}{L} \left\{ \begin{array}{lllll} a & 1 & 0 & \dots & 1  \\ & & & & \\
  1 & b & 1 & \dots & 0 \\ & & & & \\ 0 & 1 &  a & 1 & 0  \\ & & & & \\ 0 & 0 & 1 & b & 1 \dots
 \\ & & & & \\\dots & \dots & \dots & \dots & \dots \\ & & & & \\
1 & 0 & 0 & \dots & b \end{array} \right\}$$. \end{prop}

\begin{proof} A routine calculation gives
$$ \left\{ \begin{array}{l}  \frac{\partial^2 {\mathcal L}_+}{\partial x_j^2}(0) =
 2 ( \frac{1}{L} + w_+(j) f_{w_+(j)}''(0))
 ,\;\;\;\; \\ \\
\frac{\partial^2 {\mathcal L}_+}{\partial x_{j} \partial x_{j+1}}
(0) =
  \frac{- 1}{L} \end{array}\right.$$
  for $r \geq 2$. In the case of $r = 1$, the length functional is
  $2 ||(x_1, f_+(x_1)) - (x_2, f_-(x_2))||$.
 Note that for $r \geq 2$, there  are two terms of ${\mathcal L}_+$ contributing
to each diagonal matrix element and one to each off-diagonal
element, accounting for the additional factor of $2$ in the
diagonal terms. Also note that $f_{w_{+}(j)}(0) - f_{w_{+}(j +
 1)}
 (0)  = w_+(j) L $ and that $f_+''(0) = \frac{-1}{R_A}, f_-''(0) =
 \frac{1}{R_B}. $

\end{proof}

We remark that the Hessian in Cartesian coordinates in Proposition
\ref{KT} differs from that in angular coordinates in \cite{KT} in
that the off-diagonal entries differ in sign. This is because the
graph parametrization gives the opposite orientation to the
tangent $T_A \partial \Omega$ than the angular parametrization and
the same orientation at $T_B \partial \Omega$. The angular Hessian
$H^a_{2r}$ is related to the Cartesian Hessian $H_{2r}$   by
$H_{2r}^a = J H_{2r} J^t$ where $J = diag(1, -1, 1, -1, \dots, 1,
-1)$ is the change of basis matrix. Clearly, the determinants of
the two Hessians agree. Since $b_j = \frac{-1}{L}$, we obtain from
Proposition \ref{KTa} the following:

\begin{cor} \label{DETP} As above, let $H_{2r}$ denote the Hessian of ${\mathcal L}_+$ in Cartesian coordinates
at the $r$th iterate $\gamma^r$  of a  bouncing ball orbit
$\gamma$ of length $2L$. Then
$$\det (I - P_{\gamma^r}) = - L^{2r} \det ( H_{2r}).$$
\end{cor}

 The determinant   $\det H_{2r}$ is a polynomial in $\cos  \frac{\alpha}{2}$ (elliptic case),
 resp. $\cosh \frac{\alpha}{2}$ (hyperbolic case) of degree $2r$.
 In the following we restrict to the elliptic case.
\begin{prop}  \label{DetHESS}We have
$$\det H_{2r} =  - L^{- 2r} (2 - 2  \cos r \alpha) . $$

\end{prop}

\begin{proof} Let $\lambda_{r}, \lambda_{r}^{-1}$ be the eigenvalues of $P_{\gamma^r}$, so that
$det (I - P_{\gamma^r}) = 2 - (\lambda_r + \lambda^{-1}_{r}).$
Now, if the eigenvalues of $P_{\gamma}$ are $\{e^{\pm i \alpha}\}$
(in the elliptic case) then those of $P_{\gamma^r}$ are $\{e^{\pm
i r \alpha}\}$, hence  $ det (I - P_{\gamma^r}) = 2- 2 \cos r
\alpha.$ Similarly for the hyperbolic case. The formulae then
follows form Corollary \ref{DETP}.

\end{proof}

We now consider the inverse Hessian ${\mathcal H}_+ =
H_{2r}^{-1}$, which will be important in the calculation of wave
invariants. We denote its matrix elements by $h^{pq}_+$. We also
denote by ${\mathcal H}_-$ the matrix in which the roles of $a, b$
are interchanged; it is the inverse Hessian of ${\mathcal L}_-$.

\begin{prop} \label{HESSPAR}  The diagonal matrix elements $h^{pp}_{+}$ are
constant when the parity of $p$ is fixed, and we have:
$$\begin{array}{llll} p \;\; \mbox{odd} \;\; \implies & h^{pp}_+ = h^{11}_+,\;\;\; &
 p \;\; \mbox{even} \;\; \implies & h^{pp}_+ = h^{22}_+
 \\ & & & \\
p \;\; \mbox{odd} \;\; \implies & h^{pp}_- = h^{11}_-,\;\;\; &  p \;\; \mbox{even} \;\; \implies & h^{pp}_- = h^{22}_-, \\ & & & \\
h^{11}_+ = h^{22}_-, & h^{22}_+ = h^{11}_- & & .
\end{array}$$
\end{prop}

\begin{proof}

Indeed, let us  introduce the cyclic shift operator on $\R^{2r}$
given by $P e_j = e_{j + 1}$, where $\{e_j\}$ is the standard
basis, and where $P e_{2r} = e_1.$ It is then easy to check that $
P {\mathcal H}_+ P^{-1} = {\mathcal H}_-,$ hence that $ P
{\mathcal H}_+^{-1} P^{-1} = {\mathcal H}_-^{-1}.$ Since $P$ is
unitary, this says
$$h^{pq}_- = \langle {\mathcal H}_- ^{-1} e_p, e_q \rangle =  \langle  P {\mathcal H}_+^{-1} P^{-1} e_p, e_q \rangle =
 \langle   {\mathcal H}_+^{-1} P^{-1} e_p, P^{-1} e_q \rangle = h_+^{p -1, q - 1}.$$
 It follows that the matrix ${\mathcal H}_{\pm}$ is
invariant under even powers of the shift operator, which shifts
the indices $j \to j + 2 k$ ($k = 1, \dots, r$). Hence, diagonal
matrix elements of like parity are equal.
\end{proof}

\section{\label{WTEI} Resolvent trace invariants}

We now formulate the
  key results (Theorems
 \ref{SETUPA}-\ref{SETUPA}) expressing localized wave traces as   oscillatory
 integrals over the boundary with special phases and amplitudes.  We then tie these statements  together with the statements in
 Theorem 1.1 (v) of \cite{Z5}.

 First, we state a general  result, largely contained in \cite{Z4, Z5} which expresses the
 localized resolvent trace as a finite sum of special oscillatory
 integrals. For simplicity we only state it for the $r$th iterate
 of a bouncing ball orbit.

\begin{theo} \label{MAINCOR2} Suppose that $r L_{\gamma}$ is the only length
in the support of $\hat{\rho}$. Then  for each  order $k^{-R}$ in
the trace expansion of Corollary (\ref{MAINCOR}), we have
$$ Tr 1_{\Omega} R_{B \rho}^{\Omega} (k + i \tau) = \sum_{\pm}\;
  \sum_{M: 2r \leq M \leq R + 2r}\; \sum_{ \sigma: |\sigma| \leq R,  M - |\sigma|= 2 r}  I_{M,
\rho}^{\sigma, w_{\pm} }(k)\; + O(k^{-R}),$$ where $\sigma$ runs
over all maps $\sigma: \{ 1, \dots, M\} \to \{0, 1\}$, and where
$I_{M, \rho}^{\sigma, w_{\pm}}(k)$ are oscillatory integrals of
the form
\begin{equation}\label{EXPRESSION} \begin{array}{lll} I_{M, \rho}^{\sigma, w_{\pm}}(k)&
= &
 \!\!\!\int_{[-\epsilon, \epsilon]^{2r}}\!\!\!\!\;\;\;\;

e^{i k  {\mathcal L}_{w_{\pm}}(  x_1,  \dots, x_{2r})}\hat{\rho}({\mathcal L}_{w_{\pm}}(  x_1,  \dots, x_{2r})) \\ & & \\
 & \times &   a_{M, \rho}^{\sigma,
w_{\pm}}
  (k,   x_1,  x_2, \dots,  x_{2r})
d x_1 \cdots dx_{2r}.\end{array} \end{equation} Here, ${\mathcal
L}_{w_{\pm}}$ is given in (\ref{LPM}) and $a_{M, \rho}^{\sigma,
w_{\pm}}$ are certain semi-classical amplitudes (cf.
(\ref{AMPS})). The asymptotics are negligible unless $M -
|\sigma|= 2 r$ and then the order of $I_{M, \rho}^{\sigma,
w_{\pm}}(k)$ equals $- |\sigma|$.

 \end{theo}

 It follows that only  a finite number of terms $  I_{M, \rho}^{\sigma, w_{\pm}}(k)$ contribute to each order in
$k$ in the expansion in Corollary \ref{MAINCOR}:

 \begin{cor} We have:
$$\begin{array}{l}\sum_{\pm}  \sum_{M: 2 r \leq M  \leq R + 2r }  I_{M, \rho}^{\sigma, w_{\pm}}(k)
 \sim  {\mathcal D}_{B, \gamma}(k + i \tau) \; \sum_{j =
0}^{R} B_{\gamma; j} \;  k^{-j} + O(k^{-R})\end{array},
$$ where $B_{\gamma; j}$ are the Balian-Bloch invariants of the union of
the periodic orbits $\gamma$, and ${\mathcal D}_{B, \gamma}(k + i
\tau)$ is the symplectic pre-factor of (\ref{PR}).

\end{cor}

\subsection{Proof of Theorem \ref{MAINCOR2}}

As mentioned above, most of the proof is contained in \cite{Z4,
Z5}. For the sake of completeness, we sketch the key elements of
the proof.

 We follow  the path originated by Balian-Bloch and followed in many physics articles (see e.g. \cite{BB1, BB2, AG}). It starts from
 the exact formula (of Fredholm-Neumann),
\begin{equation} \label{POT} R_B^{\Omega}(k + i \tau) =
 R_0(k + i \tau) - 2 \; {\mathcal D} \ell(k + i \tau) (I +  N(k + i \tau))^{-1} r_{\Omega}
 {\mathcal S} \ell^{tr}(k + i \tau) \end{equation}
 for the resolvent with given boundary conditions.
 Here,
 ${\mathcal D} \ell (k + i \tau)$ (resp. ${\mathcal S} \ell (k + i \tau)$) is the double
 (resp. single) layer potential, $ {\mathcal S}^{tr}(k + i \tau)$ is the
 transpose,   and
$N(k + i \tau)$ is the boundary integral operator  on
$L^2(\partial \Omega)$ induced by  ${\mathcal D} \ell (k + i
\tau)$. Also, $R_0(k + i \tau)$ is the free resolvent on $\R^2,$
and $r_{\Omega}$ is the restriction to the boundary.
 The Schwartz kernel of the boundary integral operator is given
 by plus (in the Dirichlet case) or minus (in the Neumann case)
\begin{equation} \label{blayers} N(k + i \tau)f(q) =  2 \int_{\partial
\Omega}\frac{\partial}{\partial \nu_y} G_0(k + i \tau, q, q')
f(q') ds(q'),
\end{equation} where $G_0(\lambda, x, y)$ is the free Green's function (resolvent kernel)
on $\R^2$, where   $ds(q)$ is the arc-length measure on $\partial
\Omega$, where $\nu$ is the interior unit normal to $\Omega$, and
where $\partial_{\nu} = \nu \cdot \nabla$. The free Green's kernel
has an exact formula in terms of Hankel functions (\ref{NHANK}),
which gives a WKB approximation to $N(k + i \tau)$ away from the
diagonal. Its phase is the boundary distance function
$d_{\Omega}(q, q')$, indicating that $N(k + i \tau)$ is the
quantization of the billiard map.

But as  discussed extensively in \cite{Z5, Z4, HZ}, $N(k + i
\tau)$ is not a classical Fourier integral operator, but  is
rather a non-standard kind of hybrid Fourier integral operator.
Near the diagonal,  it is a  homogeneous pseudo-differential
operator of order $-1$ (in dimension two it is actually
 of order
$-2$ as proved in \cite{Z5}, Proposition 4.1), while away from the
diagonal it is
 a semi-classical Fourier integral operator of
order $0$ which quantizes the billiard map. To separate out these
two Lagrangian submanifolds (which intersect along tangent vectors
to the boundary),  we introduce a cutoff $\chi(k^{1 - \delta} |q -
q'| )$ to the diagonal, where $\delta > 1/2$ and where $\chi \in
C_0^{\infty}(\R)$ is a cutoff to a neighborhood of $0$. We then
put
\begin{equation} \label{N01}  N(k + i \tau) = N_0(k + i \tau) + N_1(k + i
\tau), \;\; \mbox{with}  \end{equation}
\begin{equation} \label{N01DEF} \left\{\begin{array}{l}  N_0(k + i
\tau, q, q') = \chi(k^{1 - \delta} |q - q'| ) \;  N(k + i \tau, q,
q'), \\ \\ N_1(k + i \tau, q, q') = (1 - \chi(k^{1 - \delta} |q -
q'| ))\; N(k + i \tau, q, q'). \end{array} \right.
\end{equation}
As proved in \cite{Z5, Z4, HZ},   $ N_1((k + i \tau), q, q') ) $
is a semiclassical Fourier integral operator of order $0$ with
phase equal to the boundary distance function $d_{\partial
\Omega}(q,q')$. The diagonal part $N_0$ is of order $-1$ (in fact,
of order $-2$ \cite{Z5}) and therefore plays a secondary role.

We now relate the expansion (\ref{PR}) of the regularized
resolvent trace to that for $\log \det N(k + i \tau)$. This
relation
 has  already been proved in \cite{EP, C,
Z4} in somewhat different ways.

The clearest proof is to
 combine the interior boundary  problem $\Delta^{\Omega}_B$ with a complementary exterior
boundary problem $\Delta_{B'}^{\Omega^c}$. Since we are only
dealing here with Dirichlet or Neumann boundary conditions, we do
not define the term `complementary' but only use the term to
indicate the special cases  $B = D, B' = N$ or $B = N, B' = D.$ We
therefore introduce the exterior {\it Green's kernel}
$G_{B'}^{\Omega^c}(k + i \tau, x, y) \in \mathcal D'(\Omega^c
\times \Omega^c)$ with boundary condition $B$, namely the kernel
of the exterior resolvent and is the unique solution of the
boundary problem:
\begin{equation}\label{GREENKEXT}
\left\{ \begin{array}{l}  - (\Delta_{B'}^{\Omega^c} + (k + i
\tau)^2) G_{B'}^{\Omega^c}(k + i \tau, x, y) =
\delta_{y}(x),\;\;\;
(x, y \in \Omega^c) \\ \\
 B' G_{B'}^{\Omega^c}(k + i \tau, x, y) = 0, \;\;\; x \in
\partial \Omega^c
\\ \\
 \frac{\partial G_{B'}^{\Omega^c}(k \!+\! i \tau, x, y)}{\partial r} - i (k
\!+\! i \tau) G_{B'}^{\Omega^c}(k \!+\! i \tau, x, y) =
o(\frac{1}{r}) ,\;\; \mbox{as}\;\; r \to \infty . \end{array}
\right.
\end{equation}

 We now  combine the interior and exterior operators with complementary boundary
 conditions $B, B'$  into the direct sum  $R_B^{\Omega}(k + i \tau)
\oplus R^{\Omega^c}_{B'}(k + i \tau).  $ For simplicity, we only
consider $B = D, B' = N$.  For $\hat{\rho} \in
C_0^{\infty}(\R^+)$, we put
\begin{equation}
R_{\rho B}^{\Omega}(k + i
\tau) \oplus R_{\rho B'}^{\Omega^c}(k + i \tau) \\[6pt]
 \hspace{20pt} = \int_{\R} \rho(k - \mu) (\mu + i \tau) \left[ R_B^{\Omega}(\mu +
i \tau)\oplus R_{ B'}^{\Omega^c}(\mu + i \tau) \right] d\mu.
\end{equation}
The purpose of combining the interior/exterior resolvents is
revealed in the following proposition, which equates the
 trace of the direct sum resolvent to the Fredholm determinant of the
boundary integral operator. It is proved in \cite{Z4} and closely
related statements are proved in \cite{EP, C}. The operator $N$ is
defined in (\ref{blayers}) in the Dirichlet case. In general it
depends on the boundary conditions $B, B'$. We follow the notation
of \cite{T} except that we multiply the $N$ of \cite{T} by
$\frac{1}{2}$ to simplify some notation.

\begin{prop} \label{MAIN} For any $\tau > 0$, the  operator $ (I + N(k + i \tau))$ has a well-defined Frehdolm determinant $\det (I + N(\lambda +
i \tau))$, and we have:
$$\begin{array}{l}
Tr_{\R^2} [R_{ \rho D}^{\Omega}(k + i \tau) \oplus
R_{ \rho N}^{\Omega^c}(k + i \tau)  - R_{0 \rho}(k + i \tau)]\\[6pt]
\hspace{43pt} = \int_{\R} \rho(k - \lambda)  \frac{d}{d \lambda}
\log \det (I + N(\lambda + i \tau)) d\lambda. \end{array}$$
Further, for  $\tau
> 0, \log \det (I + N(k + i \tau))$ is  differentiable
in $k$,  $ (I + N(k + i \tau))^{-1} N'(k + i \tau)$ is of trace
class and we have:
$$\frac{d}{d k} \log \det(I + N(k + i \tau)) = Tr_{\partial \Omega} (I
+ N(k + i \tau))^{-1} N'(k + i \tau). $$

\end{prop}

This proposition reduces wave trace expansions to the boundary.
Indeed, the direct sum resolvent is related to the direct sum wave
groups as in (\ref{RTWT}):
\begin{equation} \label{RW} R_{\rho B}^{\Omega}(k + i \tau)  \oplus R_{\rho B'}^{\Omega^c}(k + i \tau)
= \int_0^{\infty} \hat{\rho}(t) e^{i (k + i \tau) t} \left[
E_{B}^{\Omega}(t)\oplus E_{B'}^{\Omega^c}(t) \right] dt.
\end{equation} The trace of the direct sum wave group
$E_{B}^{\Omega^c}(t)\oplus E_{B'}^{\Omega}(t)$ has a singularity
expansion as in (\ref{Wavetrace}) which sums over interior and
exterior periodic orbits. As in (\ref{PR}), it may be restated in
terms of the direct sum resolvent: Let $\gamma$ be a
non-degenerate interior billiard trajectory whose length
$L_{\gamma}$ is isolated and of multiplicity one in $Lsp
(\Omega)$.  Let $\hat{\rho} \in C_0^{\infty}(L_{\gamma} -
\epsilon, L_{\gamma} + \epsilon)$, equal to one on $(L_{\gamma} -
\epsilon/2, L_{\gamma} + \epsilon/2)$ and with no other lengths in
its support.  Then the interior trace $Tr R_{B \rho}^{\Omega}(k +
i \tau)$ and the exterior trace $Tr [R_{B'\rho}^{\Omega^c}(k + i
\tau) - R_{0 \rho}(k + i \tau)]$ admit complete asymptotic
expansions of the form
\begin{equation} \label{ASYMP} \left\{ \begin{array}{l}
 Tr [R_{B'\rho}^{\Omega^c}(k + i \tau) - R_{0 \rho}(k + i
\tau)] \sim  {\mathcal D}_{B, \gamma}(k + i \tau)\;
 \sum_{j = 0}^{\infty}
B_{\gamma, j}\;
k^{-j} \\[9pt]
 Tr R_{B \rho}^{\Omega}(k + i \tau) \sim {\mathcal D}_{B, \gamma}(k + i \tau) \; \sum_{j
= 0}^{\infty} B_{\gamma, j} \; k^{-j},
\end{array} \right. \end{equation}
 whose  coefficients
$B_{\gamma; j}$ are the Balian-Bloch resolvent trace  invariants
 of periodic (internal, resp. external) billiard
orbits. We can therefore sum the two expansions to produce one for
the direct sum. The coefficients depend on the choice of boundary
condition but we do not indicate this in the notation.

Combining the results, we get:

\begin{cor} \label{MAINCOR} Suppose that $L_{\gamma}$ is the only length
in the support of $\hat{\rho}$. Then,
$$\begin{array}{l} \int_{\R} \rho(k -
\lambda) \frac{d}{d \lambda} \log \det (I + N(\lambda + i \tau)) d
\lambda \\  \\
= \int_{\R} \rho(k -
\lambda) Tr_{\partial \Omega} (I + N(\lambda + i \tau))^{-1} N'(\lambda + i \tau) d\lambda \\
\\ \sim  {\mathcal D}_{B, \gamma}(k + i \tau)\;
 \sum_{j
= 0}^{\infty} B_{\gamma, j} \; k^{-j}
\end{array},
$$ where as above $B_{\gamma; j}$ are the Balian-Bloch invariants of the union of
the periodic orbits $\gamma$ of length $L_{\gamma}$ of the
interior and exterior problems in (\ref{ASYMP}).

 \end{cor}

 In proving
the remainder estimate and the expansion in Proposition \ref{ASY},
we further microlocalize the result to the (interior) orbit
$\gamma$. This will select out the wave invariants of the desired
interior orbit $\gamma$. A periodic orbit of the billiard flow
corresponds to a periodic point of the billiard map $\beta$.
  To microlocalize to this
periodic orbit we introduce a semiclassical pseudodifferential
cutoff operator $\chi_{0} (\phi, k^{-1} D_{\phi})$. In the case of
a bouncing ball orbit, it has
 complete symbol
$\chi(\phi, \eta)$ supported in $V_{\epsilon}:= \{(\phi, \eta):
|\phi|, |\eta| \leq \epsilon\}$.

 \begin{prop}\label{CUTOFF2}  Suppose that $\gamma$ is a bouncing ball orbit, whose length  $L_{\gamma}$ is the only length
in the support of $\hat{\rho}$. Let $\chi_0$ be a cutoff operator
to the endpoints of $\gamma$. Then,
  $$\begin{array}{l}Tr \rho *  (I + N(k + i \tau ))^{-1} \circ \frac{d}{dk} N(k
+ i \tau) \\ \\ \sim Tr \rho *  (I + N(k + i \tau ))^{-1}  \circ
\frac{d}{dk} N(k + i \tau) \circ \chi_0(k).\end{array}$$

\end{prop}

 We will use the formula in Corollary \ref{MAINCOR}, as modified in Proposition \ref{CUTOFF2},   to calculate
  the $B_{\gamma; j}$ modulo remainders which are
 inessential for the inverse spectral problem. To do so, we now
 express the left hand side (for each order of singularity $k^{-j}$) as a finite  sum of oscillatory
 integrals $I_{M, \rho}^{\sigma, w}$ (see (\ref{EXPRESSION})) plus a remainder which is of lower order than $k^{-j}.$

 To define the oscillatory integrals $I_{M, \rho}^{\sigma, w}$,  we first
 expand  $(I +  N(\lambda + i \tau))^{-1}$ in a finite geometric
series plus remainder,
\begin{equation} \label{GS} (I \!+\! N(\lambda + i \tau))^{-1} = \sum_{M = 0}^{M_0} (-1)^M \; N(\lambda + i\tau)^M
+
 (-1)^{M_0 + 1} \; N(\lambda + i \tau)^{M_0 + 1}  (I
\!+\! N(\lambda + i \tau))^{-1},  \end{equation}
 and  prove that, in calculating a given
order of Balian-Bloch invariant $B_{\gamma, j}$, we may neglect a
sufficiently high remainder.

\begin{prop} \label{ASY}  For each  order $k^{-J}$ in the trace expansion
of Corollary (\ref{MAINCOR}) there exists $M_0(J)$ such that $$
\begin{array}{ll} (i) & \sum_{M = 0}^{M_0} (-1)^M  Tr  \int_{\R}
\rho(k - \lambda)\; N(\lambda + i \tau)^M N'(\lambda + i \tau) d \lambda \\
& \\  &  =  {\mathcal D}_{B, \gamma}(k + i \tau) \;
  \sum_{j = 0}^{J} B_{\gamma, j} \; k^{-j} + O(k^{-J - 1}) ,
\\&
\\(ii) & Tr \int_{\R} \rho(k - \lambda)  N(\lambda + i \tau)^{M_0 + 1}  (I
\!+\! N(\lambda \!+\! i \tau))^{-1} N'(\lambda + i \tau) d \lambda
= O(k^{-J - 1}).
\end{array} $$
The same holds after composition with $\chi_0(k)$.
\end{prop}

The proof of this Proposition is one of the principal results in
\cite{Z5, Z4}. In \cite{Z5} the result is stated in Theorem 1.1
(iii), while the remainder trace is estimated in \S 8. The version
stated in Proposition \ref{ASY} is proved in \S 5 of \cite{Z4}. It
is
 simpler than Theorem 1.1 (iii) of \cite{Z5}  because the interior integral analyzed in \S 7
of that paper is eliminated in the  reduction to the boundary.

It simplifies the formula somewhat to integrate the derivative by
parts onto $\hat{\rho}$,  since it eliminates the derivative in
the special factor $N'(\lambda + i \tau).$

\begin{cor} \label{ASYCOR}  For each  order $k^{-J}$ in the trace expansion
of Corollary (\ref{MAINCOR}) there exists $M_0(J)$ such that $$
\begin{array}{ll} (i)
&  \sum_{M = 0}^{M_0} \frac{(-1)^{M} }{M + 1}    Tr  \int_{\R}
\rho'(k - \lambda) N(\lambda + i \tau)^{M+1} d
 \lambda \\ & \\  &  =  {\mathcal D}_{B, \gamma}(k + i \tau) \;  \sum_{j =
0}^{J} B_{\gamma, j} \; k^{-j} + O(k^{-J - 1}) ,
\\&
\\(ii) & Tr \int_{\R} \rho(k - \lambda)  N(\lambda + i \tau)^{M_0 + 1}  (I
\!+\! N(\lambda \!+\! i \tau))^{-1} N'(\lambda + i \tau) d \lambda
= O(k^{-J - 1}).
\end{array} $$
The same holds after composition with $\chi_0(k)$.
\end{cor}

The next step is to prove that  the terms in Proposition
\ref{ASY}(i) may be expressed as oscillatory integrals (see
(\ref{EXPRESSION})). This is not obvious, as mentioned above,
since the $N$ operator is not a Fourier integral kernel. As
indicated in (\ref{N01})-(\ref{N01DEF}), we handle this problem by
breaking up $N$ as a sum $N = N_0 + N_1$ of two terms, where $N_0$
has the singularity on the diagonal of a pseudodifferntial
operator of order $-2$ (cf. \cite{Z5}, Proposition 4.1), and where
$N_1$ is manifestly an oscillatory integral operator of order $0$
with phase $|q(\phi) - q(\phi')|$. As mentioned above, and as
discussed in detail in \cite{Z4, HZ}, the phase is a generating
function of the billiard map, so the $N_1$ term is a quantization
of $\beta.$

We thus write,
\begin{equation} \label{BINOMIAL} (N_0 + N_1)^M = \sum_{\sigma: \{ 1, \dots, M\} \to \{0, 1\}}
N_{\sigma(1)} \circ N_{\sigma(2)} \circ \cdots \circ
N_{\sigma(M)}.
\end{equation} In \cite{Z5} \S 6, we regularized the terms by proving a composition law for
products $N_0 \circ N_1, N_1 \circ N_0$. The main technical point
is that the amplitudes of $N_0, N_1$ belong to  the  symbol class
$S^p_{\delta}({\bf T})$ where ${\bf T}$ is the unit circle
parameterizing $\partial \Omega$,  consisting of symbols $a(k,
\phi)$ which satisfy:
\begin{equation}\label{SYMBOL}
|(k^{-1} D_{\phi})^{\alpha} a(k, \phi)| \leq C_{\alpha} |k|^{p -
\delta |\alpha|},\;\;\; (|k| \geq 1).  \end{equation} This follows
from the classical formula (see e.g. \cite{Z5} \S 4; \cite{AG},
(2.2))
\begin{equation}
\label{NHANK} \begin{array}{lll}  N(k + i \tau, q(\phi_1),
q(\phi_2)) & = & - \frac{i}{4} (k + i \tau) H^{(1)}_1 ((k + i
\tau) |q(\phi_1) - q(\phi_2)|) \\ & & \\ & \times & \cos
\angle(q(\phi_2) -q (\phi_1), \nu_{q(\phi_2)}) \end{array}
\end{equation}
 for $N$ in terms of Hankel functions and from the
asymptotics of Hankel function $H^{(1)}_1$. We recall that
 the Hankel function of index $\nu$ has the integral representations (\cite{T}, Chapter 3, \S 6)
\begin{equation} \label{HANKEL} \begin{array}{lll} H^{(1)}_{\nu} (z)
 & = &  (\frac{2}{\pi z})^{1/2} \frac{e^{i(z - \pi \nu/2 -  \pi/4)}}{ \Gamma(\nu + 1/2)}
\int_0^{\infty} e^{-s} s^{-1/2} (1 - \frac{s}{2 i z})^{\nu -1/2}
ds,
\end{array} \end{equation}
from which  it follows that
 $H^{(1)}_1$ admits an asymptotic expansion as its argument tends to
 infinity of the form
\begin{equation} \label{HANKASYMPT} H^{(1)}_{1}(t) \sim  e^{it - \frac{3 \pi i}{4}}
t^{-\frac{1}{2}} \sum_{j=0}^\infty c_j t^{ - j}, \; (t \to \infty)
\end{equation}
where $c_0 = \sqrt{2/\pi}$. Moreover, the expansion can be
differentiated term by term. We set:
\begin{equation} \label{DEFa1} a_1(t) =  \frac{\sqrt{\frac{2}{\pi}}}{\Gamma(\frac{3}{2})} \int_0^{\infty} e^{-s} s^{-1/2} (1 - \frac{s}{2 i
t})^{1/2} ds, \end{equation} so that
\begin{equation} \label{HANKASYMPT2}
H^{(1)}_{1}(t) \sim  e^{it - \frac{3 \pi i}{4}} t^{-\frac{1}{2}}
a_1(t)
\end{equation}
We observe that   $ a_1$ is a complex valued semi-classical symbol
of order $0$ of $z \in \R_+$
 in the sense that (cf. (\ref{SYMBOL}))
 $$ (1 - \chi(k^{1 - \delta} z)) a_1((k + i \tau)z) \in
S^{0}_{\delta} (\R_z).  $$  We then have
\begin{equation} \label{HSYMBOL}\!\!  (k \!+\! i \tau) H^{(1)}_1( (k \!+\! i \tau) z) =
(\frac{k \!+\! i \tau}{z})^{ \frac{1}{2}} e^{i(k \!+\! i \tau) z}
a_1((k \!+\! i \tau)z),  \hspace{-15pt}\end{equation} hence

 $$\begin{array}{l}  N_1(k + i \tau, q(\phi_1), q(\phi_2)) = (1 - \chi(k^{1 -
\delta} (\phi_1 - \phi_2) )) \\ \\ \cdot (\frac{k + i
 \tau}{|q(\phi_1) - q(\phi_2)|})^{\frac{1}{2}} a_1(k + i \tau, q(\phi_1),
q(\phi_2)) e^{i (k + i \tau) |q(\phi_1) - q
(\phi_2)|}\end{array}$$ with \begin{equation} a_1(k + i \tau,
q(\phi_1), q(\phi_2)) :=  a_1 ((k + i \tau)|q(\phi_1) -
q(\phi_2)|) \cos \vartheta_{1,2} \in S^{0}_{\delta}({\bf T}^2),
\end{equation} where $\vartheta_{1,2} = \angle q(\phi_2) -
q(\phi_1), \nu_{q(\phi_2)} ). $

 The main conclusion
is that  $N_0 N_1$ and $N_1 N_0$ are semiclassical Fourier
integral operators with the same phase as $N_1$, but with an
amplitude of one lower degree in $k$. This allowed us to remove
all of the factors of $N_0$ from each of these terms except for
the term $N_0^M$.  Each remaining term except for $N_0^M$ is a
Fourier integral operator on ${\bf T}^m$ for some $m \leq M$, with
phase given by the length functional (\ref{LENGTH}) and with
amplitude in the symbol class $S^p_{\delta}({\bf T}^m)$ for some
$p$, consisting of symbols $a(k, \phi_1, \dots, \phi_m)$ which
satisfy the analogue of (\ref{SYMBOL}):
\begin{equation}\label{SYMBOL2}
|(k^{-1} D_{\phi})^{\alpha} a(k, \phi)| \leq C_{\alpha} |k|^{p -
\delta |\alpha|},\;\;\; (|k| \geq 1).  \end{equation} Because each
removal of $N_0$ drops the order by one,  the term $N_1^M$ is of
the highest order in the sum. A later estimate on traces shows
that $N_0^M$ does not contribute to the trace asymptotics (see
\cite{Z5}, \S 9.0.7).

We summarize the result as follows. Let us rewrite the terms of
(\ref{BINOMIAL}) as
 \begin{equation} N_{\sigma} : =  N_{\sigma(1)} \circ
N_{\sigma(2)} \circ \cdots \circ N_{\sigma(M)}
\end{equation}  and set
\begin{equation}    |\sigma| = \; \#
 \sigma^{-1} (0) = \; \mbox{the number of $N_0$ factors occurring in}\; N_{\sigma}.\end{equation}
In \cite{Z5}, Propositions 6.1, we show that the regularized
compositions are semiclassical Fourier integral kernels.

\begin{prop}  \label{NPOWERFINAL} We have: \medskip

 \noindent{\bf (A)} Suppose that $N_{\sigma}$  is not of the form $N_0^M$.
 Then for any integer  $R >0$, $N_{\sigma} \circ \chi_0 (k + i \tau)$ may be expressed as the sum
 $$N_{\sigma}
= F_{\sigma}(k, \phi_1, \phi_2) + K_R,  $$  where $F_{\sigma}$ is
a semiclassical Fourier integral kernel of order $- |\sigma|$
associated to $\beta^{M - |\sigma|}$ of the form
\begin{equation} \label{F}  F_{\sigma}(k, \phi_1, \phi_2) =
e^{i (k + i \tau) |q(\phi_{1}) - q(\phi_{2 })|} A_{\sigma}(k,
\phi_1, \phi_2),
\end{equation} where  $ A_{\sigma}(k,
\phi_1, \phi_2)$ is a semi-classical amplitude,  and where the
remainder $K_R$ is a bounded smooth kernel which is uniformly of
order $k^{-R}$.
\medskip

\noindent{\bf (B)} $N_0^M \circ \chi_0 \sim N_{0M} \circ \chi_0,$
where $N_{0M}$ is a semiclassical pseudodifferential operator of
order $-M$. (For the notation $\chi_0$ see Proposition
(\ref{CUTOFF2}).)

\end{prop}

As a corollary of Proposition \ref{NPOWERFINAL}, we obtain the
following preliminary form for the trace as a sum of oscillatory
integrals. It is a simplification of \cite{Z5},  Lemma 9.2 in that
we do not need any interior integrals.

\begin{cor} \label{SOCINT}
 $Tr\;  \rho ' \;*\;
N_{\sigma} \circ \chi_0$ is an  oscillatory integral of the form
$$ \begin{array}{l} I_{M, \rho}^{\sigma} (k) =  k^{ (M - |\sigma| + 3)/2} \int_{\R} \int_{\R} \int_{{\bf T}^{M -
|\sigma|}} e^{i k [(1 - \mu) t + \mu {\mathcal L}_{\sigma}
(q(\phi_1), \dots, q(\phi_{M - |\sigma|} ))]} e^{- \tau \log k
{\mathcal
L}_{\sigma} (q(\phi_1), \dots, q(\phi_{M - |\sigma|}))} \\ \\
\chi(\overline{q(\phi_1) - q(\phi_2)}, \phi_1) A_M^{\sigma}(k \mu,
\phi_1, \dots, \phi_{M - |\sigma|}) \hat{\rho'}(t) dt d\mu d\phi_1
\cdots d \phi_{M - |\sigma|},
\end{array} $$
where  $\chi(\overline{q(\phi_1) - q(\phi_2)}, \phi_1)$ is the
value at the vector $(q(\phi_1), \overline{q(\phi_1) -
q(\phi_2)})$  of  a cutoff $\chi$  to a microlocal neighborhood in
$B^*\partial \Omega$  of the direction of the bouncing ball orbit,
where
$${\mathcal L}_{\sigma} (q(\phi_1), \dots, q(\phi_{M - |\sigma|}
))
= |q(\phi_1) - q(\phi_2)| + \cdots + |q(\phi_{M - |\sigma|}) -
q(\phi_1)|, $$ and where $A_M^{\sigma} (k, \phi_1, \dots, \phi_{M
- |\sigma|}) \in S^{-|\sigma|  }_{\delta}.$
\end{cor}

\subsubsection{Completion of proof of Theorem \ref{MAINCOR2}}

We now complete the proof of Theorem \ref{MAINCOR2}. To obtain our
final form for the oscillatory integrals, we make some further
simplifications. For simplicity of exposition, and because it is
our main application,
  we specialize to a bouncing ball
orbit. In view of Propositions \ref{MAIN} and \ref{ASY}, it
suffices to prove:

\begin{prop} \label{MAINCOR2PROP} Suppose that $r L_{\gamma}$ is the only length
in the support of $\hat{\rho}$.  Then  for each  order $k^{-R}$ in
the trace expansion of Corollary (\ref{MAINCOR}), we have
$$\int_{\R} \rho(k -
\lambda) \frac{d}{d \lambda} \log \det (I + N(\lambda + i \tau)) d
\lambda \sim  \sum_{\pm}  \sum_{M: 2r \leq M \leq R + 2r} \sum_{
\sigma: |\sigma| \leq R,  M - |\sigma|= 2 r}  I_{M, \rho}^{\sigma,
w_{\pm} }(k)\; + O(k^{-R}),$$ where the oscillatory integrals
$I_{M, \rho}^{\sigma, w_{\pm}}(k)$ are as in Theorem
\ref{MAINCOR2}.

 \end{prop}

 \begin{proof}

  The first observation is that the
regularized integral $ I_{M, \rho}^{\sigma}(k + i \tau) $ of
Corollary \ref{SOCINT} has no critical points   unless $M -
|\sigma|= 2 r$ (where $r L_{\gamma}$ is the unique length in the
support of $\hat{\rho}$). We will refer to these oscillatory
integrals as contributing. Since each $T_{\epsilon}$ has two
pieces, each  contributing integral can be written as  a sum of
$2^{2r}$ terms  $ I_{M, \rho}^{\sigma, w}(k + i \tau) $,
corresponding to a choice of an element $w$ of
$$\{\pm\}^{2r}
 := \{w : \Z_{2r}
  \to \{\pm \} \}.$$
The length functional in Cartesian coordinates for a given
assignment $w$ of signs is given by
\begin{equation}\label{LW}  {\mathcal L}_{w} (x_1, \dots, x_{2r}) = \sum_{j =
1}^{2r } \sqrt{(x_{j \!+\! 1} \!-\! x_j)^2 \!+\! (f_{w(j \!+\! 1)
}(x_{j \!+\!1}) \!-\! f_{w(j)}(x_j))^2}.
\end{equation}
Here, $x_{2r
 + 1} = x_1.$

 We further observe that $ I_{M, \rho}^{\sigma, w}(k + i \tau) $ has
 no critical points  unless $w(j)$ alternates between $+$ and $-$ as $j$ increases.
Otherwise, $ I_{M, \rho}^{\sigma, w}(k + i \tau) $ is negligible
as $k \to \infty$. Thus, only two $w$ count asymptotically, which
we denote by $w_{\pm}.$ The corresponding length functionals are
given in (\ref{crit}) and their Hessians are given in Proposition
\ref{KT}.

In these remaining oscillatory integrals, we then
 eliminate the $(t, \mu)$
variables in the integral displayed in Corollary \ref{SOCINT} by
stationary phase. The Hessian in these variables is easily seen to
be non-degenerate, and the Hessian operator equals $-
\frac{\partial^2}{k \partial t
\partial \mu}.$   The
amplitude depends on $t$ only in the factor $\hat{\rho'}(t).$
Since $\hat{\rho'}(t) = t \hat{\rho}(t)$ and since $\hat{\rho}$ is
assumed  to be constant in some interval $(r L_{\gamma} -
\epsilon, r L_{\gamma} + \epsilon)$, $t \hat{\rho}(t)$ is locally
linear and therefore  only the zeroth order and $(-1)$st order
terms
$${\mathcal L} \hat{\rho}({\mathcal L}) A_M^{\sigma}(k, x) + \frac{k}{i k}
\hat{\rho}({\mathcal L}) \frac{\partial A_M^{\sigma} (k,
x)}{\partial k}$$
 in the stationary phase expansion are non-zero.  In the second term, the $k$
 in the denominator comes from the Hessian operator and the $k$ in the numerator comes from
 the $\mu$- derivative of the amplitude. After replacing the $dt d\mu$
 integral by this stationary phase expansion,  we arrive at the final form of
the oscillatory integrals (\ref{EXPRESSION}) given in the Theorem,
with amplitude
\begin{equation}  \label{AMPS} a_{M}^{\sigma,
w_{\pm}}
  (k,   x) = {\mathcal L}_{w_{\pm}} A_M^{\sigma}(k, x) + \frac{1}{i}
 \frac{\partial A_M^{\sigma}(k,
x)}{\partial k}(k, x). \end{equation}

\end{proof}

\section{\label{PRIN} Principal term of the Balian-Bloch trace}

In this section, we state and begin the proof of a  key result for
the proof of Theorems \ref{ONESYM} and \ref{DISYM}. It singles out
a single oscillatory integral (the principal term) from Theorem
\ref{MAINCOR2} which generates all terms of the wave trace (or
Balian-Bloch) expansion which contain maximal number of
derivatives of the boundary defining function per power of $k$
(i.e. order of wave invariant). As mentioned in the introduction,
the other terms will turn out to be redundant for domains in our
symmetry classes.

To clarify this notion of generating all the highest derivative
terms, we define it formally. Below, ${\mathcal J}^s$ denotes the
$s$-jet.

\begin{defn} \label{MODULO} Let $\gamma$ be an $m$-link periodic
reflecting ray, and let $\hat{\rho} \in C_0^{\infty}(\R)$ be a
cut off
 satisfying supp $\hat{\rho} \cap Lsp(\Omega) = \{ r L_{\gamma}
 \}$ for some fixed $r \in \N$. Given an oscillatory integral
 $I(k)$, we write
 $$Tr 1_{\Omega} R_{B \rho}^{\Omega} (k + i \tau) \equiv I(k)\;\;
 \mbox{mod} \;\; {\mathcal O}( \sum_{j} k^{-j} ({\mathcal J}^{2j - 2} \kappa))$$
 if
 $$Tr 1_{\Omega} R_{B \rho}^{\Omega} (k + i \tau) - I(k)$$
 has a complete asymptotic expansion of the form (\ref{PR}), and
 if the coefficient of $k^{-j}$ depends on $\leq 2j - 2$
 derivatives of the curvature $\kappa$ at the reflection points.
\end{defn}

For the sake of clarity, we state the next result only in the
simplest case of a bouncing ball orbit. The statement is similar
for any non-degenerate $m$-link periodic reflecting ray. The
description of the properties of phase and amplitude are repeated
from \cite{Z4} for the sake of self-completeness. For terminology
concerning billiard trajectories, we refer to \S \ref{BACKGROUND}.

 \begin{theo}\label{SETUPA}  Let $\gamma$ be a primitive
non-degenerate $2$-link periodic
 reflecting ray, whose reflection points are points of non-zero curvature of $\partial \Omega$,
 and let $\hat{\rho} \in C_0^{\infty}(\R)$ be a  cut off
 satisfying supp $\hat{\rho} \cap Lsp(\Omega) = \{ r L_{\gamma}
 \}$ for some fixed $r \in \N$.  Orient $\Omega$ so that $\gamma$ is
 the vertical segment $\{x = 0\} \cap \Omega$, and so that
  $\partial \Omega$ is a union of two graphs over $[-\epsilon, \epsilon]$. Then in the sense of Definition \ref{MODULO},  we have
\begin{equation}\label{EXPRESSIONA} \begin{array}{l} Tr 1_{\Omega} R_{B \rho}^{\Omega} (k + i \tau)
\equiv  \!\!\!\;\;\; \sum_{\pm} \int_{[-\epsilon, \epsilon]^{2
r}}\!\!\!\!\;\;\; e^{i (k + i \tau) {\mathcal L}_{\pm}(  x_1,
\dots, x_{ 2r})}
  \hat{\rho}  ({\mathcal L}_{\pm}(  x_1,  \dots, x_{ 2r})) a_{\pm, r}^{pr}(k,   x_1,  x_2, \dots,  x_{2r})
d x_1 \cdots dx_{2 r},\\ [6pt]  \end{array}
\hspace{-5pt}\end{equation}
 where the phase $
 {\mathcal L}_{\pm}
  (x_1, \dots, x_{2r})$
 is given in (\ref{LPM}), and where the amplitude
 is given by:
 $$a_{\pm, r}^{pr}(k, x_1, \dots, x_{2r}) =  {\mathcal L}_{w_{\pm}} A^{pr}_{\pm, r}(k, x_1, \dots, x_{2r})
+ \frac{1}{i}
 \frac{\partial}{\partial k}  A^{pr}_{\pm, r} (k,  x_1,  \dots, x_{2r}), $$

where

\begin{equation}\label{AMPLSINGL}
\begin{array}{lll}
 A^{pr}_{\pm, r}(k, x_1, \dots, x_{2r}) & = &     \Pi_{p = 1}^{2 r }\;
 (
 \frac{a_1((k \!+\! i \tau) \sqrt{(x_{p }
\!- \!x_{p\!+\!1})^2 \!+\! (f_{w_{\pm}(p)}(x_{p }) \!-\!
f_{w_{\pm}(p \!+ \!1)}(x_{p\!+\!1} ))^2}}{\left( (x_{p } \!-
\!x_{p\!+\!1})^2 \!+\! (f_{w_{\pm}(p)}(x_{p }) \!-\! f_{w_{\pm}(p
\!+ \!1)}(x_{p\!+\!1} )^2 \right)^{1/4}} \\ && \\ & \times &
 \ \frac{(x_{p} -
x_{p+1}) f'_{w_{\pm}(p)} (x_{p}) -  ( f_{w_{\pm}(p)} (x_{p}) -
f_{w_{\pm}(p+1)}(x_{p+1}) )}{ \sqrt{(x_{p} - x_{p+1})^2 +
(f_{w_{\pm}(p)} (x_{p}) - f_{w_{\pm}(p+1)}(x_{p+1}) )^2}} )
\end{array} \end{equation} where $a_1$ is the Hankel amplitude in
(\ref{HSYMBOL}). Here, as above, $x_{2r + 1} = x_1$. \end{theo}

Theorem \ref{SETUPA} is a crucial ingredient in the proof of
Theorem \ref{ONESYM}. It gives explicit formulae for the phase and
amplitude of the principal oscillatory integrals that determine
the highest order jet of $\Omega$ in each wave invariant. The
notation $A^{pr}_r, a^{pr}_r$ refers to the amplitude of the
principal terms of the $2r$th integral; these amplitudes contain
terms of all orders in $k$ and principal here does not refer to
the principal symbol, i.e. the leading order term in the
semi-classical expansion.  The calculation of the highest
derivative terms of the Balian-Bloch wave invariants uses only
some key properties of the phase and principal amplitude which may
be derived directly from the formulae in Theorem \ref{SETUPA}.
They are detailed in \S \ref{KEYPROPS}.

The proof of theorem \ref{SETUPA} requires two main steps:
\begin{enumerate}

\item Identification  of two main terms in Theorem \ref{MAINCOR2},
the {\it principal terms}, which generate the highest derivative
data, and proof that the amplitude and phase have the stated form.

\item Proof that non-principal terms contribute only lower order
derivative data.

\end{enumerate}

We now  define the principal terms. In \S \ref{KEYPROPS}, Lemma
\ref{AMPPROPS}, we prove that their phases and amplitudes have the
stated form. We further  describe the properties of the phase and
amplitude which will be used in the proof of Theorem \ref{ONESYM},
and   tie the statement of Theorem \ref{SETUPA} together with the
corresponding statement in \cite{Z5}. The fact that non-principal
terms do not contribute highest order derivative data to a given
Balian-Bloch invariant requires the analysis of the stationary
phase expansions in the next section and is  given in \S
\ref{NONPRINCIPAL}.

\begin{defn}\label{PRINCIPAL}  Let $\gamma$ be a $2$-link periodic orbit.
The {\em principal   terms} are the completely regular terms
$I_{2r, \rho}^{\sigma_0, w_{\pm}}$ coming from $N_1^{2r}$,i.e.
with $M = 2r$ and  with $\sigma_0(j) = 1$ for all $j$.  The two
terms correspond to the two possible orientations $w_{\pm}(j)$, of
the  $2r$th iterate of the bouncing ball orbit.
\end{defn}

In other words, the  principal terms are simply those coming from
the term
\begin{equation} \label{NONEK} Tr\; \rho*\;
 N_1^{2 r}(k)  \circ N_1'(k) \circ \chi(k) \end{equation}
 in the expansion (\ref{BINOMIAL}).

We observe that in fact, the two principal terms are equal. This
is not surprising, since a bouncing ball orbit is reciprocal.

\begin{prop} \label{PM} We have: $I_{2r, \rho}^{\sigma_0, w_+} (k) = I_{2r, \rho}^{\sigma_0, w_-}
(k).$ \end{prop}

\begin{proof} We permute the  variables $x_j$ according
to the cyclic permutation $s$ of their indices:
$$s = \left( \begin{array}{lllll} 1 & 2 & \cdots & r - 1 & r \\ &&&&\\
2 & 3 & \cdots & r  & 1 \end{array} \right) $$ in the integral in
(\ref{EXPRESSION}).  Since $w_+(s(j)) = w_-(j)$, this takes
${\mathcal L}_{-} \to {\mathcal L}_{+}$ and  $a^0_- \to a^0_ +$ in
(\ref{AMPLSINGL}). Indeed, ${\mathcal L}_{\pm}$ (resp.
$a^0_{\pm}$) are sums (resp. products) of terms of the form $F(x_p
- x_{p + 1}, f_{w_{\pm}(p)}(x_p) - f_{w_{\pm}(p+1)}(x_{p + 1})). $
Cyclically shifting the index by one moves each term (resp.
factor) to the next except that it does change the index
$w_{\pm}(p)$. Hence, it changes the sum (resp. product)  only by
shifting  $w_+$ to $w_-$ (and vice-versa).

\end{proof}

Henceforth, we often  omit $ I_{2r, \rho}^{\sigma_0, w_-} (k)$ and
multiply  $I_{2r, \rho}^{\sigma_0, w_+} (k)$ by $2$.

\subsection{\label{KEYPROPS} Key properties of the principal amplitude and phase}

We first prove that the phase and amplitude of the principal
oscillatory integrals have the form stated in Theorem
\ref{SETUPA}, and establish a few consequences. After that, we
assemble all of the properties used in the proof of Theorem
\ref{ONESYM}. In the following, we abbreviate $\lcal_+ =
\lcal_{w_+}$. We  use the notation $D_{x_p} =
\frac{\partial}{\partial x_p}$ and multi-index notation for its
powers.

\begin{lem}\label{AMPPROPS}  The phase and principal amplitude of the
principal oscillatory integrals $I_{2r, \rho}^{\sigma_0, w_{\pm}}$
have the following properties:

$$\begin{array}{l} (i)~~~  \mbox{In  its dependence on the boundary
defining functions}\; f_{\pm}, \; \mbox{ the amplitude}\;
a^{pr}_{+, r} \; \mbox{
 has the
form}
\\[2pt]
\hspace{21pt}  \alpha_r (k, x,  f_{\pm}, f_{\pm}').
\\  \\
(ii)~~  \mbox{As above, in its dependence on}\;  x\;   \\ [4pt]
 \qquad a^{pr}_{+,r
 } (k, x_1, \dots, x_{2r})   =  {\mathcal L}_{+} A^{pr}_{+, r}(k, x_1, \dots, x_{2r})
+ \frac{1}{i}
 \frac{\partial }{\partial k} A^{pr}_{+, r} (k,   x_1,  \dots, x_{2r}), \; \mbox{where} \\ [4pt] A^{pr}_{+,
r} (k, x_1, \dots, x_{2r}) \; = \; \Pi_{p = 1}^{2r} A_p (x_{p},
x_{p + 1})\;\;\;\; (2 r + 1 \equiv 1)
 \\ \\
(iii)\mbox{At the critical point, the principal amplitude has the
asymptotics}\\ \\  \; a^{pr}_{+, r}(k, 0) \sim (2 r L) L^{-r}
\acal_r(0) + O(k^{-1}), \;\; \mbox{where}\; \acal_r(0)\;
\mbox{depends only on}\; r\; \mbox{and not on} \; \Omega;
\\ \\\\ \\
(iiia)  \; \frac{a^{pr}_{+, r}(k, 0) e^{i (k + i \tau) {\mathcal
L_+}(0) + i \pi/4 sgn Hess {\mathcal L_+}(0)} }{\sqrt{\det Hess
{\mathcal L}_+}} \sim (2 r L)\;{\mathcal A}_r(0)\; {\mathcal
D}_{B, \gamma}(k
+ i \tau) (1 +  O(k^{-1})) \; (cf. \ref{PR});\\ \\
(iv)  \nabla a^{pr}_{+, r} (k, x_1, \dots, x_{2r})|_{x = 0} = 0.\\ \\
(v)~     D_{x_p}^{(2j - 1)} {\mathcal L}_+ |_{x=0} \equiv 2 w_+(p)
f_{w_+(p)}^{(2j - 1)}(0)\;\;\; \mbox{mod} \;\; R_{2r} ({\mathcal
J}^{2j - 2} f_+(0), {\mathcal
J}^{2j - 2} f_-(0)), \\ \\
(v.a) \;\;D_{x_p}^{(2j)} {\mathcal L}_+ |_{x=0} \equiv 2 w_+(p)
f_{w_+(p)}^{(2j)}(0) \; \mbox{mod} \; R_{2r} ({\mathcal J}^{2j -
1} f_+(0), {\mathcal J}^{2j - 2} f_-(0)),
\end{array}  $$
where $\equiv$ in general means equality modulo lower order
derivatives of $f$.

\end{lem}

\begin{proof}

The oscillatory integrals $I_{2r, \rho}^{\sigma_0, w_{\pm}}$ have
the form (\ref{EXPRESSION}) with  the phases  ${\mathcal L}_{\pm}$
(\ref{LW}), and by Proposition \ref{PM} it suffices to consider
the $+$ term.

 Formula (ii) for the amplitude follows from the
general description of the amplitudes of all the oscillatory
integrals $I_{M, \rho}^{\sigma, w}$ in the proof of Theorem
\ref{MAINCOR2} (cf. (\ref{AMPS})). The factors $A_{\pm, r}^{pr}$
of the amplitudes of $I_{2r, \rho}^{\sigma_0, w_{\pm}}$ are given
in (\ref{AMPLSINGL}).

  The  further properties of
the phase and amplitude stated in Lemma \ref{AMPPROPS} may be read
off directly from the formula in (\ref{AMPLSINGL}). Statements
(i)-(ii)  are visible from the formula. At $x = 0$, the leading
order term of the  principal amplitude in $k$ equals $2r L $ (from
the factor ${\mathcal L}$) times $ L^{-r}$ from  the
$t^{-\frac{1}{2}}$ factor in the Hankel asymptotics
(\ref{HANKASYMPT})-(\ref{HANKASYMPT2}) times a coefficient
$\acal_r(0)$ which depends on $r$ but not on $\Omega$ and which is
due to additional factors in the asymptotics of the  free Green's
function $G_0$: namely, a product of $2r$ factors of $
\sqrt{\frac{2}{\pi}} e^{\frac{3 \pi i}{4}} $ from the principal
term of the Hankel amplitude  $a_1$ (loc. cit.),  factors of
$\frac{-i}{4}$ in the relation between the free Green's function
$G_0$  and the Hankel function (\ref{NHANK}), factors of $2$ in
the relation of $N(k + i \tau)$ and $G_0$ (\ref{blayers}).
 We do not  need to know  $\acal_r(0)$ or other universal factors
  explicitly, since they multiply all terms in
the expansion. Statement (iiia) gives the principal term in the
stationary phase expansion at $x = 0$ and relates the Hessian
determinant and $L^{-r}$ to the Poincar\'e determinant as  in
Propositions \ref{KTa} and \ref{DetHESS} (see also \cite{AG},
(3.17)). The second term is of order $k^{-1}$, so will not
contribute to the highest derivative term in a given wave
invariant.

From the fact that $x = 0$ is a critical point of $f_{\pm}$ and
$(x_j - x_{j-1})^2$ we get
\begin{equation} \label{TAB3} \left\{
\begin{array}{ll} (a) & \nabla_x \left(\sqrt{(x_{p } \!- \!x_{p\!+\!1})^2 \!+\!
(f_{w_{\pm}(p)}(x_{p }) \!-\! f_{w_{\pm}(p \!+ \!1)}(x_{p\!+\!1}
)^2}) \right)|_{x = 0} = 0 \\ & \\
(b) &  \nabla_x \left( \frac{(x_{p} - x_{p+1}) f'_{w_{\pm}(p)}
(x_{p}) - ( f_{w_{\pm}(p)} (x_{p}) - f_{w_{\pm}(p+1)}(x_{p+1}) )}{
\sqrt{(x_{p} - x_{p+1})^2 + (f_{w_{\pm}(p)} (x_{p}) -
f_{w_{\pm}(p+1)}(x_{p+1}) )^2}} \right) |_{x = 0} = 0
\end{array}\right.,
\end{equation}

which implies

\begin{equation} \label{TAB4} \begin{array}{ll}
\nabla_x\; a_{+, r}^{pr} |_{x = 0}  &=  \nabla_x\; D_k a_{+,
r}^{pr} |_{x = 0} = 0.
\end{array}
\end{equation}

Statement (v) on the phase holds because
\begin{equation} \label{ODDDER}\left\{ \begin{array}{lll} D_{x_p}^{(2j - 1)} {\mathcal L}_{+} |_{x = 0}  & \equiv & \sum_{\pm} ((x_p
- x_{p \pm 1})^2 +
 (f_{w_+(p)}(x_p) -\  f_{w_+(p \pm 1) }(x_{p \pm 1}))^2)^{-1/2} \\ & &
 \\
 &\times &
(f_{w_+(p) }(x_p)  -\ f_{w_+(p \pm 1) }(x_{p \pm 1})) f_{w_+(p)
}^{(2j - 1)}(x_p) |_{x = 0} \;  \mbox{mod} \; R_{2r} ( {\mathcal
J}^{2j - 2}
f_{\pm}(0)),\\ && \\
D_{x_p}^{(2j)} {\mathcal L}_+ |_{x = 0} &\equiv & \sum_{\pm} ((x_p
- x_{p \pm 1})^2 +
 (f_{w_+(p)}(x_p) -\  f_{w_+(p \pm 1) }(x_{p \pm 1}))^2)^{-1/2} \\ & &
 \\
 &\times &
(f_{w_+(p) }(x_p)  -\ f_{w_+(p\pm 1) }(x_{p \pm 1})) f_{w_+(p)
}^{(2j)}(x_p)|_{x = 0}  \; \mbox{mod} \; R_{2r} ({\mathcal J}^{2j
- 1} f_+(0), {\mathcal J}^{2j - 2} f_-(0)).
\end{array}\right. \end{equation}
We make the crucial  observation that the  $\pm $ terms are equal
(and especially, do not cancel !), giving the factor of $2$ in (v)
since $f_{w_+(p)}(0) -\  f_{w_+(p \pm 1)}(0) = w_+(p) L $.

\end{proof}

\subsubsection{\label{FURTHERPROPSAMP} Further properties of the amplitude and phase}

We continue the discussion of the amplitude by detailing the other
special values of the phase and amplitude at the critical point
that are  used in the  \S \ref{FD} in the course of proving
Theorem \ref{ONESYM}. Although the value of the discussion will
only become clear in \S \ref{FD}, it seems best  to  give the
details at this point.

\begin{enumerate}

\item In the proof of  Lemma \ref{LEM5.4}(i), we  use that
\begin{equation} \label{USEDINLEM5.4i}   D_{x_p}^{2j - 2} a^{pr}_{+, r} |_{x = 0} \equiv 0 \;
\mbox{mod} \; R_{2r} ( {\mathcal J}^{2j - 2} f_{\pm}(0)), \;
(\forall p = 1, \dots, 2r).
\end{equation}  Indeed, by the explicit formula of   (\ref{AMPLSINGL})
one can only obtain the higher derivative  $f_{\pm}^{2j - 1}(0)$
by applying all $2j - 2$ derivatives on the term $f'_{w_{\pm}(p)}
(x_{p})$  in
$$
 \ \frac{(x_{p} -
x_{p+1}) f'_{w_{\pm}(p)} (x_{p}) -  ( f_{w_{\pm}(p)} (x_{p}) -
f_{w_{\pm}(p+1)}(x_{p+1}) )}{ \sqrt{(x_{p} - x_{p+1})^2 +
(f_{w_{\pm}(p)} (x_{p}) - f_{w_{\pm}(p+1)}(x_{p+1}) )^2}}.$$ But
then the accompanying factors of $x_{2p} - x_{2p + 1}$ vanish at
the critical point.

\item In the proof of  Lemma \ref{LEM5.4}(ii), we  use that

\begin{equation} \label{USEDINLEM5.4ii}  D_{x_p}^{(2j - 1)} D_{x_q}{\mathcal L}
\equiv 0 \; \mbox{mod} \; R_{2r} ( {\mathcal J}^{2j - 2}
f_{\pm}(0)), \; (\forall p = 1, \dots, 2r, \forall q \not= p).
\end{equation}
Indeed, in  (\ref{ODDDER})  $D_{x_p}^{(2j - 1)} {\mathcal L}$ is
displayed as a product of two factors. Since $q \not= p$,  the
 derivative $D_{x_q}$ must be applied to the factor $$((x_p -
x_{p + 1})^2 +
 (f_{w_+(p)}(x_p) -\  f_{w_+(p+1) }(x_{p + 1}))^2)^{-1/2}
(f_{w_+(p) }(x_p)  -\ f_{w_+(p+1) }(x_{p + 1})),$$ which vanishes
at $x = 0$ for any $q$.

\item In the same Lemma \ref{LEM5.4}, we also use that the only
non-vanishing third derivatives of ${\mathcal L}$ at $x = 0$ are
pure third derivatives in one variable $D_{x_j}^3 {\mathcal L}$.
Indeed, from (\ref{crit}), we see that only mixed derivatives
using two consecutive indices (say, $x_j, x_{j + 1}$) can be
non-zero. However, we have:
\begin{equation} \label{THREEDERIVSL} D_{x_j}^2 D_{x_{j + 1}}
{\mathcal L} |_{x = 0} = 0 = D_{x_j} D_{x_{j + 1}}^2 {\mathcal L}
|_{x = 0}. \end{equation} Since the identities are similar, we
only consider the first, which is equivalent to
$$D_{x_j} D_{x_{j + 1}} \frac{(x_j - x_{j + 1}) + (f_{w_{\pm}(j)}(x_j) - f_{w_{\pm}(j +
 1)}
 (x_{j + 1})) f_{w_{\pm}(j)}'(x_j)}{\sqrt{(x_j - x_{j + 1})^2 +
(f_{w_{\pm}(j)}(x_j) - f_{w_{\pm}(j + 1)}(x_{j + 1}) )^2}} |_{x =
0} = 0. $$ We write the fraction as $\frac{F(x_j, x_{j +
1})}{G(x_j, x_{j + 1})}, $ and note that
$$D_{x_j} D_{x_{j + 1}}  \frac{F}{G} |_{x = 0} = \frac{D_{x_j} D_{x_{j + 1}}
F}{G} |_{x = 0} \;\; \mbox{if}\; F(0) = \nabla G (0) = 0. $$ When
$F = (x_j - x_{j + 1}) + (f_{w_{\pm}(j)}(x_j) - f_{w_{\pm}(j +
 1)}
 (x_{j + 1})) f_{w_{\pm}(j)}'(x_j)$, we also have
 $D_{x_j} D_{x_{j + 1}}
F |_{x = 0} = 0. $

\item Further, we use that, for all $p$, $D_{x_p}^3 \lcal_+ (0) =
2 w_+(p)  f_{w_+(p)}'''(0). $ Indeed, as in the calculation of the
higher derivatives in Lemma \ref{AMPPROPS}, there are two terms,
and each (in the notation above) has  the form $\frac{D_{x_p}^2
F(0)}{G(0)}$. To obtain a non-zero term,  the two derivatives must
fall on the factor $f_{w_{\pm}(p)}'(x_p)$, and thus we get
$$\begin{array}{lll} D_{x_p}^{(3)} {\mathcal L}_+ (0)  &= & \sum_{\pm} ((x_p - x_{p
\pm 1})^2 +
 (f_{w_+(p)}(x_p) -\  f_{w_+(p \pm 1) }(x_{p \pm 1}))^2)^{-1/2} \\ & &
 \\
 &\times &
(f_{w_+(p) }(x_p)  -\ f_{w_+(p\pm 1) }(x_{p \pm 1})) f_{w_+(p)
}^{(3)}(x_p)|_{x = 0} \\  && \\
& = & 2 w_+(p)  f_{w_+(p) }^{(3)}(0). \end{array}$$ Again, we
observe that the $x_{p \pm 1}$ terms agree; therefore they add
rather than cancel.

\end{enumerate}

\subsection{Comparison with \cite{Z5}}

For the sake of completeness, we  tie together the statement of
Theorem \ref{SETUPA} with the corresponding statement (v) of
Theorem 1.1 of \cite{Z5} and with \cite{Z4}:

\medskip

{\it Theorem 1.1 (v) of \cite{Z5}: Let  $\gamma$  be  a primitive
non-degenerate $m$-link periodic
 reflecting ray of length $ L_{\gamma}$,  and let $\hat{\rho} \in C_0^{\infty}(\R)$ be a cut off
 satisfying supp $\hat{\rho} \cap Lsp(\Omega) = \{ r L_{\gamma}
 \}$ for some fixed $r \in \N$.   Then
modulo an error term  $R_{2r} ({\mathcal J}^{2j - 2} \kappa
(a_j))$ depending only on the $(2j - 2)$-jet of curvature $\kappa$
of $\partial \Omega$ at the $m$ reflection points $a_j$ of
$\gamma$, the wave invariant $ B_{\gamma^r, j - 1} +
B_{\gamma^{-r}, j - 1}$ can be obtained by applying stationary
phase to  the  oscillatory integral
$$Tr\;  \rho*\;
 N_1^{m r} \circ \chi(k) \circ {\mathcal S} \ell(k + i
 \tau)^{tr} \circ  {\mathcal D} \ell(k + i
 \tau).$$}

In Theorems \ref{MAINCOR2} and \ref{SETUPA}, we have followed
  \cite{Z4} in combining
  the interior and exterior problems. Taking the trace then  eliminates the single
  and double layer potentials ${\mathcal S} \ell$ resp. ${\mathcal D}
  \ell$ in Theorem 1.1(v) of \cite{Z5}, allowing for the reduction of the trace to the boundary in  (\ref{NONEK}).

\section{\label{FD} Feynman diagrams in inverse spectral theory}

In this section, we use the oscillatory integrals in Theorems
  \ref{SETUPA}  to obtain explicit
formulae for the highest derivative terms of the  wave trace
invariants at a bouncing ball orbit in terms of the curvature
function of the boundary. To our knowledge, these are the first
explicit formulae. In the next section it will be proved that
lower order derivative data is redundant for domains with our
symmetries.

For simplicity we restrict to bouncing ball orbits. There are
similar results for general periodic reflecting rays (see Lemma
\ref{BGAMMAJDI} for the dihedral case).  We first state the result
for domains without symmetries, and then specialize to mirror
symmetric domains in Corollary \ref{BGAMMAJSYM}. We use the graph
parametrization rather than the curvature in the formulae. In the
following, $h_+^{pq}$ are the matrix elements of the inverse
Hessian $Hess({\mathcal L}_+)^{-1}$ of the positively oriented
length functional ${\mathcal L}_+ = {\mathcal L}_{w_+}$ of
(\ref{crit}) and (\ref{LW}) in the principal terms.

\begin{theo} \label{BGAMMAJ} Let $\Omega$ be a smooth domain with a bouncing
ball orbit   $\gamma$  of length $ r L_{\gamma}$. Then there exist
polynomials $p_{2, r, j}(\xi_1, \dots, \xi_{2j + 1}; \eta_1,
\dots, \eta_{2j + 1})$ which are
 homogeneous of degree $-j$ under the dilation $f \to \lambda f,$ which are  invariant under the substitutions
$\xi_j \iff - \eta_j$ and under $f(x) \to f(-x)$ such that:

\begin{itemize}

\item $B_{\gamma^r, j}   = p_{2, r, j}(f_-^{(2)}(0), f_-^{(3)}(0),
\cdots, f_-^{(2j + 2)}(0); f_+^{(2)}(0), f_+^{(3)}(0), \cdots,
f_+^{(2j + 2)}(0)).$

\item   In the Balian-Bloch (resolvent trace) expansion  of
Corollary \ref{MAINCOR} and in (\ref{ASYMP}),
 the data $f^{(2j)}_{\pm}(0), f^{(2j - 1)}_{\pm
 }(0)$ appear first in the $k^{-j + 1}$st order term,
 and then only in the expansion of the principal terms;

\item This coefficient has the form
$$ \begin{array}{l}  B_{\gamma^r, j - 1}  \equiv 4 r L {\mathcal A}_0(r) \{ 2 (w_{{\mathcal G}_{1,
j}^{2j, 0}}) ((h^{11}_{+, 2r} )^j f^{(2j)}_{+}(0) - (h^{22}_{+,
2r})^j f^{(2j)}_{-, 2r}(0))  \\ \\ +  4 \sum_{q, p = 1 }^{2r}
[(w_{{\mathcal G}_{2,  j + 1 }^{2j - 1, 3,
 0}}) (h^{pp}_{+})^{j
- 1} h^{qq}_{+, 2r} h^{pq}_{+, 2r} + (w_{\widehat{{\mathcal
G}}_{2, j + 1 }^{2j - 1, 3,
 0}}) (h^{pp}_{+, 2r})^{j - 2}
 (h^{pq}_{+, 2r})^3] w_+(p) w_+(q) f_{w_+(p)}^{(2j - 1)}(0)
 f_{w_+(q)}^{(3)}(0)\}
\\ \\ + R_{2r} ({\mathcal
J}^{2j - 2} f_+(0), {\mathcal J}^{2j - 2} f_-(0)),\end{array} $$
where the remainder $ R_{2r} ({\mathcal J}^{2j - 2} f_+(0),
{\mathcal J}^{2j - 2} f_-(0))$ is a polynomial in the designated
jet of $f_{\pm}.$ Here, $w_+(p) = (-1)^{p+1}$ and  as in the
introduction, $w_{\gcal} = \frac{1}{|Aut (\gcal)|}$ are
combinatorial factors independent of $\Omega$ and $r$.

\end{itemize}
\end{theo}

Where possible, we have simplified the sums using Proposition
\ref{HESSPAR}. The top even derivative term is calculated in Lemma
\ref{LEM5.3} and the top odd derivative is cacluated in Lemma
\ref{LEM5.4}.

The  methods we use to make the calculations  could be also  used
to evaluate the oscillatory integrals in Theorem  \ref{MAINCOR2}
and the wave invariants to all orders of derivatives. This could
be useful in the inverse spectral problem for general domains
without symmetry. However, we are content here to study the
highest derivative terms and apply the results to domains with
symmetry.

We prove Theorem  \ref{BGAMMAJ} by making a  stationary phase
analysis of the oscillatory integrals in Corollary \ref{MAINCOR2}.
As mentioned in the introduction, our  strategy involves a novel
feature of the stationary phase expansion, namely to separate out
the terms of the stationary each order in $k$ which have the
maximum number of derivatives of the boundary defining function or
equivalently of its curvature.

 Since the formulae (\ref{SPEXP})-
(\ref{MSPJTERM}) are very complicated, we organize the
calculations  by the diagrammatic method. Since Feynman diagrams
have not been used before in inverse spectral theory, we digress
to present the fundamentals of  the diagrammatic approach to the
stationary phase expansion;  clear expositions are given in
 \cite{A, E} (see also \cite{AG}).

\subsection{\label{SPFD} Stationary phase diagrammatics}

We consider a general  oscillatory integral
$$Z_k = \int_{\R^n} a(x) e^{ik S(x)} dx$$
where $a \in C_0^{\infty}(\R^n)$ and where $S$ has a unique
critical point in supp $a$  at $0$. We write $H$ for the Hessian
of $S$ at $0$ and $R_3$ for the third order remainder in its
Taylor expansion at $x = 0$:
$$S(x) = S(0) + \langle H x, x\rangle/2 + R_3(x).$$
The stationary phase expansion is:
$$\begin{array}{l} Z_k = (\frac{2\pi}{k})^{n/2} \frac{e^{i \pi sgn (H)/4}}{\sqrt{|det H|}} e^{i k S(0)} Z_k^{h \ell}, \;\; \mbox{where}\\ \\
Z_k^{h \ell} = [a(\frac{\partial}{\partial J}) e^{i k R_3(\frac{\partial}{\partial J})}]_{J = 0} e^{-\frac{1}{2 i k}
\langle J, H^{-1} J \rangle} \\ \\
= \sum_{I = 0}^{\infty} \sum_{V = 0}^{\infty}
[a(\frac{\partial}{\partial J}) [\frac{i k}{V!}
(R_3(\frac{\partial}{\partial J}))^V]_{J = 0} \frac{[-\frac{1}{2 i
k} \langle J, H^{-1} J \rangle]^I}{I!}. \end{array} $$ The
graphical analysis of the stationary phase expansion consists in
the identity \begin{equation} \label{GRAPHID}
[a(\frac{\partial}{\partial J}) [\frac{i k}{V!}
(R_3(\frac{\partial}{\partial J}))^V]_{J = 0} \frac{[-\frac{1}{2 i
k} \langle J, H^{-1} J \rangle]^I}{I!} = \sum_{(\gcal, \ell)  \in
G_{V, I}} \frac{I_{\ell} (\gcal)}{|Aut(\gcal)|}
\end{equation} where $G_{V, I}$ is the class of labelled graphs
$(\gcal,  \ell)$ with $V$ {\it closed vertices} of valency $\geq
3$ (each corresponding to the phase), with one {\it open vertex}
(corresponding to the amplitude), and with $I$ edges. The function
$\ell$ `labels' each end of each edge of $\gcal$ with an index $j
\in \{1, \dots, n\}.$

\begin{rem} The term ` open vertex' is equivalent to `
marked' or `external' vertex in some texts,   and is graphed here
as an unshaded circle. A ` closed' vertex is the same as an
`unmarked' or `internal'  vertex and is graphed as a shaded
circle. Also, it is non-standard  to include the labels $\ell$  in
the notation for Feynman amplitudes; we do so because in our
problems certain labels are distinguished.
\end{rem}

Above, $|Aut(\gcal)|$  denotes  the order of the automorphism
group of $\gcal$, and  $I_{\ell} (\gcal)$ denotes the `Feynman
amplitude' associated to the labelled graph  $(\gcal, \ell)$. By
definition, $I_{\ell}(\gcal)$ is obtained by the following rule:
To each edge with end labels $m,n$ one assigns a factor of
$\frac{-1}{ik} h^{mn}$ where as above $H^{-1} = (h^{mn}).$ To each
closed vertex one assigns a factor of $i k \frac{\partial^{\nu} S
(0)}{\partial x^{i_1} \cdots
\partial x^{i_{\nu}}}$ where $\nu$ is the valency of the vertex
and $i_1 \dots, i_{\nu}$ at the index labels of the edge ends
incident on the vertex. To the open vertex, one assigns  the
factor $\frac{\partial^{\nu} a(0)}{\partial x^{i_1} \dots \partial
x^{i_{\nu}}}$, where $\nu$ is its valence.   Then
$I_{\ell}(\gcal)$ is the product of all these factors.  To the
empty graph one assigns the amplitude $1$.  In summing over
$(\gcal, \ell)$ with a fixed graph $\gcal$, one sums the product
of all the factors as the indices run over $\{1, \dots, n\}$.

We note that the power of $k$ in a given term with $V$ vertices
and $I$ edges equals $k^{\chi_{\gcal'}}$, where $\chi_{\gcal'} = V
- I$ equals the Euler characteristic of the graph $\gcal'$ defined
to be $\gcal$ minus the open vertex. We thus have;
\begin{equation} \label{ZK} Z_k^{h \ell} = \sum_{j = 0}^{\infty}  \{\sum_{(\gcal, \ell):
\chi_{\gcal'} = - j} \frac{I_{\ell} (\gcal)}{|Aut(\gcal)|}\}.
\end{equation} We note that there are only finitely many graphs
for each $\chi$ because the valency condition forces $I \geq 3/2
V.$ Thus, $V \leq 2 j, I \leq 3 j.$

\subsubsection{Stationary phase formula for $I_{M, \rho}^{\sigma, w_{\pm}}$}

Since Feynman diagrams and amplitudes are unfamiliar in wave trace
calculations, we digress to give some details of the proof of
(\ref{GRAPHID}) and to tie it together with the form of the
stationary phase expansion in standard texts  in partial
differential equations (cf. \cite{Ho}I). This latter form can also
be used to corroborate the calculations below.

The  stationary phase of  (\cite{Ho}I, Theorem 7.7.5) reads:
\begin{equation} \label{SPEXP} Z_k   \sim (\frac{2\pi}{k})^{n/2} \frac{e^{\frac{i \pi}{4} sgn H } e^{i k
S(0)}}{\sqrt{|\det H |} } \sum_{j = 0}^{\infty} k^{-j} {\mathcal
P}_j a (0)\end{equation} where
\begin{equation} \label{MSPJTERM} {\mathcal P}_j a(0) = \sum_{\nu - \mu = j}
 \sum_{2 \nu \geq 3\mu} \frac{i^{-j} 2^{-\nu}}{\mu! \nu!}
\langle H^{-1} D, D \rangle^{\nu} ( a  R_3^{\mu} ) |_{x  = 0}
\end{equation}

  In diagrammatic terms,  the pair $(\mu, \nu)$ correspond to
  graphs
with $\nu = I$  edges and $\mu = V$  closed vertices, hence of
Euler characteristic $\mu - \nu = - j$. We note that the factor
$i^{-j}$ is  common to all graphs of Euler characteristic $-j$ and
in our analysis we absorb  into the prefactor. To tie
(\ref{MSPJTERM}) together with (\ref{GRAPHID}),  we sketch the
proof of the latter, following the exposition in \cite{E} in the
case where the amplitude is $\equiv 1$. We outline the procedure
following the notes of Etingof \cite{E}  This special case turns
out to be the most important for the applications in this paper,
since terms with derivatives of the amplitude will not contribute
to the highest order jets in the wave invariants. The notes of
Axelrod \cite{A} give a clear discussion (as above) of the
contribution of the amplitude to the Feynman amplitude.

\begin{prop} \label{TWOVERSIONS} We have:
$$\frac{2^{-\nu}}{\mu! \nu!} \langle H^{-1} D, D \rangle^{\nu} (   R_3^{\mu}
) |_{x = 0} =  \sum_{(\gcal, \ell) \in G_{\nu, \mu}}
\frac{I_{\ell}(\gcal)}{|Aut(\gcal)|}. $$ \end{prop}

\begin{proof}

We need to re-write the left side  as a sum over graphs in
$G_{\nu, \mu}$ (the class of graphs with $\nu$ edges, $\mu$ closed
vertices of valency $\geq 3$).

Let ${\bf n} = (n_0, n_1, \dots)$ be a sequence of non-negative
integers, of which all but a finite number are zero, and let
$G({\bf n})$ denote the set of graphs with $n_0$ $0$-valent
vertices, $n_1$ $1$-valent vertices etc. We are only considering
the case where the amplitude equals one,  one so there are no
external vertices.

We write  $R_3(x)  =  \sum_{m \geq 3}  B_m(x, \dots, x)/m!, $
where $B_m = d^m S(0),$ as a sum of its homogeneous terms. Change
variables $x \to \sqrt{k} x$, write $e^{i k
R_3(\frac{x}{\sqrt{k}}))} = \Pi_m e^{i k
B_m(\frac{x}{\sqrt{k}})/m!} $ and Taylor expand each exponential
to obtain
\begin{equation} \label{NSUM} \begin{array}{l} Z_k = \sum_{{\bf n}} Z_{{\bf n}}, \;\;
\mbox{with}\\ \\Z_{{\bf n}} = \int_{\R^n} e^{i  H(y,y)/2} \Pi_m
\frac{1}{(m!)^{n_m} n_m!} ((i k)^{-\frac{m}{2} + 1} B_m(y,\cdots,
y))^{n_m} dy. \end{array} \end{equation} The integral may be
calculated by Wick's formula. The diagrammatric interpretation
attaches to each factor $iB_m$ a `flower' of valency $m$, i.e. a
closed vertex with $m$ outgoing edges. Thus, the index ${\bf n}$
prescribes a set of $n_m$ flowers of valency $m$. Let  $T$ be the
set of the ends of the outgoing edges of all of the flowers. For
each pairing $\sigma$ of the ends one obtains a graph $\gcal_{{\bf
n}, \sigma}$.

 Associated to each graph is  its Feynman  amplitude
$F_{{\bf n}, \sigma}$. As described above, one labels each end of
each edge of the graph by indices in $\{1, \dots, n\}$,  assigns a
factor of $\frac{-1}{ik} h^{mn}$ to an edge with end labels $m,n$
and  flower (closed vertex) of valency $i$ with end labels
$(x_{n_1}, \dots, x_{n_i})$  one assigns a factor of $i k
\frac{\partial^{i} S (0)}{\partial x^{n_1} \cdots
\partial x^{n_{i}}}$. One multiplies these expressions over all
edges and closed vertices and then sums over all labelings. One
then has
$$Z_{{\bf n}} = \frac{(2 \pi)^{n/2}}{\sqrt{\det H}}  \Pi_m
\frac{1}{(m!)^{n_m} n_m!} k^{-n_m(\frac{m}{2} + 1)} \sum_{\sigma}
F_{{\bf n}, \sigma}. $$

By comparison, in  (\ref{MSPJTERM}), one Taylor expands the full
factor  $e^{R_3}$ to obtain
$$\begin{array}{l} e^{i k R_3(\frac{x}{\sqrt{k}})} =  \sum_{\mu} \frac{1}{\mu!} \left(
i \sum_m k^{-m/2 + 1}  B_m/m!\right)^{\mu} = \sum_{\mu}
\frac{i^{\mu}}{\mu!} \sum_{{\bf n}: |{\bf n}| = \mu} \Pi_m
k^{-n_m(\frac{m}{2} + 1)} {\mu \choose {\bf n}} \frac{
B_m^{n_m}}{(m!)^{n_m} }.
\end{array}$$
Since
\begin{equation} \label{MUFAC}  \frac{1}{\mu!} \sum_{{\bf n}: |{\bf n}| = \mu} {\mu
\choose {\bf n}} \Pi_m  \frac{ B_m^{n_m}}{(m!)^{n_m} } =
 \sum_{{\bf n}: |{\bf n}| = \mu} \Pi_m
\frac{ B_m^{n_m}}{(m!)^{n_m} (n_m)!}, \end{equation} it follows
that
\begin{equation} \label{MUFAC2}
 \frac{2^{-\nu}}{\mu! \nu!} \langle H^{-1} D, D \rangle^{\nu} (   R_3^{\mu}
) |_{x = 0} =  \frac{2^{-\nu}}{ \nu!} \langle H^{-1} D, D
\rangle^{\nu}  \sum_{{\bf n}: |{\bf n}| = \mu} \Pi_m \frac{
B_m^{n_m}}{(m!)^{n_m} (n_m)!}. \end{equation}

For each fixed  ${\bf n}$, the term on the right side for this
${\bf n}$ is the $\nu$th term in the expansion of $Z_{{\bf n}}$
when (as in the proof in \cite{Ho}) one applies the Plancherel
formula  to the integral (\ref{NSUM}) for  $Z_{{\bf n}}$ and
Taylor expands $e^{i H^{-1}(y,y)/2}. $ The $\nu$th   term can be
sifted out  by replacing $H \to \lambda H$ and finding the term of
order $\lambda^{- \nu}$ on each side. Note that $(\mu, \nu)$ are
determined by ${\bf n}$ : Indeed, $\mu = \sum_m n_m$, and since
each outgoing vertex is paired with exactly one other outgoing
vertex to form an edge,  $\nu = \frac{1}{2} \sum_m m n_m. $ We
write $\mu({\bf n}), \nu({\bf n})$ for the these values.  The
$\lambda^{- \nu}$ terms in the sum over ${\bf n}$ with $|{\bf n}|
= \mu$  run over those ${\bf n}$ for which $\nu({\bf n}) = \nu$,
and thus we have
$$  \frac{2^{\nu}}{\nu!} \langle H^{-1} D, D
\rangle^{\nu} \sum_{{\bf n}: |{\bf n}| = \mu} \Pi_m \frac{
B_m^{n_m}}{(m!)^{n_m} (n_m)!} =   \Pi_m \frac{1}{(m!)^{n_m} n_m!}
\sum_{{\bf n}:  |{\bf n}| = \mu, \nu({\bf n}) = \nu, \sigma}
F_{{\bf n}, \sigma}. $$

Finally, as  explained in \cite{E},
$$\sum_{{\bf n}, \sigma} F_{{\bf n}, \sigma} = \sum_{\gcal, \ell} \frac{\Pi_m (m!)^{n_m}
n_m!}{|Aut(\gcal)|} I_{\ell}(\gcal).$$ The same identity holds if
we restrict to  pairings and graphs with $\mu$ vertices and $\nu$
edges. Cancelling common factors, we get
\begin{equation} \label{BASICFD}  \frac{( 2^{-\nu})}{\nu!} \langle
H^{-1} D, D \rangle^{\nu} \sum_{{\bf n}: |{\bf n}| = \mu} \Pi_m
\frac{ B_m^{n_m}}{(m!)^{n_m} (n_m)!} = \sum_{(\gcal, \ell) \in
G(\mu, \nu)} \frac{I_{\ell}({\gcal})}{|Aut(\gcal)|}.
\end{equation} Combining with (\ref{MUFAC2}) completes the proof.

\end{proof}

\subsection{Maximal derivative terms }

We now apply the diagrammatic stationary phase method  to the
oscillatory integrals $I_{M, \rho}^{\sigma, w_{\pm}}$
(\ref{EXPRESSION}). Further, we consider the additional aspect of
extracting from the stationary phase expansion the terms which
involve the highest number of derivatives of the boundary defining
function $f_{\pm}$ in each power of $k^{-1}$. Such terms with the
maximal number of derivatives arise only from special graphs and
from special  terms in the corresponding Feynman amplitudes with
{\it special labelings} of the vertices. This is a non-standard
feature of diagrammatic analysis and indeed depends on the very
special phase and amplitudes in $I_{M, \rho}^{\sigma, w_{\pm}}$. A
further key issue is the dependence on the number of iterates $M$
of the bouncing ball orbit.

For emphasis, we state our objective as follows:

\begin{itemize}

\item Enumerate the diagrams of each Euler characteristic whose
amplitudes contain the maximum number of derivatives of $\partial
\Omega$ among diagrams of the same Euler characteristic. Determine
which vertex labellings produce the maximum number of derivatives.
Then determine the corresponding ``maximal derivative Feynman
amplitudes", i.e. the sums of monomials containing the highest
number of derivatives. We denote them by $I^{\max}(\gcal).$
\end{itemize}

As we will see, only the principal oscillatory integrals of
Definition \ref{PRINCIPAL} give rise to terms in
$I^{\max}(\gcal).$ We use the following notation for the class of
labelled graphs which give rise to two  types of maximal
derivative terms.

\begin{itemize}

\item $G_{\nu, \mu}^{a, b, c} \subset G_{\nu, \mu}$ are the (not
necesssarily unique)  labelled graphs  whose Feynman amplitude
contains terms of the form $f^{(a)}(0) f^{(b)}(0) a_0^{(c)}(0).$
In fact, we will show that $c = 0$ for all labelled graphs
contributing to the highest number of derivatives of $f$ in a
given order of wave invariant.

\end{itemize}

We denote by  ${\mathcal J}^p$ the operation of extracting the
terms with $p$ derivatives. That is,  ${\mathcal J}^p$ applied to
a monomial in derivatives of the phase is equal to the monomial if
it contains a factor with $p$ derivatives of the phase and zero
otherwise. From Proposition \ref{TWOVERSIONS}, we can evaluate the
combinatorial  coefficients of Feynman amplitudes with a specified
number of derivatives.

\begin{cor} We have:
$$ {\mathcal J}^{p} \frac{2^{-\nu}}{\mu! \nu!}
{\mathcal H}_{\pm}^{\nu} (   R_3^{\mu} ) |_{x_0 = x_1 = \cdots =
x_{2m} = 0}  = \sum_{(\gcal, \ell) \in G(\mu, \nu)}
\frac{{\mathcal J}^p I_{\gcal, \ell}}{|Aut(\gcal)|}. $$
\end{cor}

\subsection{The principal terms}

Our first step is to analyze the stationary phase expansions of
the  {\it principal terms} $I_{2r, \rho}^{\sigma_0, w_{\pm}}(k)$
in the sense of Definition \ref{PRINCIPAL}. By Proposition
\ref{PM} it suffices to consider $w_+$.  We show that the
non-principal terms only contribute lower order derivative data to
the Balian-Bloch invariants $B_{\gamma, j}$. In the next section,
this data will be proved redundant  in the case of the symmetric
domains of this article. As mentioned in the introduction, we only
use the attributes of the phase and amplitude described in Theorem
\ref{SETUPA}. We now use this information to determine where the
data $f_{\pm}^{2j}(0), f^{(2j - 1)}_{\pm}(0)$ first appears in the
stationary phase expansion for the oscillatory integrals.

The only critical point occurs where $x = 0$. We denote by
${\mathcal H}_{\pm}$  the Hessian operator in the variables $(
x_1, \dots, x_{2r})$ at the critical point $ x = 0$ of the phase
${\mathcal L}_{\pm}.$ That is ${\mathcal H}_{\pm} = \langle
Hess({\mathcal L}_{\pm})^{-1} D, D \rangle,$ where $D$ is short
for $(\frac{\partial}{\partial x_{1}}, \cdots,
\frac{\partial}{\partial x_{2r}}).$

\subsubsection{The principal term: The data $f_{\pm}^{2j}(0)$}

We first claim that  $f_{\pm}^{(2j)}(0)$ appears first in the
$k^{-j + 1}$ term in the stationary phase expansion of $I_{2r,
\rho}^{\sigma_0, w_{+}}$. This is because any labelled graph
$(\gcal, \ell) $ for which $I_{\ell}(\gcal)$ contains the factor
$f_{\pm}^{(2j)}(0)$ must have a closed vertex of valency $\geq
2j$,  or the open vertex must have  valency $\geq 2j - 1.$ The
minimal absolute Euler characteristic $|\chi(\gcal
')|$ in the
first case is $j - 1$. Since the Euler characteristic is
calculated after the open vertex is removed, the minimal absolute
Euler characteristic in the second case is $j$ (there must be at
least $j$ edges.) Hence such graphs do not have minimal absolute
Euler characteristic.  More precisely, we have:

\begin{lem} \label{LEM5.3} In the stationary phase expansion of $I_{2r, \rho}^{\sigma_0, w_{+}}$, the
 only  labelled graph $(\gcal, \ell)$ with $-\chi(\gcal') = j - 1$ with
$I_{\ell} (\gcal)$ containing  $f_{\pm}^{(2j)}(0)$ is  given by:
\begin{itemize}

\item ${\mathcal G}_{1, j}^{2j, 0, 0} \in G_{1, j} $  (i.e.  $\mu
= V = 1, I = \nu =  j)$. There is a   unique graph in this class.
It has no open vertex, one closed vertex and $j$ loops at the
closed vertex.

\item
 The only labels producing the desired data are those $\ell_p$ which assign
all  endpoints
of all edges labelled  the same index $p$.

\end{itemize}

The ${\mathcal J}^{2j}$th part of the Feynmann amplitude is
$$I^{\max} ({\mathcal G}_{1, j}^{2j, 0}) = 4 r L\; (w_{{\mathcal
G}_{1, j}^{2j, 0}})\; {\mathcal A}_0(r) \{ \;(h^{11}_+)^j
f^{(2j)}_{+}(0) - (h^{11}_-)^j f^{(2j)}_{-}(0) \},$$  where we
neglect terms with  $\leq 2j - 1$ derivatives.

\end{lem}

We are also interested in the $f_{\pm}^{(2j - 1)}(0)$ terms, but
postpone the calculation of the $f_{\pm}^{(2j - 1)}(0)$ - terms
arising from the diagram ${\mathcal G}_{1, j}^{2j, 0}$ until Lemma
\ref{LEM5.4}(ii) (they turn out to vanish).

\begin{figure}\label{FIGFIVE}
\centerline{\includegraphics{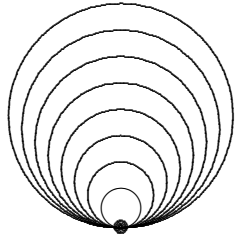}} \caption{${\mathcal
G}_{1, j}^{2j, 0, 0} (- \chi = j - 1,  V = 1, I = j),$ $j$ loops
at one closed vertex. All labels the same. Form of  Feynman
amplitude: $(h_+^{pp})^j D_{x_p}^{(2j)} {\mathcal L}_{+} \equiv
(h_+^{pp})^j f^{(2j)}(0)$}
\end{figure}

\begin{proof}
By (\ref{MSPJTERM}),   the data $f_{\pm}^{2j}(0)$ only occurs in
the term $\mu = 1, \nu = j$ of (\ref{MSPJTERM}). To see this, we
note that
 the Hessian operator ${\mathcal H}^{\nu}_{+}$ associated to ${\mathcal L}_{+}$  has the form
$${\mathcal H}_{+}^{\nu} = \sum_{(i_1, j_1, \dots, i_{\nu}, j_{\nu})} h^{i_1 j_1}_{+} \cdots
h^{i_{\nu} j_{\nu}}_{+} \frac{
\partial^{2\nu}}{\partial x_{i_1} \partial x_{j_1} \dots \partial x_{i_{\nu}} \partial x_{j_{\nu}}}.$$
Any term $(h^{pp}_+ D_{x_p}^2 )^j$ applied to $R_3$ produces a
$f^{(2j)}_{\pm}(0)$ term.

We can also argue non-diagrammatically that no  $\nu_j \geq 2(j + 1)$, i.e.   the power $k^{-j + 1}$ is the greatest power of $k$ in which $f_{\pm}^{(2j)}(0)$ appears.
 Indeed,  it requires $3 \mu$ derivatives to remove the zero
of $R_3^{\mu}$.  That leaves $2 \nu - 3 \mu = 2j -2 - \mu$ further
derivatives to act on one of the terms $D^3 R_3$, or  $2j -2 -
\mu$  derivatives to act on the amplitude. The only possible
solutions of $(\nu, \mu)$ are $(j - 1, 0), (j, 1).$ Referring to
statement (i) of Theorem \ref{SETUPA} and to (\ref{AMPLSINGL}), we
see that the principal symbol of the amplitude depends only on
$f_{\pm}, f_{\pm}'$, so there is no way to differentiate the
amplitude $2j - 2$ times to produce the datum $f_{\pm}^{(2j)}(0).$
 Hence,
$(\nu, \mu) = (j, 1)$ and  the only possibility of producing $f_{\pm}^{(2j)}(0)$
is to throw all $2j$ derivatives on the phase.

Now let us determine $I_{\ell_p}^{\max}(\gcal)$ for the labelled
graphs $(\gcal, \ell)$  above. The terms with maximal number $2j$
of derivatives in the  Feynman amplitude (apart from the overall
universal factor in (\ref{PR})) are given for some non-zero
constant $C_{\gcal}$ by
\begin{equation} \begin{array}{lll}  I_{\ell_p}^{\max}(\gcal) &
=& C_{\gcal} (4 r L) {\mathcal A}_0(r)\sum_{p = 1}^{2r}
(h_+^{pp})^j D_{x_p}^{2j}
{\mathcal L}_+ (0) \\ & & \\
& = & C_{\gcal} (4 r L) {\mathcal A}_0(r)  \sum_{p = 1}^{2r}
(h_+^{pp})^j w_+(p) f_{w_+(p)}^{(2j)}(0).
\end{array}
\end{equation}
The factor $(4 r L) {\mathcal A}_0(r)$ comes from the leading
value of the amplitude (cf. Lemma \ref{AMPPROPS}). By  Proposition
\ref{TWOVERSIONS},  $C_{\gcal} = \frac{1}{|Aut(\gcal)|} =
w_{\gcal}$.

Indeed,  to obtain $f_{\pm}^{(2j)}(0)$, all labels at all
endpoints of all edges must be the same index, or otherwise put
 only the `diagonal terms' of ${\mathcal H}_{+}^j$, i.e. those involving only
 derivatives in a single variable $\frac{\partial}{\partial x_k}$, can produce the factor $f^{(2j)}_{\pm}(0)$.
We then use Lemma \ref{AMPPROPS} (va)  to complete the evaluation.
 The part of the $p$th term $ (h_+^{pp})^j
D_{x_p}^{2j} {\mathcal L}_+ (0)$ of the sum which involves
$f^{(2j)}_{w_+(p)}(0)$  equals
$$\begin{array}{l}  (h_+^{pp})^j |
 (f_{w_+(p)}(0) -\  f_{w_+(p+1) }(0))|^{-1}
(f_{w_+(p) }(0)  -\ f_{w_+(p+1) }(0))  f_{w_+(p)}^{(2j)}(0)\\ \\
= (h_+^{pp})^j w_+(p)  f_{w_+(p)}^{(2j)}(0).
\end{array}
$$
by  (\ref{ODDDER}).

We then break up the sums over $p$ of even/odd parity and use
Proposition \ref{HESSPAR} to replace the odd parity Hessian
elements by $h_+^{11}$ and the even ones by $h_+^{22}$. Taking
into account that $w_{+}(p) = 1 (-1)$ if $p$ is odd (even), we
conclude that
\begin{equation}\label{F2JFINAL}  B_{\gamma^r, j-1}  \equiv 8 r L\;  (w_{{\mathcal G}_{1, j}^{2j,
0}})\; {\mathcal A}_0(r) \{  \;(h^{11}_+)^j f^{(2j)}_{+}(0) -
(h^{11}_-)^j f^{(2j)}_{-}(0) \} + \cdots, \end{equation} where
again $\cdots$ refers to terms with $\leq 2j - 1$ derivatives. We
observe that, as claimed, the result is invariant under the
up-down symmetry $f_+ \iff - f_-$ and under the right left
symmetry $f_{\pm}(x) \to f_{\pm}(-x).$

\end{proof}

Thus, we have obtained the even derivative terms in Theorem
\ref{BGAMMAJ}.

\subsubsection{The principal term: The data $f_{\pm}^{(2j - 1)}(0)$}

We now consider the trickier odd-derivative data $f_{\pm}^{(2j -
1)}(0)$ in the stationary expansion of $I_{2r, \rho}^{\sigma_0,
w_{\pm}}$, which will require the attributes of
 the amplitude (\ref{AMPLSINGL}) detailed in Theorem \ref{SETUPA}.

We again claim that the Taylor coefficients $f_{\pm}^{(2j -
1)}(0)$ appear first in the term of order $k^{-j + 1}. $ Further,
only five graphs can produce such a factor, and of these only two
contribute a non-zero Feynman amplitude.  These two  graphs are
illustrated in the figures. In the following section,  we will
show that $f_{\pm}^{(2j - 1)}(0)$ can only occur in  higher order
terms in $k^{-1}$  also in the singular trace terms.

To prove  this, we  first enumerate the labelled graphs $\gcal$ in
the stationary phase expansion of $I_{2r, \rho}^{\sigma_0,
w_{\pm}}$  whose Feynman amplitude $I_{\ell}(\gcal)$  contains a
 factor of $f_{\pm}^{(2j - 1)}(0)$ in the term of order $k^{-j + 1}$,  and we show
that this data does not appear in terms of lower order in
$k^{-1}.$

 We recall
that $\equiv$ means equality modulo $ R_{2r} ({\mathcal J}^{2j -
2} f_+(0), {\mathcal J}^{2j - 2} f_-(0))$.

\begin{lem} \label{LEM5.4} In the stationary phase expansion of $I_{2r,
\rho}^{\sigma_0, w_{\pm}}$,

\noindent{\bf (i)} There are no labelled graphs $\gcal$ with
$-\chi'(\gcal):= - \chi(\gcal') < j - 1$  for which
$I_{\ell}(\gcal)$ contains
the factor $f_{\pm}^{(2j - 1)}(0)$.  \\
\noindent{\bf (ii)} There are exactly two types of labelled
diagrams $(\gcal, \ell)$ with
 $\chi(\gcal') = -j + 1$ such that  $I_{\ell}(\gcal)$ is non-zero and contains the
 factor
  $f_{\pm}^{(2j - 1)}(0)$. They are given by (see figures):

\begin{itemize}

\item ${\mathcal G}_{2, j + 1 }^{2j - 1, 3, 0} \subset {\mathcal
G}_{2, j + 1 }$ with $V = 2, I = j + 1$: Two closed vertices, $j -
1$ loops at one closed  vertex, $1$ loop at the second closed
vertex, one edge between the closed vertices; no open vertex.
Labels $\ell_{p,q}$: All labels at the closed vertex with valency
$2j - 1$ must be the same index $p$ and all at the second closed
vertex must the be same index $q$. Form of  Feynman amplitude:
$(h^{pp}_{+})^{j - 1} h^{qq}_{+} h^{pq}_{+}  D_{x_p}^{2j - 1}
{\mathcal L}_{+} D_{x_q}^3 {\mathcal L}_{+} \equiv (h^{pp}_{+})^{j
- 1} h^{qq}_{+} h^{pq}_{+} f_{\pm}^{(2j - 1)}(0) f_{\pm}^{(3)}(0).
$ Thus, this graph contributes
$$I^{\max}_{{\mathcal G}_{2, j + 1 }^{2j - 1, 3, 0}} = 8 r L {\mathcal A}_r(0) \; (w_{{\mathcal G}_{2,  j + 1
}^{2j - 1, 3,
 0}}) \sum_{p, q = 1}^{2r} (h^{pp}_{+})^{j
- 1} h^{qq}_{+} h^{pq}_{+} w_+(p) w_+(q) f_{w_+(p)}^{(2j - 1)}(0)
f_{w_+(q)}^{(3)}(0).
$$

\item $\hat{{\mathcal G}}_{2, j + 1 }^{2j - 1, 3, 0} \subset
{\mathcal G}_{2, j + 1 } $ with $ V = 2, I = j + 1$: Two closed
vertices, with $j - 2$ loops at one closed vertex, and with three
edges between the two closed vertices; no open vertex. Labels
$\ell_{p,q}$: All labels at the closed vertex with valency $2j -
1$ must be the same index $p$ and all at the second closed vertex
must the be same index $q$; $(h^{pp}_{\pm})^{j - 2}
(h^{pq}_{\pm})^3 D_{x_p}^{2j - 1} {\mathcal L}_{\pm} D_{x_q}^3
{\mathcal L}_{\pm} \equiv (h^{pp}_{\pm})^{j - 2} (h^{pq}_{\pm})^3
f_{\pm}^{(2j - 1)}(0) f_{\pm} ^{(3)}(0). $ Thus, this graph
contributes
$$I^{\max}_{\hat{{\mathcal G}}_{2, j + 1 }^{2j - 1, 3, 0} } = 8
 r L {\mathcal A}_r(0) \; (w_{\hat{{\mathcal G}}_{2,  j
+ 1 }^{2j - 1, 3,
 0}}) \sum_{p, q = 1}^{2r}  (h^{pp}_{\pm})^{j - 2} (h^{pq}_{\pm})^3 w_+(p) w_+(q)  f_{w_+(p)}^{(2j -
1)}(0) f_{w_+(q)} ^{(3)}(0). $$

\item In addition, there are three other graphs whose Feynman
amplitudes contain factors of $ f_{\pm}^{(2j - 1)}(0)$. But for
our special phase and amplitude, the corresponding amplitudes
vanish.

\end{itemize}

\end{lem}

\begin{proof}  It will be seen in the course of the proof that only connected graphs
can contribute highest order derivative data (the amplitude for a
disconnected graph is the product of the amplitudes over its
components).  Connected labelled graphs $(\gcal, \ell)$ with
$-\chi' \leq j - 1$
 for which $I_{\ell}(\gcal)$ contains the factor
$f_{\pm}^{(2j - 1)}(0)$ as a factor must satisfy the following
constraints:
\begin{itemize}

\item (a) $\gcal$ must contain a distinguished vertex (either open
or closed). If it is closed it must have
 valency $\geq 2j - 1.$ If it is open, it must have valency $2j - 2$.
We denote by $\ell$ the number of loops at this vertex and by $e$
the number of non-loop edges at this vertex.

\item (b) $- \chi(\gcal') = I - V  \leq j - 1$

\item (c) Every closed vertex has valency $\geq 3$; hence $2 I \geq 3 V$.

\end{itemize}

We  distinguish   two overall classes of graphs: those for which
the distinguished vertex is open and those for which it is closed.
Statement (a) follows from  the attributes of the amplitude in
Theorem \ref{SETUPA}: In the first case, $2j - 2$ derivatives must
fall on the amplitude (i.e. the open vertex) to produce
$f_{\pm}^{(2j - 1)}(0)$. In the second case, $2j - 1$ derivatives
must fall on the phase (i.e. the closed vertex).

We first  claim that $V \leq 2$ under constraints (a) - (c). When
the distinguished vertex is open, then  $V = 0$ if $-\chi' = j -
1$ (as noted above),  and there are no possible graphs with
$-\chi' \leq j - 2.$  So assume the distinguished vertex is
closed. Let us consider the `distinguished flower' $\Gamma_0$
consisting just of this vertex and of the edges incident on it.
Denoting the number of loops in $\Gamma_0$ by $\ell$, we must have
$2 \ell + e \geq 2j - 1 $ edges in $\Gamma_0$ to produce
$f_{\pm}^{(2j - 1)}(0).$ We then  complete $\Gamma_0$ to a
connected graph $\gcal$ with $-\chi' \leq j - 1$. We may add one
open vertex, $V - 1$ closed vertices and $N$ new edges.

Suppose that there is no open vertex. We then have:
\begin{equation} \left\{ \begin{array}{l} (i) \; 2 \ell + e \geq 2j - 1 \\ \\
(ii) \;  \ell + e - V + N = j - 1 \\ \\
(iii) \; e + 2 N \geq 3 (V - 1) \end{array} \right. \end{equation}
The last inequality follows from the facts that each new vertex
has valency at least three, and that each of the $r$ edges begins
at the distinguished vertex.  Solving for $V$ in (ii)  and
plugging into  (iii) we obtain $N \leq 3j - 3 \ell - 2 e$.
Plugging back into (ii) we obtain
 $V \leq 2j - 2 \ell - e + 1 \leq 2j + 1 - (2j - 1) = 2,$ by (i). Thus the claim is proved.

Now suppose that $\gcal$ contains one open vertex and $V$ closed
vertices.   Then
 (i) and (ii) remain the same since the $\chi(\gcal')$ is computed without counting  the open vertex.
 On the other hand, (iii) becomes $e + 2N \geq 3(V - 1) + 1,$
since the open vertex has valence at least one. This simply subtracts one from
 the previous computation, giving
$V \leq 1.$ Thus, the distinguished vertex is the only closed vertex.

Now we bound $N$ in the connected component of the distinguished
constellation. First suppose that $V = 1$. There is nothing to
bound unless the graph also contains one open vertex, in which
case $N$ counts the number of loops at the open vertex. We claim
that $N = 0$ in this case. Indeed, we have  $\ell + e + N = j$.
Substituting in (i), we obtain $2N + e \leq 1.$ The only solution
is $N = 0, e = 1.$

Next we consider the case $V = 2$. As we have just seen, no open
vertex occurs. From (i) + (ii) we obtain $2 N + e \leq 3,$ hence
the only solutions are $N = 1, e = 1$ or $N = 0, e = 3.$

We tabulate these results as follows:

\center \begin{tabular}{|c|c|c|c|c|} \hline
\multicolumn{5}{|c|}{\bf Graph parameters} \\ \hline
 V &  $\ell $&   e &  N &  O \\
\hline
 0 &  j-1 &  0 &  0 &  1 \\
\hline
 1 &   j  &  0 & 0 &   0 \\
\hline
 1 &  j-1 &  1 &  0  &  1 \\
\hline
 2 &  j-1 &  1 &  1  &  0 \\
\hline
2 & j-2 &  3 & 0  &  0 \\ \hline
\end{tabular}
\bigskip

We now determine the Feynman amplitudes for each of the associated
graphs.   As we will see, the amplitudes vanish for the first
three lines of the table, and do not vanish for the last two. The
non-vanishing diagrams are pictured in the figures (Figures 6 and
7).

\begin{itemize}

\item (i)    The only possible graph
 with $V = 0$ is:  ${\mathcal G}_{0, j - 1}^{0, 2j - 2},  V = 0, I = j - 1$: $j - 1$ loops at the open vertex.
Taking into account the structure of the amplitude in Theorem
\ref{SETUPA}, in  order to produce  $f^{(2j - 1)}(0)$, all labels
at the open vertex must be the same index $p$. We claim that the
Feynman amplitude vanishes:
 \begin{equation} \label{VAN1} I_{{\mathcal G}_{0, j - 1}^{0, 2j - 2}}^{\max} = (Const.) \sum_{p = 1}^{2r} (h^{pp})^{j - 1}  D_{x_p}^{2j - 2}
 {\mathcal A} \equiv 0 \times f_{\pm}^{(2j - 1)}(0) = 0.
 \end{equation}
Indeed, this is the  case $(\mu, \nu) = (j - 1, 0)$ of
(\ref{MSPJTERM}), which  corresponds to applying all derivatives
$D_{x_p}^{2j - 2}$ on the principal symbol $a^0$ of the amplitude
for some $p = 1, \dots, 2r,$ and it is proved in \S \ref{KEYPROPS}
(\ref{USEDINLEM5.4i}) that it vanishes.

\item (ii)  ${\mathcal G}_{1, j }^{2j, 0} \subset {\mathcal G}_{1,
j },  V = 1, I = j$: $j$ loops at the closed vertex. This is the
graph which produced $f^{(2j)}(0)$, and we now verify that it does
not produce an amplitude containing $f^{(2j - 1)}(0)$. To produce
$f^{(2j - 1)}(0)$, all but one label must be the same ($p$), the
last label different ($q \not= p$). Feynman amplitude:
$$I_{{\mathcal G}_{1, j }^{2j, 0} }^{\max} = (Const.) \sum_{p, q = 1}^{2r}(h^{pp})^{j - 1} h^{pq}
  D_{x_p}^{(2j - 1)} D_{x_q}{\mathcal L}  \equiv
(h^{pp})^{j - 1} h^{pq}  f_{\pm}^{(2j - 1)}(0) f_{\pm}'(0) = 0.$$
The vanishing is verified in  \S \ref{KEYPROPS}
(\ref{USEDINLEM5.4ii}).

\item (iii)  ${\mathcal G}_{1, j }^{2j-1, 1} \subset {\mathcal
G}_{1, j },  V = 1, I = j$: $j - 1$ loops at the closed vertex,
one edge between the open and closed vertex. To produce  $f^{(2j -
1)}(0)$,  all labels at the closed vertex must be the same index
$p$. We claim that again the Feynman amplitude vanishes:
$$I_{{\mathcal G}_{1, j }^{2j-1, 1}}^{\max} = (Const.) \sum_{p, q = 1}^{2r} (h^{pp})^{j - 1} h^{pq} D_{x_p}^{2j - 1} {\mathcal
L}  D_q a^0 \equiv 0 \times f^{(2j - 1)}(0) = 0.
$$ Indeed, exactly one derivative is thrown on the amplitude. To
check this, we note that this is the  case $(\mu, \nu) = (j, 1)$
of (\ref{MSPJTERM}) in which
 $ {\mathcal H}^j_{\pm}$ is applied to $a^0_+ R_3.$  To produce the data $f_{\pm}^{(2j - 1)}(0)$,  the operators
 $  D_{x_p}^{2j - 1} D_{x_q}$ contribute by applying $D_{x_p}^{2j - 1}$ to $R_3$ ($p = 1, \dots, 2r)$,
and by applying the final derivative $D_{x_q}$ to the amplitude.
But $\nabla a_+^0(0) = 0$ by (\ref{TAB3}).

\item (iv)   ${\mathcal G}_{2, j + 1 }^{2j - 1, 3, 0} \subset
{\mathcal G}_{2, j + 1 } (- \chi = j - 1; V = 2, I = j + 1$): Two
closed vertices, $j - 1$ loops at one closed  vertex, $1$ loop at
the second closed vertex, one edge between the closed vertices;
the open vertex has valency $0$.
 All labels at the closed vertex with valency $2j - 1$ must be the same
index $p$ and all at the closed vertex must the be same index $q$.
Since there are no derivatives of the amplitude, we extract its
principal term and obtain
$$\begin{array}{lll} I_{{\mathcal G'}_{2, j + 1 }^{2j - 1, 3, 0} }^{\max} &= & 2 r L {\mathcal A}_r(0)
 C_{{\mathcal G}_{2, j + 1}^{2j-1, 3, 0}}   \sum_{p, q = 1}^{2r}
 (h^{pp}_{+})^{j - 1} h^{qq}_{+} h^{pq}_{+}  D_{x_p}^{2j - 1} {\mathcal L}_{+} D_{x_q}^3 {\mathcal
 L}_{+}\\ &&\\
&\equiv &  8 r L {\mathcal A}_r(0)
 C_{{\mathcal G}_{2, j + 1}^{2j-1, 3, 0}}    \sum_{p, q = 1}^{2r} (h^{pp}_{+})^{j - 1} h^{qq}_{+}
h^{pq}_{+} w_+(p) w_+(q) f_{w_+(p)}^{(2j - 1)}(0)
f_{w_+(q)}^{(3)}(0). \end{array}
$$ The calculation of the coefficients is similar to that in
(iii), except that now we have two factors of the phase. The
factor containing $2j-1$  derivatives of $\lcal$ is evaluated in
(iv) - (v) of the table in Lemma \ref{AMPPROPS} and the third
derivative factor is evaluated in \S \ref{FURTHERPROPSAMP} (4).
Again the combinatorial constant is evaluated in Proposition
\ref{TWOVERSIONS}.

\item (v)  There is a second graph  $ \widehat{{\mathcal G}}_{2, j
+ 1 }^{2j - 1, 3, 0} \subset {\mathcal G}_{2, j + 1 } (-\chi= j -
1; V = 2, I = j + 1$): It has two closed vertices, with $j - 2$
loops at one closed vertex, and  three edges between the two
closed vertices; the open vertex has valency $0$. Labels
$\ell_{p,q}$: All labels at the closed vertex with valency $2j -
1$ must be the same index $p$ and all at the closed vertex must
the be same index $q$. Again, there are no derivatives on the
amplitude, and we get
 $$\begin{array}{lll} I_{\widehat{{\mathcal G}}_{2, j + 1 }^{2j - 1, 3, 0}}^{\max} &= &   2 r L {\mathcal A}_r(0)
 C_{\widehat{{\mathcal G}}_{2, j + 1 }^{2j - 1, 3, 0}}   \sum_{p, q = 1}^{2r} (h^{pp}_{+})^{j - 2} (h^{pq}_{+})^3 D_{x_p}^{2j - 1} {\mathcal L}_{+} D_{x_q}^3
 {\mathcal L}_{+} \\ & & \\
 &\equiv&  2 r L {\mathcal A}_r(0)
 C_{\widehat{{\mathcal G}}_{2, j + 1 }^{2j - 1, 3, 0}}   \sum_{p, q = 1}^{2r}
 (h^{pp}_{+})^{j - 2} (h^{pq}_{+})^3 w_+(p) w_+(q) f_{w_+(p)}^{(2j -
 1)}(0)\;
 f_{w_+(q)}^{(3)}(0). \end{array} $$
 As noted above (cf. \S \ref{KEYPROPS} (\ref{THREEDERIVSL})), other (mixed) third derivatives
of ${\mathcal L}$ vanish on the critical set. The combinatorial
constant is evaluated in Proposition \ref{TWOVERSIONS}.

We now combine the terms in (iv) and (v) and evaluate the
coefficients  to obtain
\begin{equation} \label{LPLUS} \begin{array}{l} 2 r L {\mathcal A}_r(0) \; \{(w_{{\mathcal G}_{2,  j + 1
}^{2j - 1, 3,
 0}}) \sum_{q, p = 1 }^{2r}
[(h^{pp}_{+})^{j - 1} h^{qq}_{+} h^{pq}_{+} \\ \\  +
(w_{\widehat{{\mathcal G}}_{2,  j + 1 }^{2j - 1, 3,
 0}})(h^{pp}_{+})^{j - 2}
 (h^{ pq}_{+})^3] w_+(p) w_+(q) f_{w_+(p)}^{(2j - 1)}(0) f_{w_+(q)}^{(3)}(0). \end{array} \end{equation}
 We obtain the  expression stated in Theorem
(\ref{BGAMMAJ}) by breaking up into indices of like parity and
using Proposition \ref{HESSPAR}.

\end{itemize}

\end{proof}

We pause to review the sources of the various  constants and to
check that sums  over the several $\pm$ signs do not cancel. In
particular, it is crucial that the coefficient of
$I^{\max}_{\hat{{\mathcal G}}_{2, j + 1 }^{2j - 1, 3, 0} }$ is
non-zero, since it is this term which determines odd Taylor
coefficients and allows us to decouple even and odd derivative
terms.

\begin{rem} The constants and sums over $\pm$ are of the following
kinds:

\begin{itemize}

\item The factor of $\lcal$ in the amplitude produces $2 r L$.

\item The following $\pm$ signs arise (with some redundancy):
$\gamma^{\pm}$, $f_{\pm}$, $w_{\pm}$ or equivalently
$\lcal_{\pm}$,  $p$ even (odd), and the two terms of $\lcal$ which
depend on a given index $x_p$ (\ref{ODDDER}). Proposition \ref{PM}
shows that the two possible choices of $w_{\pm}$ produce the same
data. Since $\gamma = \gamma^{-1}$ there is no question of
cancellation between $B_{\gamma^{\pm}}$.

\item The odd derivative monomials with maximal derivatives of $f$
have the form
$$f_{+}^{(2j - 1)}(0) f_+^{(3)}(0),\; f_{+}^{(2j - 1)}(0)
f_-^{(3)}(0), \; f_{-}^{(2j - 1)}(0) f_+^{(3)}(0), \; f_{-}^{(2j -
1)}(0) f_-^{(3)}(0). $$ By Theorem \ref{BGAMMAJ}, the  wave
invariants are invariant under  $f_+ \to - f_-, f_- \to - f_+$,
hence the only possible cancellation could occur between
$f_{+}^{(2j - 1)}(0) f_+^{(3)}(0)$ and $ f_{+}^{(2j - 1)}(0)
f_-^{(3)}(0)$. However, no such cancellation occurs,  as noted
after the calculation in (\ref{ODDDER}), or in Theorem
\ref{BGAMMAJ} where it is noted that the monomials always occur in
the form $$w_+(p) w_+(q) f_{w_+(p)}^{(2j - 1)}(0)
 f_{w_+(q)}^{(3)}(0).$$ In fact, the $\pm$ sum in each factor
$D_{x_p}^{2j} \lcal, D_{x_p}^{2j-1} \lcal, D_{x_p}^3 \lcal$ gives
rise to a factors of $4$ in odd derivative terms, and factors of
$2$ in even derivative terms. For the same reason,  no
cancellations occur between the sum over $p$ even versus $p$ odd.

\end{itemize}

\end{rem}

\begin{figure}\label{NONZERO1}
\centerline{\includegraphics{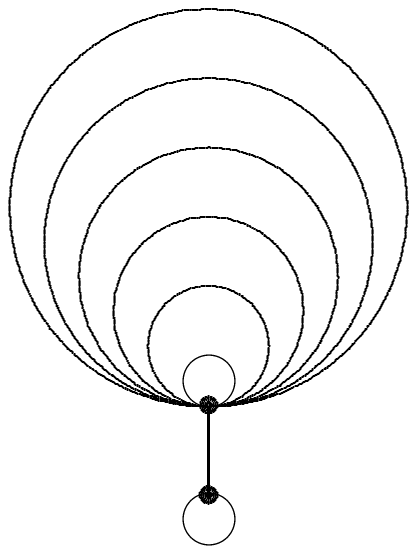}} \caption{$(iv):
{\mathcal G}_{2, j + 1 }^{2j - 1, 3, 0} \subset {\mathcal G}_{2, j
+ 1 } (- \chi = j - 1; V = 2, I = j + 1$)}:
\end{figure}

\begin{figure}\label{NONZERO2}
\centerline{\includegraphics{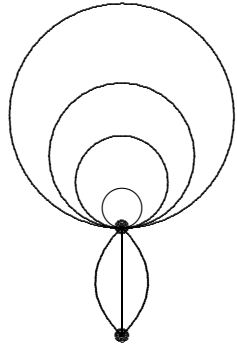}} \caption{$(v):
\widehat{{\mathcal G}}_{2, j + 1 }^{2j - 1, 3, 0} \subset
{\mathcal G}_{2, j + 1 } (-\chi= j - 1; V = 2, I = j + 1$)}:
 \end{figure}

\subsection{\label{NONPRINCIPAL} Non-principal terms}

To complete the proof of Theorems \ref{SETUPA} and \ref{BGAMMAJ},
it suffices to show the non-principal oscillatory integrals $I_{M,
\rho}^{\sigma, w}$ with $M > 2r$  do not contribute the data
$f_{\pm}^{(2j)}(0), f_{\pm}^{(2j - 1)}(0)$ to the coefficient of
the $k^{-j + 1}$- term (or to the $k^{- m}$ term for any $m \leq j
-1$).

 We recall from Proposition \ref{MAINCOR2PROP} that
$I_{M, \rho}^{\sigma, w}$ can only have a critical point if $M
\geq 2r$ and  $M - |\sigma| = 2r$. In the non-principal terms
where $M
> 2r$,  the
 oscillatory integral $I_{M, \rho}^{\sigma, w}$ is obtained by
 regularizing the kernel of $N_{\sigma}$ in Proposition
 \ref{NPOWERFINAL}, which is an oscillatory integral with a
 singular phase and amplitude (cf. \cite{Z5}, \S 6).

 The regularization produces the oscillatory
described in Corollary \ref{SOCINT}. In the case where $M -
|\sigma| = 2r$ it is an integral over ${\bf T}^{2r}$ with the same
phase as in the principal terms but with an amplitude  of order
$-|\sigma|$. The sum over $M$ in Proposition \ref{ASY} and over
$\sigma$ in (\ref{BINOMIAL}) can thus be seen as the construction
of an oscillatory integral expression for the trace of Proposition
\ref{ASY}, with an  amplitude  obtained by regularizing the sum of
singular oscillatory integrals.

The stationary phase analysis of the sub-principal terms $I_{M,
\rho}^{\sigma, w}$ is  therefore almost essentially the same as
for the principal term.  The  only additional feature is the
following description of the amplitude:

\begin{lem} \label{NONPRIN} The amplitude $A_{\sigma}(k, \phi_1, \phi_2)$ of $N_{\sigma}$
in Proposition 3.8 is a semi-classical amplitude of order
$-|\sigma|$. In its semi-classical  expansion $A_{\sigma}(k,
\phi_1, \phi_2) \sim \sum_{n = 0}^{\infty} k^{-|\sigma| - n}
A_{\sigma, n}(\phi_1, \phi_2)$, the term $A_{\sigma, n}$ depends
at most $n + 2$ derivatives of $f$. In particular, the value
$D_{\phi}^{\alpha} A_{\sigma, n}|_{\phi^0}$ of its $\alpha$th
derivative at the critical point depends at most on $n + 2 +
|\alpha|$ derivatives of $f$ at $x = 0$.
\end{lem}

\begin{proof} The
algorithm for calculating $A_{\sigma}(k, \phi_1, \phi_2)$ is given
in \cite{Z5} \S 6 (see also \cite{AG}). We  briefly review the
algorithm in order to prove that the amplitude has the stated
properties.

The algorithm  consists in successively removing factors of $N_0$
from compositions of $N_0$ and $N_1$ in $N_{\sigma}$  (cf. \S 3).
The first step consists in expressing the compositions $N_0 \circ
N_1$ and $N_1 \circ N_0$ as oscillatory integrals of one lower
order (cf. Lemma 6.2 of \cite{Z5}). From the explicit formula for
the composition (cf. (74) of \cite{Z5}), the new  amplitude $A(k +
i \tau, \phi_1, \phi_2)$ has the form \begin{equation}
\label{COMPOSITION} A(k + i \tau, \phi_1, \phi_2) = \int_{\R}
\chi(k, u, \phi_1, \phi_2) G(k + i \tau, u, \phi_1, \phi_2) |u|
H^{(1)}_1((k + i \tau) |u|) e^{i k a u} du, \end{equation} where
$\chi$ is a suitable cutoff and $G$ is a semi-classical amplitude
constructed from the amplitude of $N_1$ (cf. (78)-(79) of
\cite{Z5}). Also, $a = \sin \langle (q(\phi_2) - q(\phi_1),
\nu_{q(\phi_2)}). $

The amplitude    $G$ is constructed as follows: From $N_0$ one
obtains a contribution of $ H^{(1)}_1
 ((k \mu + i \tau) |q(\phi_{3 }) - q(\phi_{1})|)
\cos \angle (q(\phi_{3 }) - q(\phi_{1}), \nu_{q(\phi_{3 })} )$,
while from $N_1$ one obtains a semi-classical amplitude.  One
changes variables by putting
\begin{equation}  u: = \left\{ \begin{array}{ll} |q(\phi_{3})
- q(\phi_{1})| , & \phi_{1} \geq \phi_{3} \\ & \\
- |q(\phi_{3}) - q(\phi_{1})|, & \phi_{1} \leq \phi_{3}
\end{array} \right., \end{equation}
under which the amplitude of $N_1$ is transformed to a smooth
amplitude of the same order in $(\phi_2, u)$, while   the factor
of $\cos \angle (q(\phi_{3 }) - q(\phi_{1}), \nu_{q(\phi_{3 })} )$
changes to $|u| K(\phi_1,  u)$ where $K$ is smooth in $u$. A
simple calculation shows that $K(\phi_1,  0) = -\frac{1}{2}
\kappa(\phi_1)$. The full amplitude $G$ is a product of these two
factors. One sees that it depends analytically  on $f, f', f''$
with $f''$ coming from the cosine factor.

One then  Taylor expands $G$ in $u$ and verifies that it produces
a semi-classical expansion of $A(k + i \tau, \phi_1, \phi_2)$. The
$du$ integrals can be explicitly evaluated using the cosine
transform of the Hankel function (\cite{Z5}, Proposition 4.7; see
also \cite{AG}). The $|u| du$ in the cosine transform  gives rise
to a factor of $k^{-2}$, and the factor of $N_0$ carries a factor
of $k$, so that the removal of $N_0$ introduces a net factor of
$k^{-1}$.  This factor is responsible for the lowering of the
order by one for each removal of $N_0$.

 The coefficient of $k^{-1 - n}$ in the final amplitude thus
derives from the $n$th term in the Taylor expansion of $G(k, u,
\phi)$ in $u$  and in particular depends on the same number of
derivatives of $f$. Since $G$ is an analytic function  of $f, f',
f''$, it follows that the $k^{-1 -n}$ term depends at most on $n +
2$ derivatives of $f$.

The process then repeats as another factor of  $N_0$ is removed
from the resulting composition. The same argument shows that each
elimination of $N_0$ introduces a new factor of $k^{-1}$ which is
unrelated to Taylor expansions of $G$. We now verify that after
$r$ repetitions of the algorithm, the new amplitude is
semi-classical and its $k^{-r - n}$ term depends on only $n + 2$
derivatives of $f$.

We argue by induction, the case $r = 1$ having been checked above.
After $r - 1$ steps, we obtain an oscillatory integral operator
with an amplitude $A_{r-1}$ satisfying the hypothesis and with the
phase of  $N_1$. We then apply the algorithm for the composition
of $N_0$ with this  oscillatory integral operator. It has the form
of (\ref{COMPOSITION}) except that now $G = G_{r}$ is constructed
using $A_{r-1}$ and $N_0$. The algorithm is to multiply $A_{r-1}$
by the cosine factor above, to change variables to $u$, to Taylor
expand the cosine factor to one order to obtain $|u| K$ and to
define $G_r = K A_{r-1} J$ where $J$ is the Jacobian. The Taylor
expansion producing $K$ is responsible for the initial increase in
the number of derivatives of $f$ to $f''$. After that point, it is
only the Taylor expansion of $G_r$ in $u$ which  produces further
derivatives of $f$. Thus, the number of derivatives of $f$ in the
term of order $k^{-r - n}$ is $n + 2$.

It follows that,  after removing all $|\sigma|$ factors of $N_0$,
one obtains an amplitude which is of order $-|\sigma|$ and whose
$k^{- |\sigma| - n}$ term involves at most $n$ derivatives of
$f''$.

\end{proof}

\begin{lem} \label{SUBPRINNOCONT} The non-principal terms do not contribute the data
$f_{\pm}^{2j}(0), f_{\pm}^{2j - 1}(0)$ to the term of order $k^{-1
- j}$. \end{lem}

\begin{proof}

We consider  the diagrammatic analysis of $I_{M, \rho}^{\sigma,
w}$ along the same lines as for  the principal term. The only new
aspect is   the amplitude. Since it now has order $- |\sigma| <
0$, the terms where one differentiates the phase to the maximal
degree now have order $k^{- j + 1 - |\sigma|}$ and thus do not
occur in the $k^{-1 - j}$ term.

The only remaining possibility  is that the data could occur in
terms  where one differentiates the amplitude to the maximal
degree. By Proposition \ref{NONPRIN}, the term of order $k^{-
|\sigma| - n}$ contains at most $n + 2$ derivatives of $f$. To
obtain a term of order $-j + 1$, one needs $ |\sigma| + n \leq j -
1$ and one can take only $2 (j - 1 - |\sigma| - n)$ further
derivatives in the $k^{-j + 1}$ term. This produces a maximum of
$2j - 2|\sigma -n $ derivatives of $f$. The maximum occurs when $n
= 0$, in which case there are $\leq 2j - 2 |\sigma| \leq 2j - 2$
derivatives of $f$.

\end{proof}

For emphasis, we determine the lowest order term in which such
data do occur:

\begin{sublem} In the stationary phase expansion of the non-principal
term $I_{M, \rho}^{\sigma, w}$, the  data $f_{\pm}^{2j}(0),
f_{\pm}^{2j - 1}(0)$ appear first in the $ k^{1 - j - |\sigma|}$-
term.
\end{sublem}

\begin{proof}

To determine the power of $k^{-1}$ in which this data first
appears, we need to minimize $|\sigma| + \nu - \mu$ subject to the
constraint that $2\nu - 3\mu \geq 2j - 3$. This is $|\sigma| $
plus the constrained minimum of $\nu - \mu$. The sole change to
the principal case is that the constraint is $2\nu - 3\mu \geq 2j
- 3$ in the top order term of the amplitude rather than $2\nu -
3\mu \geq 2j - 2.$ Since the solutions must be non-negative
integers, it is easy to check that again $\nu \geq j - 1$ and that
$(\mu, \nu) = (0,j - 1), (1, j) $ achieve the minimum of $\nu -
\mu = j - 1$.
 If there are $r$
drops in the symbol order, we need to minimize $|\sigma| +r + \nu
- \mu$ subject to the constraint that $2\nu - 3\mu \geq 2j - 3 -
r$. The minimizer produces the result stated in the Sublemma.

\end{proof}

This completes the proof of Theorems \ref{SETUPA} and
(\ref{BGAMMAJ}).

\subsection{Appendix: Non-contributing diagrams}

In figures (6)-(7), we displayed the diagrams which contribute
non-zero amplitudes to the leading order derivative terms.  For
the sake of completeness, we also include diagrams which do {\it
not} contribute because the corresponding amplitudes vanish.
Figures (8) - (9)  are labelled consistently with the discussion
above. Figure (5) is also a `non-contributing diagram' to the
leading order odd derivative term.

\begin{figure}\label{one}
\centerline{\includegraphics{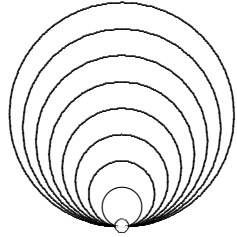}} \caption{$(i):
{\mathcal G}_{0, j - 1}^{0, 2j - 2}$.}
\end{figure}

\begin{figure}\label{two}
\centerline{\includegraphics{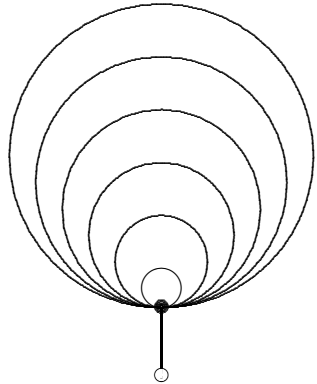}} \caption{$(ii): {\mathcal
G}_{1, j }^{2j-1, 1} \subset {\mathcal G}_{1, j }$).}
\end{figure}


\newpage

\subsection{Balian-Bloch invariants at bouncing ball orbits of
up-down symmetric domains}

We now simplify the expression in Theorem \ref{BGAMMAJ} in the
case of $\Z_2$-symmetric domains. The  following result,  stated
in (\ref{BGAMMAJSYMpre}), is essentially a corollary of Theorem
\ref{BGAMMAJ}. It uses one simplification which will be proved in
Proposition \ref{QBB}.

\begin{cor} \label{BGAMMAJSYM} Suppose that $(\Omega, \gamma)$
 is invariant under an isometric involution $\sigma$, and that
$\gamma$ is a periodic $2$-link reflecting ray which is reversed
by $\sigma.$ Then, modulo the error term  $R_{2r} ({\mathcal
J}^{2j - 2} f(0))$, $ B_{\gamma^r, j - 1} $ is given by the
expression (\ref{BGAMMAJSYMpre}).

\end{cor}

\begin{proof} Using that $f_- = - f_+$, we can cancel the signs in
the formula of Theorem \ref{BGAMMAJ} and add the top and bottom to
obtain,

$$ \begin{array}{lll}B_{\gamma^r, j - 1}  &\equiv &
4 r L\; {\mathcal A}_0(r)\;   \{ (w_{{\mathcal G}_{1, j}^{2j,
0}})\;   \sum_{p = 1}^{2r }  (h^{pp})^j   f^{(2j)}(0)\\ & & \\
&  +& \;  4 \sum_{q, p = 1}^{2 r } [(w_{{\mathcal G}_{2,  j + 1
}^{2j - 1, 3,
 0}}) (h^{pp})^{j - 1} h^{pq} h^{qq} \\&&\\ &+ & 4 (w_{\hat{{\mathcal G}}_{2,  j + 1
}^{2j - 1, 3,
 0}})
(h^{pp})^{j - 2} (h^{pq})^3] \}
   f^{(3)}(0) f^{(2j - 1)}(0).\end{array} $$
Further, in this $\Z_2$-symmetric case, all of the coefficients
$h^{pp}$ are clearly equal. The sum $\sum_{q = 1}^{2r} h^{pq}$ is
 independent of $p$ and is evaluated in Proposition
(\ref{QBB}), leaving
 the stated expression.  \end{proof}

\section{\label{PFONESYM} Proof of Theorem (\ref{ONESYM})}

We now prove the inverse spectral result for simply connected
analytic plane domains with one special symmetry that reverses the
endpoints of a bouncing ball orbit. The method is to recover  the
Taylor coefficients of the boundary defining function from the
Balian-Bloch invariants at this orbit.

As simple warm-up for the proof, we give a new proof that
centrally symmetric convex  analytic domains whose shortest orbit
is the unique orbit of its length (up to time-reversal) are
spectrally determined within that class:
\medskip

\noindent{\bf Proof of Corollary \ref{ONESYMCOR}}: Consider the
wave invariants of the shortest orbit as  given in Theorem
\ref{BGAMMAJ}. They are spectral invariants since the shortest
length is a spectral invariant. By Ghomi's theorem \cite{Gh}, the
shortest orbit is a bouncing ball orbit.  The orbit must be
invariant under the two symmetries up to time-reversal since its
length is of multiplicity one.  Hence, the two symmetries imply
that $f_+ = - f_- := f$ and that $f^{(2j + 1)}(0) = 0$ for all
$j$.  It follows that $f^{(2j)}(0)$ are spectral invariants for
each $j$, and thus the domain is determined. QED

The same proof shows that simply connected analytic  domains with
the symmetry of an ellipse and with one axis  of prescribed length
$L$ are spectrally determined in that class.

\subsection{ Completion of the proof of Theorem \ref{ONESYM}}
We now complete the proof of Theorem \ref{ONESYM}. Thus, we
assume that $(\Omega, \gamma)$ is up-down symmetric, i.e. is
invariant under an isometric involution $\sigma$, and that
$\gamma$ is a periodic $2$-link reflecting ray which is reversed
by $\sigma.$

There are two overall steps in the proof. First, and foremost, we
study the expressions in Corollary \ref{BGAMMAJSYM}. The key point
is that the Hessian of the length function is a circulant matrix
in the symmetric case, and that allows us to analyze  the Hessian
sums which occur as coefficients in the Balian-Bloch wave
invariants. In particular, we decouple even and odd derivatives
using the behavior of the Hessian sums under iterates $\gamma^r$.
After that, a simple inductive argument shows that all Taylor
coefficients of $f_+$ may be determined from the Balian-Bloch
invariants.

We now begin the analysis of the Hessian sums.

\subsection{Circulant Hessian at $\Z_2$-symmetric bouncing ball orbits}

In the case of $\Z_2$-symmetric domains in the sense of Theorem
\ref{ONESYM},  $R_A = R_B: = R$ and
\begin{equation} \label{FLOQUET} \cos \alpha/2 = 2 (1 - \frac{L}{R}) \;\;(\mbox{elliptic
case}),\;\;\;\cosh \alpha/2 = 2 (1 - \frac{L}{R})
\;\;(\mbox{hyperbolic case}). \end{equation} We put:
\begin{equation} \label{adef} a = - 2 \cos \alpha/2\;\; (\mbox{elliptic
case}),\;\;\;a = - 2 \cosh \alpha/2  \;\;(\mbox{hyperbolic case}).
\end{equation} By  \ref{COSALPHA} and Proposition \ref{KT},
 the Hessian of the Length function  in Cartesian graph
coordinates simplifies to:
\begin{equation} \label{HBB2} H_{2r} = \frac{-1}{L} \left\{ \begin{array}{lllll} a & 1 & 0 & \dots & 1  \\ & & & & \\
  1 & a & 1 & \dots & 0 \\ & & & & \\ 0 & 1 &  a  & 1 & 0  \\ & & & & \\ 0 & 0 & 1 & a & 1 \dots
 \\ & & & & \\\dots & \dots & \dots & \dots & \dots \\ & & & & \\
1 & 0 & 0 & \dots &  a \end{array} \right\}.
\end{equation}
  We observe that
(\ref{HBB2})  is a symmetric {\it circulant } matrix (or simply
circulant) of the form
\begin{equation} (-L ) H_{2r} = C(a, 1, 0, \dots, 0, 1), \end{equation}
where   a circulant is a matrix of the form (cf. \cite{D})
\begin{equation} \label{CIRC}  C(c_1, c_2, \dots, c_n) = \left\{ \begin{array}{llll} c_1 & c_2&   \dots & c_n
 \\ & & &  \\
  c_n& c_1 & \dots & c_{n-1}
 \\ & & &  \\\dots & \dots &  \dots & \dots \\ & & &  \\
c_2& c_3&  \dots &  c_1 \end{array} \right\}. \end{equation}
Circulants are diagonalized by the finite Fourier matrix $F$ of
rank $n$ defined by
\begin{equation}\label{CIRCDIAG}  F^* = n^{-1/2} \left\{ \begin{array}{llll} 1 & 1&   \dots & 1
 \\ & & &  \\
  1& w & \dots & w^{n-1}
 \\ & & &  \\\dots & \dots &  \dots & \dots \\ & & &  \\
1& w^{n-1}&  \dots &  w^{(n-1)(n-1)} \end{array} \right\},\;\;\;\;
w = e^{\frac{2 \pi i }{n}} \end{equation} Here, $F^* = (\bar{F})^T
= \bar{F}$ is the adjoint of $F$. By \cite{D}, Theorem 3.2.2, we
have $C = F^* \Lambda F$ where
\begin{equation} \Lambda = \Lambda_C = \mbox{diag}\; (p_C(1), \dots, p_C(w^{n-1}), \;\; \mbox{with}\;\;
p_C(z) = c_1 + c_2 z + \dots + c_n z^{n-1}. \end{equation} Here,
by $\mbox{diag}$ we mean the diagonal matrix with the exhibited
entries.

\subsection{Diagonalizing $H_{2r}^{-1}$}

Applying the above to  $C = H_{2r}$:

\begin{prop} \label{INVHESS} We have: $$ \begin{array}{l}
 H_{2r}^{-1}
= -  L\; F^* (\mbox{diag}\;(\frac{1}{a + 2} , \dots,\frac{1}{ a +
2 \cos \frac{(2r - 1) \pi}{r})} ) F,
\end{array}$$
where $a$ is defined in (\ref{adef}).
\end{prop}

\begin{proof}

We use the notation  $p_{a, r} (z)$ for $p_C(z)$ in the case where
$C = C(a, 1, 0, \dots, 0, 1)$. Thus,
\begin{equation} \label{PALPHA}
p_C(z) : = p_{a, r} (z) : =  a + z + z^{2r-1}.
\end{equation}
By (\ref{CIRCDIAG})  we have,
\begin{equation} H_{2r} = \frac{-1}{L} \;  F^*\mbox{diag}\; (p_{a,  r}(1), \dots,
p_{a,  r}(w^{2r-1}))F, \;\;\;\;\; (w = e^{\frac{i \pi}{r}}).
\end{equation}
Since
\begin{equation} \label{PALPHAR}
p_{a, r}(w^k) : = a + w^k + w^{-k}, \;\;(w = e^{\frac{i \pi}{r}})
\end{equation}
we have
\begin{equation} \label{CIRCHESS} H_{2r} = \frac{-1}{L}\;  F^* \mbox{diag}\;(a + 2,
 \dots, a + 2 \cos \frac{(2r - 1) \pi}{r}) F,
\end{equation}
and inverting gives the statement.
\end{proof}

\subsection{Matrix elements of $H_{2r}^{-1}$  at a $\Z_2$-symmetric bouncing ball orbit}

We will need explicit formulae for the matrix elements
$h^{pq}_{2r}$ of $H_{2r}^{-1}$.  The diagonalization of
$H_{2r}^{-1}$ above gives one kind of formula. We also consider a
second approach to inverting $H_{2r}$ (due to \cite{K}) via finite
difference equations.  The two approaches give quite different
formulae for the inverse Hessian sums and have different
applications in the inverse results. In several of the
calculations in this section, we assume for simplicity of
exposition that $\gamma$ is elliptic;  the hyperbolic case is
easier and  all formulae analytically continue from the elliptic
to the hyperbolic cases.

For our purposes it will suffice to know the formulae for the
elements $h^{1 q}_{2r}.$ To emphasize the fact that the matrix
elements depend on, and only on,  $(r, a)$ we denote them by
$h^{pq}_{2r}(a)$.  The first formula comes directly from the
diagonalization above.

\begin{prop}  \label{FOURIER} With the above notation, we have
$$ h^{1 q}_{2r} (a) = \; \frac{- L  }{2r}  \sum_{k = 0}^{2r - 1} \frac{w^{ (q - 1) k }
}{p_{a, r}(w^k)}, \;\;\;\; (w = e^{\frac{i \pi}{r}}) $$ where the
denominators are defined in (\ref{PALPHA})-(\ref{PALPHAR}).

\end{prop}

The second, finite difference,  approach expresses the inverse
Hessian matrix elements $h_{2r}^{pq}$ in terms of Chebychev
polynomials $T_n,$ resp. $U_n$,  of the first, resp. second, kind.
They are defined by:
$$T_n(\cos \theta) = \cos n \theta,\;\;\;\;\; U_n(\cos \theta) = \frac{\sin (n + 1) \theta}{\sin \theta}.$$

\begin{prop} \label{CHEB} \cite{K} (p. 190)
 With the above notation,

$$\begin{array}{l} (-L)^{-1} h^{pq}_{2r}(a) = \frac{1}{2[ 1 - T_{2r}(- a/2)]} [U_{2r - q + p -1}(-a/2) +
 U_{q - p - 1}(-a/2)],\;\;\;\; 1 \leq p \leq q \leq 2r)  \end{array}$$ \end{prop}
We note that $h^{pq} = h^{qp}$ so this formula determines all of
the matrix elements.

The special cases $r = 1, 2$ are already very helpful in the
inverse problem. We recall that
$$\begin{array}{l} T_1(x) = x, \;  T_2(x) = 2 x^2 - 1, \; T_3(x) = 4 x^3 - 3x, \; T_4(x) = 8 x^4 - 8x^2 + 1; \;\;\\
\\U_1(x) = 2x, \;
U_2(x) = 4x^2 - 1, \; U_3(x) = 8 x^3 - 4 x, \; U_4 = 16 x^2 - 12
x^2 + 1,\end{array}$$ from which we calculate:

\begin{equation} \label{H2INV} H_2^{-1} = \frac{- L}{a^2 - 4} \begin{pmatrix}  a & - 2 \\ & \\
- 2 &  a \end{pmatrix}, \end{equation}

and

\begin{equation} \label{H4INV} H_4^{-1} = \frac{-L}{a^4 - 4 a^2} \begin{pmatrix} a^3 - 2 a &
-
a^2 & 2a & -a^2 \\ \\
-a^2 & a^3 - 2a & - a^2 & 2a \\ \\
2a & -a^2 & a^3 - 2a & -a^2 \\ \\
-a^2 & 2a & -a^2 & a^3 - 2a \end{pmatrix}. \end{equation}

In terms of Floquet angles, we have (in the  elliptic case),
\begin{equation}  h_{2r}^{pq} = \frac{-L }{2[ 1 - T_{2r}(\cos \alpha/2)]} [U_{2r - q + p -1}(-\cos \alpha/2) + U_{q - p - 1}(-\cos \alpha/2)],\;\;\; (1 \leq p \leq q \leq 2r)), \end{equation} hence
\begin{equation}\label{HPQ} \begin{array}{l} (- L)^{-1}\; h_{2r}^{pq}
= \left\{ \begin{array}{ll}  \frac{(-1)^{p - q}}{2[ 1 - \cos r \alpha)]} [ \frac{\sin (2r - q + p) \alpha/2}{\sin \alpha/2}  +  \frac{\sin ( q - p ) \alpha/2}{\sin \alpha/2} ] & (1 \leq p \leq q \leq 2r) \\ & \\
\frac{(-1)^{p - q}}{2[ 1 - \cos r \alpha]} [ \frac{\sin (2r - p +
q) \alpha/2}{\sin \alpha/2}  +  \frac{\sin ( p - q )
\alpha/2}{\sin \alpha/2} ] & (1 \leq q \leq p \leq 2r) \end{array}
\right. \end{array} \end{equation}

We note that  the expression in Proposition (\ref{FOURIER}) is the
Fourier inversion formula for (\ref{HPQ}).

\begin{cor} \label{SQUARES} We have:
$(- L)^{-1}\;  h_{2r}^{11}  = \frac{U_{2r - 1}(-\frac{a}{2})}{2 (1
- T_{2r}(- \frac{a}{2}))} = \frac{ \sin r \alpha}{2 (1 - \cos r
\alpha) \sin \frac{\alpha}{2}}  =  \frac{1}{2 \sin
\frac{\alpha}{2}} \cot \frac{r \alpha}{2}.$
\end{cor}

\subsection{Linear sums}

We now complete the proof of Corollary \ref{BGAMMAJSYM} by summing
the  matrix elements in the first row $[H_{2r}^{-1}]_1 = (h^{11},
\dots, h^{1 (2r)})$ (or column) of the inverse. As a check on the
notation and assumptions, we calculate it in two different ways:

\begin{prop} \label{QBB} Suppose that $\gamma$ is a $\Z_2$-symmetric bouncing ball orbit. Then, for any $p$,
  $\sum_{q = 1}^{2 r }  h_{2r}^{pq} = \frac{-L}{a + 2} = \frac{- L }{ 2 -  2\cos \alpha/2} $.   \end{prop}

\begin{proof} Because $H_{2r}^{-1}$ is a circulant matrix, the
column sum is the same for all columns. Hence we only need to
consider the first column.

\noindent{\bf (i)} By Proposition \ref{FOURIER}, we have
$$\begin{array}{l}  \sum_{q = 1}^{2r} h^{p q}_{2r} = \sum_{q = 1}^{2r} h^{1 q}_{2r} = \frac{-L }{2r} \sum_{q = 1}^{2r}
 \sum_{k = 0}^{2r - 1} \frac{w^{ (q - 1) k } }{p_{a, r}(w^k)} \\ \\
=  (-L) \; \sum_{k = 0}^{2r - 1} \frac{\delta_{k0} }{p_{a,
r}(w^k)} = \frac{-L }{p_{a, r}(1)} = \frac{-L}{2 + a} =
\frac{-L}{2 - 2\cos \alpha/2}.
\end{array}$$

\noindent{\bf (ii)}  Since  $\sum_{q = 1}^{2r} h^{1 q}_{2r} =
\sum_{q = 1}^{2r} h^{p q}_{2r}$ for any $p = 1, \dots, 2r$, we can
set $p = 1$ in the sum over $q$ to obtain,
\begin{equation}  1 = \sum_{p, q = 1}^{2r} h_{p q'} h^{p q} = [\sum_{p = 1}^{2r} h_{p q'}][ \sum_{q = 1}^{2 r }  h^{pq}]. \end{equation}
It then follows from (\ref{COSALPHA}) and Proposition (\ref{KT})
that  $(-L)^{-1}\; \sum_{p = 1}^{2r} h_{p q'} = 2 + a =  2 - 2
\cos \alpha/2$.

\end{proof}

\subsection{Decoupling Balian-Bloch invariants}

Corollary (\ref{BGAMMAJSYM}) expresses  $B_{\gamma^r, j-1 } $ in
terms of inverse Hessian matrix elements. To prove Theorem
\ref{ONESYM}, it is essential  to show that we can separately
determine the two terms
 \begin{enumerate}

\item $
 (h^{11}_{2r}(a) )^2 \{ 2(w_{{\mathcal G}_{1, j}^{2j,
0}}) f^{(2j)}(0) + 4 \frac{(w_{{\mathcal G}_{2,  j + 1 }^{2j - 1,
3,
 0}})}{2 +a }
 f^{(3)}(0) f^{(2j - 1)}(0)\},$

 \item  $4 (w_{\hat{{\mathcal G}}_{2,  j + 1
}^{2j - 1, 3,
 0}}) \sum_{q = 1}^{2r} (h^{1 q}_{2r}(a))^3\} f^{(3)}(0) f^{(2j - 1)}(0).$ \end{enumerate}

To decouple the terms we prove that they have behave independently
under iterates $r$ of the bouncing ball orbit.  We use the simple
observation:
\begin{lem} \label{DECOUPLE}  Let $F_3(r, a) =  \sum_{q = 1}^{2r} (h^{1 q}_{2r}(a))^3$.
If $(h^{11}_{2r}(a))^{-2} F_3(r, a)$  is non-constant in $r = 1,
2, 3, \dots$,  then both terms (1)-(2)  can be determined from
their sum as $r$ ranges over ${\bf N}.$ \end{lem}

\begin{proof}  Put $$A = \{2(w_{{\mathcal G}_{1, j}^{2j,
0}})f^{(2j)}(0) + 4\frac{(w_{{\mathcal G}_{2,  j + 1 }^{2j -1, 3,
 0}})}{2 + a} f^{(3)}(0) f^{(2j -
1)}(0)\},\;\; B = 4 (w_{\hat{{\mathcal G}}_{2,  j + 1 }^{2j - 1,
3,
 0}}) f^{(3)}(0) f^{(2j - 1)}(0).$$ It is assumed
that we know $ (h^{11}_{2r}(a))^2 A + F_3(r, a) B$ for all $r \in
{\bf N}$. To determine $A, B$ it is clearly sufficient that the
matrix
$$\left( \begin{array}{ll} (h^{11}_{2r}(a))^2 & F_3(r, a) \\ & \\ (h^{11}_{2s}(a))^2 & F_3(s, a) \end{array} \right)$$
is invertible for some integers $r \not= s$. But this says
precisely that
 $(h^{11}_{2r}(a))^{-2} F_{3}(r, a) \not= (h^{11}_{2s}(a))^{-2} F_{3}(s, a)$ for some integers $r \not= s.$
\end{proof}

\subsection{Cubic Hessian sums}

We now prove that $(h^{11}_{2r}(a))^{-2}F_3(r, a)$ is indeed
non-constant for all but finitely many  $a$.

\begin{prop} \label{FINITE} The `bad' set ${\mathcal B}$ of
(\ref{BAD}) consists of $\{0, - 1, \pm 2\}.$  \end{prop}

\begin{proof}

We will give two   different proofs of the finiteness of
${\mathcal B}$. In both, we consider the sets $${\mathcal B}_{r,
s} =
 \{ a \in \R: (h^{11}_{2r}(a))^{-2} F_3 (r, a) = (h^{11}_{2s}(a))^{-2} F_3(s, a) \}.$$

\subsubsection{First proof of Proposition \ref{FINITE}: Dedekind
sums}

The first is based on an explicit calculation of $F_3(r, a)$ as a
Dedekind sum. It is not very efficient in bounding the cardinality
of ${\mathcal B}_{r,s}$ but gives a clear proof that this set is
finite.

\begin{lem} \label{NQ}  We have :  $$ F_3(r, a)=
\frac{(-L)^3 }{(2r)^2}  \sum_{k_1, k_2 = 0}^{2r - 1} \frac{1}{(a+
 2 \cos \frac{k_1 \pi}{r}) (a + 2 \cos \frac{k_2 \pi}{r})(a + 2 \cos \frac{(k_1 +k_2) \pi}{r})} .$$
 In the hyperbolic case, we obtain a similar result with $\cos$
 replaced by $\cosh.$
\end{lem}

\begin{proof} Using Proposition \ref{FOURIER},  we have (with $w = e^{\frac{\pi i}{r}}$, and $\equiv$ equal to congruence
modulo $2r$),

\begin{equation}\begin{array}{l} \label{CUBES} (2r)^3 \frac{-1}{L^3} \sum_{q = 1}^{2r} (h^{1 q}_{2r}(a))^3
=  \sum_{q = 1}^{2r} \{\sum_{k = 0}^{2r - 1} \frac{w^{ (q - 1) k } }{p_{a, r}(w^k)}\}^3  \\ \\
= \sum_{q = 1}^{2r} \{\sum_{k_1, k_2, k_3  = 0}^{2r - 1} \frac{w^{ (q - 1) (k_1 + k_2 + k_3) } }{p_{a, r}(w^{k_1})
 p_{a, r}(w^{k_2}) p_{a, r}(w^{k_3})}\} \\ \\
=  2r  \sum_{0 \leq k_i \leq 2r-1; k_1 + k_2 +  k_3 \equiv 0}
\frac{1 }{p_{a, r}(w^{k_1}) p_{a, r}(w^{k_2}) p_{a, r}(w^{k_3})}
\\ \\
 = 2r  \sum_{0\leq k_i \leq 2r-1; k_1 + k_2 +  k_3 \equiv 0}\frac{1 }{(a + 2 \cos \frac{k_1 \pi}{r})
 (a + 2\cos \frac{k_2 \pi}{r})(a + 2 \cos \frac{k_3 \pi}{r})}
\\ \\
= 2r  \sum_{k_1, k_2 = 0}^{2r - 1} \frac{1 }{(a + 2 \cos \frac{k_1
\pi}{r}) (a + 2 \cos \frac{k_2 \pi}{r})(a + 2 \cos \frac{(k_1
+k_2) \pi}{r})}.
\end{array}
 \end{equation}

\end{proof}

We now complete the proof of Proposition \ref{FINITE}. By
Corollary \ref{SQUARES},  $(h_{2r}^{11}(a))^{-2} F_3(r, a)$ is the
rational function $\left(\frac{U_{r-1}(-\frac{a}{2}) }{2 (1 -
T_r(-\frac{a}{2}))} \right)^{-2} F_3(r, a)$, where as above, $T_n,
U_n$ are the Chebychev polynomials.

We now observe that for $r \not= s$,
$\left(\frac{U_{r-1}(-\frac{a}{2}) }{2 (1 - T_r(-\frac{a}{2}))}
\right)^{-2} F_3(r, a)$ and $\left(\frac{U_{s-1}(-\frac{a}{2}) }{2
(1 - T_s(-\frac{a}{2}))} \right)^{-2} F_3(s, a)$ are independent
rational functions. Indeed, the poles for given $r$ are the values
$a = - 2 \cos \frac{\alpha}{2}$ where $\alpha = \frac{2 \pi k}{r}$
for some $k = 1, \dots, 2r.$ Hence, there can exist only finitely
many solutions of the equation
\begin{equation} \label{RATFUNEQ}  \left(\frac{U_{r-1}(-\frac{a}{2}) }{2 (1 - T_r(-\frac{a}{2}))}
\right)^{-2} F_3(r, a)= \left(\frac{U_{s-1}(-\frac{a}{2}) }{2 (1 -
T_s(-\frac{a}{2}))} \right)^{-2} F_3(s, a) \end{equation} for any
$r\not= s$, i.e.  ${\mathcal B}_{r, s}$ is finite.

\end{proof}

It is interesting to observe that the sums  above are generalized
Dedekind sum, i.e. the sum $ \sum_{\zeta \in D_r} I_3(\zeta, z) $
of the function
 $$I_3(x; z) = \frac{1 }{(z +
\cos  x_1) (z+ \cos  x_2)(z + \cos (x_1 + x_2))}$$ over the set
$D_{2r}$ of $2r$th roots of unity $\frac{\pi k}{r} \mbox{mod}\; 2
\pi \Z^2$ with $k = (k_1, k_2) \in [0, 2r - 1] \times [0, 2r - 1]$
of the torus.  The summand is
 is  a continuous periodic
function of $(x_1, x_2) \in [0, 1] \times [0, 1]$   for $z
\notin[-1, 1]$. In fact,  $I_3(x, z)$ is also symmetric under
inversion and reflection across the diagonal and the  sum has
additionally the form of a multiple Dedekind sum
$$s_2(1, 1;
2r) =  \sum_{k_1,k_2 (\text{mod}\,2r)}f(k_1, r) f(k_2, r)f(k_1 +
k_2, r)
$$ of
two variables in the sense of L. Carlitz \cite{Ca}, with $f (k, r)
= \frac{1 }{(z + \cos 2\pi k/r)}.$

We remark that under the non-degeneracy  assumption that
$\alpha/\pi \notin \Q$,
 $\cos \alpha/2$ is never a pole of $F_3(r, z)$ for any $r$.
  In the hyperbolic case,
  it is obvious that $\cosh
\alpha$ is never a pole of $F_3(r, z)$.

\subsubsection{Second proof: Explicit inversion of the Hessian}

We now give a second (and quite  elementary) method of determining
${\mathcal B}$  by simply using the formulae for $H_2^{-1}$
(\ref{H2INV}) and $H_4^{-1}$ (\ref{H4INV}). This calculation is
due to the referee and to H. Hezari.

From the explicit formula for $H_2^{-1}$ we have:
$$\sum_{q = 1}^2 (h_2^{1q}(a))^3 = \left( \frac{- L}{a^2 - 4} \right)^3
(a^3 - 8).$$ Further, $h_{2}^{11} = \frac{-a L }{a^2 - 4}$. From
the explicit formula for $H_4^{-1}$ we have
$$\sum_{q = 1}^4 (h_4^{1q}(a))^3 = \left( \frac{- L}{a^4 - 4a^2} \right)^3
(a^9 - 6a^7 -2 a^6 + 12 a^5).$$ Further, $h_4^{11} = (-L)
\frac{a^3 - 2a}{a^4 - 4 a^2}. $

Thus, ${\mathcal B}_{1,2}$ is the set of solutions $a$ of the
equation
$$\begin{array}{l} \frac{a^3 - 8}{(a^2 - 4)^3} \frac{(a^2 -
4)^2}{a^2} = \frac{(a^4 - 4a^2)^2}{(a^3 - 2a)^2} \frac{a^9 - 6 a^7
- 2a^6 + 12 a^5}{(a^4 - 4a^2)^3} \\ \\
\iff (a^3 - 2a)^2 (a^3 - 8) = a^9 - 6 a^7 - 2a^6 + 12 a^5
\end{array}$$
A little bit of cancellation reduces the  equation to degree $6$.
The distinct roots are  $\{0, -1, 2, -2\}$.   QED

\subsection{Final step in proof of Theorem \ref{ONESYM}: Inductive determination of Taylor coefficients}

  We now prove
by induction that
 on $j$ that $f^{2j}(0), f^{(2j - 1)}(0)$ are wave trace
 invariants, hence spectral invariants of the Laplacian among
 domains in ${\mathcal D}_{1, L}$.

It is clear for $j = 1$ since $(1 - L f^{(2)}(0) = \cos \alpha/2$
(resp. $\cosh \alpha/2$) and $\alpha$ is a Balian-Bloch (wave
trace) invariant at $\gamma$ (see \cite{Fr}). In the case $j = 2$,
 the Balian-Bloch invariants have the form (\ref{BGAMMAJSYMpre}).
 Using that $\alpha$ is a Balian-Bloch invariant and
 the decoupling argument of Lemma \ref{DECOUPLE} and
 Proposition \ref{FINITE}, $(f^{(3)}(0))^2$ is a spectral
 invariant.
 By reflecting the domain across
the bouncing ball axis if necessary, we may assume with no loss of
generality that $f^{(3)}(0)
> 0$, and we have then determined $(f^{(3)}(0))$ from the
sequence of Balian-Bloch invariants. Using again that $\alpha$ is
determined by the  Balian-Bloch invariants, it  follows that $
f^{(4)}(0)$ is determined.

 We now carry forward the argument by induction. As $j \to j + 1$, we may assume that ${\mathcal J}^{2j -2} f(0)$ is
 known. The terms denoted $R_{2r} {\mathcal J}^{2j -2} f(0)$ in
 Theorem \ref{BGAMMAJ} are universal polynomials in the data ${\mathcal J}^{2j -2}
 f(0)$, hence are also known. Thus, it suffices
to determine $f^{(2j)}(0), f^{(2j - 1)}(0)$ from
(\ref{BGAMMAJSYMpre}). By the decoupling argument, we can
determine $(f^{(3)}(0))(f^{(2j - 1)}(0))$, hence $(f^{(2j -
1)}(0))$, as long as $(f^{(3)}(0)) \not= 0.$ But then we can
determine $f^{(2j)}(0).$ By induction, $f$ is determined and hence
the domain.

This completes the proof of Theorem (\ref{ONESYM}). QED

\begin{rem} From this argument it is only necessary  that the coefficients
$w_{{\gcal}}$ etc. are  non-zero and universal. It is not
necessary to know the precise values  of the coefficients of
$f^{(2j)}(0), f^{(2j - 1)}(0)$.

\end{rem}

\subsection{The case where $f^{(3)}(0) = 0$}

If  $f^{(3)}(0) = 0$,  the inductive  argument clearly  breaks
down. There is a natural analogue of it as long as $f^{(5)}(0)
\not= 0$. We only sketch the analogue to make it seem plausible,
but do not provide a complete proof.

Instead of inductively determining
 $f^{(2j)}(0), f^{(2j - 1)}(0)$, we inductively determine  $f^{(2j)}(0), f^{(2j - 3)}(0)$ by
a similar argument. Since $f^{(3)}(0) = 0$, the terms $ f^{(2j -
1)}(0)$ have zero coefficients, and each  new `odd' term as $j \to
j + 1$ now  has the form $[\sum_{q = 1}^r (h^{pq})^5] f^{(5)}(0)
f^{(2j - 3)}(0)$. To carry out the analogue of the previous
argument, it suffices to show that $h_{2r}^{-1} [\sum_{q = 1}^r
(h^{pq})^5]$ is a non-constant function of $r$. It should  be
plausible that this is the case, at least if we exclude a finite
number of values of the Floquet exponents.

 There then arises an infinite sequence of further
sub-cases where all odd derivatives vanish up to some $j_0 + 1.$
To handle this case, we would need to show that
$h_{2r}^{-1}[\sum_{q = 1}^r (h^{pq})^{2j_0 + 1}]$ is non-constant
for all $j_0$. This should again be plausible.

In the case where all odd derivatives vanish, the function $f_+$
is even and the proof reduces to the previously established case
of two symmetries.

\section{\label{PFDISYM} Proof of Theorem (\ref{DISYM})}

We now generalize the results from a bouncing ball orbit to
iterates of a primitive  $D_m$-invariant $m$-link reflecting ray
$\gamma$. For short, we call $\gamma$ a $D_m$-ray.

\subsection{Structure of  coefficients at a $D_m$-ray}

\subsubsection{$D_m$-rays}

In the dihedral case, we orient $\Omega$ so that the center of the
dihedral action is $(0,0)$ and so that one vertex $v_0$ of
$\gamma$ lies on the $y$-axis. We again define a  small strip
$T_{\epsilon}(\gamma)$, which intersects the boundary in $n$ arcs.
We label the one through $v_0$ by $\alpha$. We then write $\alpha$
as the graph $y = f(x)$ of a function defined on a small interval
around $(0,0)$ on the horizontal axis. Since we are only
considering $D_n$-invariant rays, the domain is entirely
determined by $\alpha$ and $f$.

We first need to choose a convenient parametrization of $\partial
\Omega \cap T_{\epsilon}(\gamma).$ Either a polar parametrization
or a Cartesian parametrization would do. For ease of comparison to
the bouncing ball case, we prefer the Cartesian one. Thus, we use
the parametrization $x \in (-\epsilon, \epsilon) \to (x, f(x))$
for the $\alpha$ piece. We then use $x \to R_{2 \pi/m}^j (x,
f(x))$ for the rotate $R_{2 \pi/m}^j \alpha.$  When considering
$\gamma^r$, we need variables $x_{j s} (j = 1, \dots, m; s = 1,
\dots, r),  x_{js} \to  R_{2 \pi/m}^j (x_{js} , x_{js}).$ We have:
 $$\begin{array}{l} R_{2 \pi/m}^{\sigma(p)}( x_p, f (x_p)) = (x_p^{\sigma(p)}, (f(x_p))^{\sigma(p)})\\ \\:=
 (\cos (2 p \pi /m) x_p + \sin(2 p  \pi/m) f(x_p), -\sin  (2  p \pi /m) x_p + \cos (2 p \pi/m) f(x_p)).\end{array}  $$
We also put $(-1, f'(x_p))^{\sigma(p)}:= R_{2 \pi/m}^{\sigma(p)}
(-1, f'(x_p)).$

We then define the length functional

\begin{equation} \begin{array}{ll} {\mathcal L}^{\sigma} (y, x_0, x_1, \dots, x_{m r}) & = |(x_0, y) -
 (x_1, f(x_1))^{\sigma(1)}|+ |(x_0, y) - (x_{rm}, f(x_{rm}))^{\sigma(rm)}|  \\ & \\ &
+ \sum_{p = 1}^{mr-1}|(x_p, f(x_p))^{\sigma(p)} - (x_{p + 1},
f(x_{p + 1}))^{\sigma(p + 1)}|  \end{array} \end{equation} We will
need a formula for its Hessian in the case of a $D_m$-ray. By
(\cite{KT}, Proposition 3), the Hessian $H_{rm}$ in  $x-y$
coordinates at the  critical point $(x_1, \dots, x_{r m })$
corresponding to $\gamma^r$ is given by the matrix (\ref{HBB})
with
  $s = \frac{2L}{R \sin \vartheta}$.

A key point in what follows (as in \cite{Z1, Z2}) is that the
reflection symmetry of $\alpha$ and $f$ implies that $ f^{(2j -
1)}(0) = 0$ for all $j$.  This eliminates the most serious
obstacle to recovering $f$ from the wave trace invariants at
$\gamma^r$, namely the fact that in the transition from the $j$th
Balian-Bloch invariant to the $(j + 1)$st, two new derivatives of
$f$ appear.

As in the $\Z_2$-symmetric case, there are principal and
non-principal terms. The principal term in the $D_m$ case,
analogously to the bouncing ball case,  equals $Tr \rho * N_1^{m
r} \circ N_1'(k) \circ \chi(k)$ for $r$ repetitions of the
dihedrally symmetric orbit.

 In analogy to Lemma (\ref{BGAMMAJ}) we
prove:

\begin{lem} \label{BGAMMAJDI} Let $\gamma$ be a $D_m$- ray, and let $\rho$ be a smooth cutoff to $t = r L_{\gamma}$
as above. Then:

\begin{itemize}

\item $B_{\gamma^r, j} = p_{m, r, j}(f^{(2)}(0), f^{(3)}(0),
\cdots, f^{(2j + 2)}(0))$  where $p_{2, r, j}(\xi_1, \dots,
\xi_{2j})$ is a polynomial. It is homogeneous of degree $-j$ under
the dilation $f \to \lambda f,$ is invariant under the
substitution $f(x) \to f(-x)$, and has degree $j + 1$ in the
Floquet data $e^{i \alpha r}.$

\item   In the expansion in Theorem (1.1) of \cite{Z5} of
  $Tr R_{ \rho} ((k + i \tau))$, $f^{(2j)}(0)$ appears first in the $k^{-j + 1}$st order term,
  and then only in the $k^{-j + 1}$st order term in the stationary
  phase expansion of the principal term $Tr \rho * N_1^{m r}
\circ N_1'(k) \circ \chi(k)$;

\item This coefficient has the form
$$ \begin{array}{l}  B_{\gamma^r, j - 1} = m r (h^{11})^j f^{(2j)} (0)
+ R_{m r} ({\mathcal J}^{2j - 2} f(0)),\end{array} $$ where the
remainder $ R_{m r} ({\mathcal J}^{2j - 2} f(0))$ is a polynomial
in the designated jet of $f.$

\end{itemize}

\end{lem}

\noindent{\bf Proof of Lemma (\ref{BGAMMAJDI})}

We use the analogue of Theorem \ref{MAINCOR2} for the case of the
dihedral ray. As in the case of a bouncing ball orbit, we have a
finite number of oscillatory integrals $I_{M \rho}^{\sigma, w}$
arising from the regularization of the trace. We express the
resulting oscillatory integrals in Cartesian coordinates of (polar
coordinates are also convenient for this calculation). We put $x =
(x_0, y_0)$. Each oscillatory integral $I_{M, \rho}^{\sigma, w}$
localizes  at critical points, we may insert a cutoff to
$T_{\epsilon}(\gamma).$ This gives $m^M$ possible terms,
corresponding to the possible choices of the arcs in the product
$(\partial \Omega \cap T_{\epsilon}(\gamma))^M.$ We put:
$$\{m^M\} := \{ \sigma: \Z_M \to \{1, \dots, m\} \},$$
and write
 $$\begin{array}{l} R_{2 \pi/m}^{\sigma(p)}( x_p, f (x_p)) = (x_p^{\sigma(p)}, (f(x_p))^{\sigma(p)})\\ \\:=
 (\cos (2 p \pi /m) x_p + \sin(2 p  \pi/m) f(x_p), -\sin  (2  p \pi /m) x_p + \cos (2 p \pi/m) f(x_p)).\end{array}  $$
We also put $(-1, f'(x_p))^{\sigma(p)}:= R_{2 \pi/m}^{\sigma(p)} (-1, f'(x_p)).$

The oscillatory integrals have the  phase functions ${\mathcal
L}^{\sigma}$
 on $ (\partial \Omega \cap T_{\epsilon}(\gamma))^{r m }$  of the form:
\begin{equation}  {\mathcal L}^{\sigma} ( x_1, \dots, x_{m r}) =
\sum_{p = 1}^{mr-1}|(x_p, f(x_p))^{\sigma(p)} - (x_{p + 1}, f(x_{p
+ 1}))^{\sigma(p + 1)}|  \end{equation}

Only $2m$ $\sigma$'s ( $2$ modulo cyclic permutations)  give length functions which have critical points with critical value $r L_{\gamma},$
 namely the  ones $\sigma_0$ where $\sigma_0 (n) = R(\pm n 2 \pi/m)$. Indeed, the only Snell polygon with this length is
$\gamma^r$ by assumption, and so $(x_1^{\sigma(1)}, \dots, x_{rm}^{\sigma(rm)})$ must correspond to the vertices
of $\gamma^{\pm r}$. Since the good length functions represent isometric situations, it suffices to consider
the case $\sigma_0 (n) = R(n 2 \pi/m)$. In this case, we denote the length function simply by $L$ and to simplify
the notation we drop the subscript in $\sigma_0.$

We now make a stationary phase analysis as in the bouncing ball
case to obtain the expressions in  Theorem  (\ref{BGAMMAJ}). As
mentioned above,  there are two principal terms: The  principal
oscillatory integrals $I_{r m, \rho}^{\sigma_0, w_{\pm}}$ are
those in  which $M = r m$ and in which no factors of $N_0$ occur,
i.e. $\sigma_0(j) = 1$ for all $j = 1, \dots, r m$. Also, there
are now $m$ components of the boundary at the reflection points,
and $w_{\pm}$ cycles around them for $r$ iterates.

\subsection{The principal terms}

They have the phase
 \begin{equation}  {\mathcal L}^{\sigma} ( x_1, \dots, x_{mr}) =
 \sum_{j = 1}^{mr-1} \sqrt{(x_{j + 1}^{\sigma(j + 1)} - x_j^{\sigma(j)})^2 + (f(x_{j + 1})^{\sigma(j + 1)} -
(f (x_j))^{\sigma(j)})^2},  \end{equation} and the amplitude

\begin{equation}\label{AMPLSINGLDI} \begin{array}{l} a^0(k, x_1, \dots, x_{mr}, y)
=
 \Pi_{p = 1}^m a_1((k + i \tau)  \sqrt{(x_{p - 1}^{\sigma(p-1)}  - (x_{p}^{\sigma(p)} )^2 + (f(x_{p - 1})^{\sigma(p-1)}  - f(x_{p})^{\sigma(p)} )^2})\\ \\ \frac{(x_{p - 1}^{\sigma(p-1)}, , f(x_{p - 1})^{\sigma(p-1)}  )   - ( x_{p}^{\sigma(p)},   f(x_{p})^{\sigma(p)}  )\cdot\nu_{x_{p}^{\sigma(p)},   f(x_{p})^{\sigma(p)}  }}{
\sqrt{(x_{p - 1}^{\sigma(p-1)}  - x_{p}^{\sigma(p)} )^2 + (f(x_{p - 1})^{\sigma(p-1)}  - f(x_{p})^{\sigma(p)} )^2}}\end{array} \end{equation}
We observe that it has the form ${\mathcal A}(x, y, f, f').$  The $f'$ dependence will
be particularly important later on.

\subsubsection{The principal term: The data $f^{2j}(0)$}

As in the bouncing ball case, by the same argument, the data
$f^{(2j)}(0)$ appears first in the term of order $k^{-j + 1}$ and
it appears linearly in the term $ a^0  {\mathcal H}^{j } R_3. $
 We now show that its coefficient is given by the formula in Lemma (\ref{BGAMMAJDI}). Due to symmetry,
it suffices to consider any axis and one endpoint of it.   We observe that only the `diagonal terms' of ${\mathcal H}^j$, i.e. those involving only
 derivatives in a single variable $\frac{\partial}{\partial x_k}$, can produce the factor $f^{(2j)}(0)$. Since
 $f'(0) = x|_{x = 0} = 0$ and since the
angle between successive links and the normal equals $\pi/m$ an
examination of (\ref{AMPLSINGL}) shows that
 the coefficient of $f^{(2j)}(0)$ equals
$$\begin{array}{l} \sum_{p = 1}^{rm }  (h^{pp})^j  (\frac{\partial}{\partial x_p})^{2j}  {\mathcal L}^{\sigma} (y; x_0, \dots, x_k, \dots x_{mr})
=   [\sum_{p = 1}^{m r}  (h^{pp})^j]  f^{(2j)}(0). \end{array} $$
 The data $f^{(2j - 1)}(0)$ vanishes due to the symmetry around each dihedral axis.

 Finally,
as in the bouncing ball case, and for the same reasons,
non-principal oscillatory integrals do not contribute to this
data.

This completes the proof of Lemma (\ref{BGAMMAJDI}). \qed \medskip

\noindent{\bf Remark} It would  also be  natural to employ polar
coordinates in the proof. In that case, we  align $\Omega$ so that
one of the reflection axes is the positive $x$-axis, and   express
$\partial \Omega$ parametrically in the form $r = r(\vartheta)$
where $\vartheta$ is the angle to the $x$-axis. Then $r(-
\vartheta) = r(\vartheta), r(\vartheta + \frac{ 2 \pi j}{m}) =
r(\vartheta).$ The goal then is to determine $r$. To do so, we
write out that $q (\vartheta) = (r(\vartheta) \cos (\vartheta),
r(\vartheta) \sin (\vartheta)$ and compute as above. We find that
$r^{(2j)}(0)$ arises first in the $k^{-1 + j}$ term with the same
coefficient as for $f^{(2j)}(0)$ above. The rest of the proof
proceeds as with Cartesian coordinates.

\subsection{Dihedral domains: Proof of Theorem (\ref{DISYM})}

We now complete the proof of  Theorem (\ref{DISYM}).

We prove   by induction on $j$ that $f^{2j}(0)$ is a  Balian-Bloch
invariant. It is clear for $j = 1$ since $(1 - L f^{(2)}(0) =
\cos(h) \alpha/2$ and $\alpha$ is a Balian-Bloch  (wave trace)
invariant at $\gamma$. In general, the eigenvalues of $P_{\gamma}$
are wave trace invariants \cite{F}.

Assuming the result for $n < j - 1$, it follows that
$p^{sub}_{r , n-1}$ is a spectral invariant. It thus suffices to extract
$f^{2j}(0)$ from $p^0_{r, j-1}$, i.e. from  $  \{  \sum_{p = 1}^{2r }  (h^{pp})^j \}  f^{(2j)}(0). $

Thus, the only missing step is to show that if
 $\gamma$ is $D_m$-ray,  then the $h^{pp}$ are  Balian-Bloch invariants of $\gamma^r$.
In other words,  that $s$ is a wave trace invariant. If $\lambda,
\lambda^{-1}$ denote the eigenvalues of $P_{\gamma^r}$, then we
have $\lambda + \lambda^{-1} = 2 + \det H_{mr}.$ Here we use that
all $b_j$ equal $1$. It follows that $s$ is a function of
$\lambda$, hence that it is a Balian-Bloch invariant.

The proof of Theorem (\ref{DISYM}) is complete.

\end{document}